12mm

# There is No Standard Model of ZFC and $ZFC_2$ with Henkin semantics

**Jaykov Foukzon[1], Elena Men'kova[2]**

[1]Department of Mathematics, Israel Institute of Technology, Haifa, Israel
Email: jaykovfoukzon@list.ru

[2]Lomonosov Moscow State University, Moscow, Russia
Email: e_menkova@mail.ru

## Abstract

**In this article we proved so-called strong reflection principles corresponding to formal theories** *Th* **which has omega-models or nonstandard model with standard part. A possible generalization of Löb's theorem is considered. Main results are**:
(i) $\neg Con(ZFC + \exists M_{st}^{ZFC})$, (ii) $\neg Con(ZF + (V = L))$, (iii) $\neg Con(NF + \exists M_{st}^{NF})$,
(iv) $\neg Con(ZFC_2)$,
(v) **let** $k$ **be inaccessible cardinal then** $\neg Con(ZFC + \exists\kappa)$.

## 6mm . 1.Introduction 3mm

### 5mm . 1.1. Main results 2mm

Let us remind that accordingly to naive set theory, any definable collection is a set. Let $R$ be the set of all sets that are not members of themselves. If $R$ qualifies as a member of itself, it would contradict its own definition as a set containing all sets that are not members of themselves. On the other hand, if such a set is not a member of itself, it would qualify as a member of itself by the same definition. This contradiction is Russell's paradox. In 1908, two ways of avoiding the paradox were proposed, Russell's type theory and Zermelo set theory, the first constructed axiomatic set theory. Zermelo's axioms went well beyond Frege's axioms of extensionality and unlimited set abstraction, and evolved into the now-canonical Zermelo–Fraenkel set theory *ZFC*. *"But how do we know that ZFC is a consistent theory, free of contradictions? The short answer is that we don't; it is a matter of faith (or of skepticism)"*— E. Nelson wrote in his paper [1]cite: Nelson11. However, it is deemed unlikely that even $ZFC_2$ which is significantly stronger than *ZFC* harbors an unsuspected contradiction; it is widely believed that if $ZFC$ and $ZFC_2$ were consistent, that fact would have been uncovered by now. This much is certain —$ZFC$ and $ZFC_2$ is immune to the classic paradoxes of naive set theory: Russell's paradox, the Burali-Forti paradox, and Cantor's paradox.

**Remark 1.1.1**. The inconsistency of the second order set theory $ZFC_2$ originally have been uncovered in [2]cite: Foukzon13 and officially announced in [3]cite: Lemhoff16, see also ref. [4], [5], [6]cite: Foulzon17cite: Foukzon15cite: Foukzon19.

**Remark 1.1.2**. In order to derive a contradiction in second order set theory $ZFC_2$ with the Henkin semantics [7]cite: Henkin50, we remind the definition given in P. Cohen handbook [8]cite: Cohen66, (see [8]]cite: Cohen66 Ch. III, sec. 1, p. 87). P. Cohen wroted: "A set which can be obtained as the result of a transfinite sequence of predicative definitions Godel called "constructible". His result then is that the constructible sets are a model for $ZF$ and that in this model GCH and AC hold. The notion of a predicative construction must be made more precise, of course, but there is essentially only one way to proceed. Another way to explain constructibility is to remark that the constructible sets are those sets which jnust occur in any model in which one admits all ordinals. The definition we now give is the one used in [9]cite: Godel68.

**Definition 1.1.1**. [8]]cite: Cohen66. Let $X$ be a set. The set $X'$ is defined as the union of $X$ and the set $Y$ of all sets $y$ for which there is a formula $A(z,t_1,\ldots,t_k)$ in $ZF$ such that if $A_X$ denotes $A$ with all bound variables restricted to $X,$ then for some $\bar{t}_i, i=1,\ldots,k.$ in $X,$

$$y=\left\{z\in X\,|\,A_X(z,\bar{t}_1,\ldots,\bar{t}_k)\right\}. \tag{1}$$

Observe $X'\subsetneqq P(x)\cup X,$ $\overline{\overline{X}}'=\overline{\overline{X}}$ if $X$ is infinite (and we assume $AC$). It should be clear to the reader that the definition of $X',$ as we have given it, can be done entirely within $ZF$ and that $Y=X'$ is a single formula $A(X,Y)$ in $ZF.$ In general, one's intuition is that all normal definitions can be expressed in $ZF,$ except possibly those which involve discussing the truth or falsity of an infinite sequence of statements. Since this is a very important point we shall give a rigorous proof in a later section that the construction of $X'$ is expressible in $ZF.$"

**Remark 1.1.3**. We will say that a set $y$ is definable by the formula $A(z,t_1,\ldots,t_k)$ relative to a given set $X.$

**Remark 1.1.4**. Note that a simple generalsation of the notion of the definability which has been by Definition 1.1.1 immediately gives Russell's paradox in second order set theory $ZFC_2$ with the Henkin semantics [7]cite: Henkin50.

**Definition 1.1.2**. [6]**cite: Foukzon19**. (i)We will say that a set $y$ is definable relative to a given set $X$ iff there is a formula $A(z,t_1,\ldots,t_k)$ in $ZFC$ then for some $\bar{t}_i\in X, i=1,\ldots,k.$ in $X$ there exists a set $z$ such that the condition $A(z,\bar{t}_1,\ldots,\bar{t}_k)$ is satisfied and $y=z$ or symbolically

$$\exists z[A(z,\bar{t}_1,\ldots,\bar{t}_k)\wedge y=z]. \tag{2}$$

It should be clear to the reader that the definition of $X',$ as we have given it, can be done entirely within second order set theory $ZFC_2$ with the Henkin semantics [7]cite: Henkin50 denoted by $ZFC_2^{Hs}$ and that $Y=X'$ is a single formula $A(X,Y)$ in $ZFC_2^{Hs}$.

(ii)We will denote the set $Y$ of all sets $y$ definable relative to a given set $X$ by $Y\triangleq\Im_2^{Hs}$.

**Definition 1.1.3**. Let $\Re_2^{Hs}$ be a set of the all sets definable relative to a given set $X$ by the first order 1-place open wff's and such that

$$\forall x(x\in\Im_2^{Hs})[x\in\Re_2^{Hs}\Leftrightarrow x\notin x]. \tag{3}$$

**Remark 1.1.5**. (a) Note that $\Re_2^{Hs}\in\Im_2^{Hs}$ since $\Re_2^{Hs}$ is a set definable by the first order 1-place open wff $\Psi(Z,\Im_2^{Hs})$ :

$$\Psi(Z,\Im_2^{Hs})\triangleq\forall x(x\in\Im_2^{Hs})[x\in Z\Leftrightarrow x\notin x], \tag{4}$$

**Theorem 1.1.1**. [6]cite: Foukzon19. Set theory $ZFC_2^{Hs}$ is inconsistent.

Proof. From (3) and Remark 1.1.2 one obtains

$$\Re_2^{Hs}\in\Re_2^{Hs}\Leftrightarrow\Re_2^{Hs}\notin\Re_2^{Hs}. \tag{5}$$

From (5) one obtains a contradiction

$$(\Re_2^{Hs} \in \Re_2^{Hs}) \wedge (\Re_2^{Hs} \notin \Re_2^{Hs}). \tag{6}$$

**Remark 1.1.6**. Note that in paper [6]**cite: Foukzon19** we dealing by using following definability condition: a set $y$ is definable if there is a formula $A(z)$ in $ZFC$ such that

$$\exists z[A(z) \wedge y = z]. \tag{7}$$

Obviously in this case a set $Y = \Re_2^{Hs}$ is a countable set.

**Definition 1.1.4**. Let $\Re_2^{Hs}$ be the countable set of the all sets definable by the first order 1-place open wff's and such that

$$\forall x(x \in \Im_2^{Hs})[x \in \Re_2^{Hs} \Leftrightarrow x \notin x]. \tag{8}$$

**Remark 1.1.7**. (a) Note that $\Re_2^{Hs} \in \Im_2^{Hs}$ since $\Re_2^{Hs}$ is a set definable by the first order 1-place open wff $\Psi(Z, \Im_2^{Hs})$ :

$$\Psi(Z, \Im_2^{Hs}) \triangleq \forall x(x \in \Im_2^{Hs})[x \in Z \Leftrightarrow x \notin x], \tag{9}$$

one obtains a contradiction $(\Re_2^{Hs} \in \Re_2^{Hs}) \wedge (\Re_2^{Hs} \notin \Re_2^{Hs})$.

In this paper we dealing by using following definability condition.

**Definition 1.1.5**. (i) Let $M_{st} = M_{st}^{ZFC}$ be a standard model of $ZFC$. We will say that a set $y$ is definable relative to a given standard model $M_{st}$ of $ZFC$ if there is a formula $A(z, t_1, \ldots, t_k)$ in $ZFC$ such that if $A_{M_{st}}$ denotes $A$ with all bound variables restricted to $M_{st}$, then for some $\bar{t}_i \in M_{st}, i = 1, \ldots, k.$ in $M_{st}$ there exists a set $z$ such that the condition $A_{M_{st}}(z, \bar{t}_1, \ldots, \bar{t}_k)$ is satisfied and $y = z$ or symbolically

$$\exists z[A_{M_{st}}(z, \bar{t}_1, \ldots, \bar{t}_k) \wedge y = z]. \tag{10}$$

It should be clear to the reader that the definition of $M'_{st}$, as we have given it, can be done entirely within second order set theory $ZFC_2$ with the Henkin semantics.

(ii) In this paper we assume for simplicity but without loss of generality that

$$A_{M_{st}}(z, \bar{t}_1, \ldots, \bar{t}_k) = A_{M_{st}}(z). \tag{11}$$

**Remark 1.1.8**. Note that in this paper we view (i) the first order set theory $ZFC$ under the canonical first order semantics (ii) the second order set theory $ZFC_2$ under the Henkin semantics [7]cite: Henkin50 and (iii) the second order set theory $ZFC_2$ under the full second-order semantics [8]cite: Cohen66, [9]cite: Godel68 , [10]cite: Rossberg , [11]cite: Shapiro00 , [12]cite: RayoUzquiano99 but also with a proof theory based on formal Urlogic [13]cite: Vaananen01 .

**Remark 1.1.9**. Second-order logic essantially differs from the usual first-order predicate calculus in that it has variables and quantifiers not only for individuals but also for subsets of the universe and variables for $n$-ary relations as well
[7]-[13]cite: Henkin50cite: Cohen66cite: Godel68cite: Rossbergcite: Shapiro00cite: RayoUzquiano99
deductive calculus $\mathbf{DED}_2$ of second order logic is based on rules and axioms which guarantee that the quantifiers range at least over definable subsets [7]cite: Henkin50. As to the semantics, there are two types of models: (i) Suppose $\mathbf{U}$ is an ordinary first-order structure and $\mathbf{S}$ is a set of subsets of the domain $A$ of $\mathbf{U}$. The main idea is that the set-variables range over $\mathbf{S}$, i.e.
$\langle \mathbf{U}, \mathbf{S} \rangle \models \exists X\Phi(X) \Leftrightarrow \exists S(S \in \mathbf{S})[\langle \mathbf{U}, \mathbf{S} \rangle \models \Phi(S)]$.

We call $\langle \mathbf{U}, \mathbf{S} \rangle$ a Henkin model, if $\langle \mathbf{U}, \mathbf{S} \rangle$ satisfies the axioms of $\mathbf{DED}_2$ and truth in $\langle \mathbf{U}, \mathbf{S} \rangle$ is preserved by the rules of $\mathbf{DED}_2$. We call this semantics of second-order logic the Henkin semantics and second-order logic with the Henkin semantics the Henkin second-order logic. There is a special class of Henkin models, namely those $\langle \mathbf{U}, \mathbf{S} \rangle$ where $\mathbf{S}$ is the set of all subsets of $A$.

We call these full models. We call this semantics of second-order logic the full semantics and second-order logic with the full semantics the full second-order logic.

**Remark 1.1.10**. We emphasize that the following facts are the main features of second-order logic:

**1. The Completeness Theorem**: A sentence is provable in $\mathbf{DED}_2$ if and only if it holds in all Henkin models [7] -

[13]cite: Henkin50cite: Cohen66cite: Godel68cite: Rossbergcite: Shapiro00cite: RayoUzquiano99

**2**. **The Löwenheim-Skolem Theorem**: A sentence with an infinite Henkin model has a countable Henkin model.

**3**. **The Compactness Theorem**: A set of sentences, every finite subset of which has a Henkin model, has itself a Henkin model.

**4**. **The Incompleteness Theorem**: Neither $\mathbf{DED}_2$ nor any other effectively given deductive calculus is complete for full models, that is, there are always sentences which are true in all full models but which are unprovable.

**5**. Failure of the Compactness Theorem for full models.

**6**. Failure of the Löwenheim-Skolem Theorem for full models.

**7**. There is a finite second-order axiom system $\mathbb{Z}_2$ such that the semiring $\mathbb{N}$ of natural numbers is the only full model of $\mathbb{Z}_2$ up to isomorphism.

**8**. There is a finite second-order axiom system $RCF_2$ such that the field $\mathbb{R}$ of the real numbers is the only full model of $RCF_2$ up to isomorphism.

**Remark 1.1.11**. For let second-order $ZFC$ be, as usual, the theory that results obtained from $ZFC$ when the axiom schema of replacement is replaced by its second-order universal closure, i.e.

$$\forall X[Func(X) \Rightarrow \forall u \exists v \forall r[r \in v \Leftrightarrow \exists s(s \in u \wedge (s,r) \in X)]], \qquad (12)$$

where $X$ is a second-order variable, and where $Func(X)$ abbreviates " $X$ is a functional relation", see [12]cite: RayoUzquiano99.

Thus we interpret the wff's of $ZFC_2$ language with the full second-order semantics as required in [12]cite: RayoUzquiano99, [13]cite: Vaananen01 but also with a proof theory based on formal urlogic [13]cite: Vaananen01.

**Designation 1.1.1**. We will denote: (i) by $ZFC_2^{Hs}$ set theory $ZFC_2$ with the Henkin semantics,

(ii) by $ZFC_2^{fss}$ set theory $ZFC_2$ with the full second-order semantics,

(iii) by $\overline{ZFC}_2^{Hs}$ set theory $ZFC_2^{Hs} + \exists M_{st}^{ZFC_2^{Hs}}$ and

(iv) by $ZFC_{st}$ set theory $ZFC + \exists M_{st}^{ZFC}$, where $M_{st}^{Th}$ is a standard model of the theory $Th$.

**Remark 1.1.12**. There is no completeness theorem for second-order logic with the full second-order semantics. Nor do the axioms of $ZFC_2^{fss}$ imply a reflection principle which ensures that if a sentence $Z$ of second-order set theory is true, then it is true in some model $M^{ZFC_2^{fss}}$ of $ZFC_2^{fss}$ [11]cite: Shapiro00.

Let $Z$ be the conjunction of all the axioms of $ZFC_2^{fss}$. We assume now that: $Z$ is true, i.e. $Con\left(ZFC_2^{fss}\right)$. It is known that the existence of a model for $Z$ requires the existence of strongly inaccessible cardinals, i.e. under $ZFC$ it can be shown that $\kappa$ is a strongly inaccessible if and only if $(H_\kappa, \in)$ is a model of $ZFC_2^{fss}$. Thus

$$\neg Con(ZFC_2^{fss}) \Rightarrow \neg Con(ZFC + \exists\kappa)). \qquad (13)$$

In this paper we prove that:

(i) $ZFC_{st} \triangleq ZFC + \exists M_{st}^{ZFC}$ (ii) $\overline{ZFC}_2^{Hs} \triangleq ZFC_2^{Hs} + \exists M_{st}^{ZFC_2^{Hs}}$ and (iii) $ZFC_2^{fss}$ is inconsistent, where $M_{st}^{Th}$ is a standard model of the theory $Th$.

**Axiom** $\exists M^{ZFC}$. [8]cite: Cohen66. There is a set $M^{ZFC}$ and a binary relation $\epsilon \subseteq M^{ZFC} \times M^{ZFC}$ which makes $M^{ZFC}$ a model for $ZFC$.

**Remark 1.1.13**. (i) We emphasize that it is well known that axiom $\exists M^{ZFC}$ a single statement in $ZFC$ see [8]cite: Cohen66, Ch. II, section 7. We denote this statement throught all this paper by symbol $Con(ZFC; M^{ZFC})$.The completness theorem says that $\exists M^{ZFC} \Leftrightarrow Con(ZFC)$.

(ii) Obviously there exists a single statement in $ZFC_2^{Hs}$ such that $\exists M^{ZFC_2^{Hs}} \Leftrightarrow Con(ZFC_2^{Hs})$. We denote this statement through all this paper by symbol $Con\left(ZFC_2^{Hs}; M^{ZFC_2^{Hs}}\right)$ and there exists a single statement $\exists M^{Z_2^{Hs}}$ in $Z_2^{Hs}$. We denote this statement throught all this paper by symbol $Con\left(Z_2^{Hs}; M^{Z_2^{Hs}}\right)$.

**Axiom** $\exists M_{st}^{ZFC}$. [8]cite: Cohen66. There is a set $M_{st}^{ZFC}$ such that if $R$ is $\{\langle x,y\rangle | x \in y \wedge x \in M_{st}^{ZFC} \wedge y \in M_{st}^{ZFC}\}$then $M_{st}^{ZFC}$ is a model for $ZFC$ under the relation $R$.

**Definition 1.1.6**. [8]cite: Cohen66. The model $M_{st}^{ZFC}$ is called a standard model since the relation $\in$ used is merely the standard $\in$- relation.

**Remark 1.1.14**. Note that axiom $\exists M^{ZFC}$ doesn't imply axiom $\exists M_{st}^{ZFC}$, see ref. [8]cite: Cohen66.

**Remark 1.1.15**. We remind that in Henkin semantics, each sort of second-order variable has a particular domain of its own to range over, which may be a proper subset of all sets or functions of that sort. Leon Henkin (1950) defined these semantics and proved that Gödel's completeness theorem and compactness theorem, which hold for first-order logic, carry over to second-order logic with Henkin semantics. This is because Henkin semantics are almost identical to many-sorted first-order semantics, where additional sorts of variables are added to simulate the new variables of second-order logic. Second-order logic with Henkin semantics is not more expressive than first-order logic. Henkin semantics are commonly used in the study of second-order arithmetic. Väänänen [13]cite: Vaananen01 argued that the choice between Henkin models and full models for second-order logic is analogous to the choice between *ZFC* and **V** (**V** is von Neumann universe), as a basis for set theory: "As with second-order logic, we cannot really choose whether we axiomatize mathematics using **V** or *ZFC*. The result is the same in both cases, as *ZFC* is the best attempt so far to use **V** as an axiomatization of mathematics."

**Remark 1.1.16**.Note that in order to deduce: (i) $\sim Con(ZFC_2^{Hs})$ from $Con(ZFC_2^{Hs})$,
(ii) $\sim Con(ZFC)$ from $Con(ZFC)$, by using Gödel encoding, one needs something more than the consistency of $ZFC_2^{Hs}$, e.g., that $ZFC_2^{Hs}$ has an omega-model $M_{\omega}^{ZFC_2^{Hs}}$ or an standard model $M_{st}^{ZFC_2^{Hs}}$ i.e., a model in which the *integers are the standard integers and the all wff* of $ZFC_2^{Hs}$, *ZFC, etc. represented by standard objects.*To put it another way, why should we believe a statement just because there's a $ZFC_2^{Hs}$-proof of it? It's clear that if $ZFC_2^{Hs}$ is inconsistent, then we won't believe $ZFC_2^{Hs}$-proofs. What's slightly more subtle is that the mere consistency of $ZFC_2$ isn't quite enough to get us to believe arithmetical theorems of $ZFC_2^{Hs}$; we must also believe that these arithmetical theorems are asserting something about the standard naturals. It is "conceivable" that $ZFC_2^{Hs}$ might be consistent but that the only nonstandard models $M_{Nst}^{ZFC_2^{Hs}}$ it has are those in which the integers are nonstandard, in which case we might not "believe" an arithmetical statement such as "$ZFC_2^{Hs}$ is inconsistent" even if there is a $ZFC_2^{Hs}$-proof of it.

**Remark 1.1.17**. Note that assumption $\exists M_{st}^{ZFC_2^{Hs}}$ is not necessary if nonstandard model $M_{Nst}^{ZFC_2^{Hs}}$ is a transitive or has a standard part $M_{st}^{Z_2^{Hs}} \subset M_{Nst}^{Z_2^{Hs}}$, see [14]cite: Friedman73 , [15]cite: MagidorShelahStavi83.

**Remark 1.1.18**. Remind that if $M$ is a transitive model, then $\omega^M$ is the standard $\omega$. This implies that the natural numbers, integers, and rational numbers of the model are also the same as their standard counterparts. Each real number in a transitive model is a standard real number, although not all standard reals need be included in a particular transitive model. Note that in any nonstandard model $M_{Nst}^{Z_2^{Hs}}$ of the second-order arithmetic $Z_2^{Hs}$ the terms $\bar{0}$, $\mathbf{S}\bar{0} = \bar{1}, \mathbf{SS}\bar{0} = \bar{2}, \ldots$ comprise the initial segment isomorphic to $M_{st}^{Z_2^{Hs}} \subset M_{Nst}^{Z_2^{Hs}}$. This initial segment is called the standard cut of the $M_{Nst}^{Z_2^{Hs}}$. The order type of any nonstandard model of $M_{Nst}^{Z_2^{Hs}}$ is equal to $\mathbb{N} + A \times \mathbb{Z}$, see ref. [16]cite: Bovykin00 , for some linear order $A$

Thus one can choose Gödel encoding inside the standard model $M_{st}^{Z_2^{Hs}}$.

**Remark 1.1.19**. However there is no any problem as mentioned above in second order set theory $ZFC_2$ with the full second-order semantics because corresponding second order arithmetic $Z_2^{fss}$ is categorical.

**Remark 1.1.20**. Note if we view second-order arithmetic $Z_2$ as a theory in first-order predicate calculus. Thus a model $M^{Z_2}$ of the language of second-order arithmetic $Z_2$ consists of a set $M$ (which forms the range of individual variables) together with a constant 0 (an element of $M$), a function $S$ from $M$ to $M$, two binary operations + and $\times$ on $M$, a binary relation $<$ on $M$, and a collection $D$ of subsets of $M$, which is the range of the set variables. When $D$ is the full powerset of $M$, the model $M^{Z_2}$ is called a full model. The use of full second-order semantics is equivalent to limiting the models of second-order arithmetic to the full models. In fact, the axioms of second-order arithmetic have only one full model. This follows from the fact that the axioms of Peano arithmetic with the second-order induction axiom have only one model under second-order semantics, i.e. $Z_2$, with the full semantics, is categorical by Dedekind's argument, so has only one model up to isomorphism. When $M$ is the usual set of natural numbers with its usual operations, $M^{Z_2}$ is called an $\omega$-model. In this case we may identify the model with $D$, its collection of sets of naturals, because this set is enough to completely determine an $\omega$-model. The unique full omega-model $M_\omega^{Z_2^{fss}}$, which is the usual set of natural numbers with its usual structure and all its subsets, is called the intended or standard model of second-order arithmetic.

## 6mm . 2.Generalized Löb's theorem. Remarks on the Tarski's undefinability theorem 3mm

### 5mm . 2.1. Remarks on the Tarski's undefinability theorem 2mm

**Remark 2.1.1**. In paper [2]cite: Foukzon13 under the following assumption

$$Con(ZFC + \exists M_{st}^{ZFC}) \tag{14}$$

it has been proved that there exists countable Russell's set $\Re_\omega$ such that the followng statement is satisfied:

$$ZFC + \exists M_{st}^{ZFC} \vdash$$
$$\exists \Re_\omega (\Re_\omega \in M_{st}^{ZFC}) \wedge (card(\Re_\omega) = \aleph_0) \wedge \left[\vDash_{M_{st}^{ZFC}} \forall x (x \in \Re_\omega \Leftrightarrow x \notin x)\right]. \tag{15}$$

From (15) it immediately follows a contradiction

$$\vDash_{M_{st}^{ZFC}} (\Re_\omega \in \Re_\omega) \wedge (\Re_\omega \notin \Re_\omega). \tag{16}$$

From (16) and (14) by reductio and absurdum it follows

$$\neg Con(ZFC + \exists M_{st}^{ZFC}) \tag{17}$$

**Theorem 2.1.1**. [17]-[19]**cite: Mendelson97cite: Takeuti13cite: Quine37**. (Tarski's undefinability theorem). Let $\mathbf{Th}_{\mathcal{L}}$ be first order theory with formal language $\mathcal{L}$, which includes negation and has a Gödel numbering $g(\circ)$ such that for every $\mathcal{L}$-formula $A(x)$ there is a formula $B$ such that $B \leftrightarrow A(g(B))$ holds. Assume that $\mathbf{Th}_{\mathcal{L}}$ has a standard model $M_{st}^{\mathbf{Th}_{\mathcal{L}}}$ and $Con(\mathbf{Th}_{\mathcal{L},st})$ where

$$\mathbf{Th}_{\mathcal{L},st} \triangleq \mathbf{Th}_{\mathcal{L}} + \exists M_{st}^{\mathbf{Th}_{\mathcal{L}}}. \tag{18}$$

Let $T^*$ be the set of Gödel numbers of $\mathcal{L}$-sentences true in $M_{st}^{\mathbf{Th}_{\mathcal{L}}}$. Then there is no $\mathcal{L}$-formula $\mathbf{True}(n)$ (truth predicate) which defines $T^*$. That is, there is no $\mathcal{L}$-formula $\mathbf{True}(n)$ such that for every $\mathcal{L}$-formula $A$,

$$\mathbf{True}(g(A)) \Leftrightarrow [A]_{M_{st}^{\mathbf{Th}_{\mathcal{L}}}}, \tag{19}$$

where the abbreviation $[A]_{M_{st}^{\mathbf{Th}_{\mathcal{L}}}}$ means that $A$ holds in standard model $M_{st}^{\mathbf{Th}_{\mathcal{L}}}$, i.e. $[A]_{M_{st}^{\mathbf{Th}_{\mathcal{L}}}} \Leftrightarrow \vDash_{M_{st}^{\mathbf{Th}_{\mathcal{L}}}} A$. Therefore $Con(\mathbf{Th}_{\mathcal{L},st})$ implies that

$$\neg \exists \mathbf{True}(x) \left(\mathbf{True}(g(A)) \Leftrightarrow [A]_{M_{st}^{\mathbf{Th}_{\mathcal{L}}}}\right) \tag{20}$$

Thus Tarski's undefinability theorem reads

$$Con(\mathbf{Th}_{\mathcal{L},st}) \Rightarrow \neg \exists \mathbf{True}(x) \left(\mathbf{True}(g(A)) \Leftrightarrow [A]_{M_{st}^{\mathbf{Th}_{\mathcal{L}}}}\right). \tag{21}$$

**Remark 2.1.2**. (i) By the other hand the Theorem 2.1.1 says that given some really consistent formal theory $\mathbf{Th}_{\mathcal{L},st}$ that contins formal arithmetic, the concept of truth in that formal theory $\mathbf{Th}_{\mathcal{L},st}$ is not definable using the expressive means that that arithmetic affords. This implies a major limitation on the scope of "self-representation." It is possible to define a formula $\mathbf{True}(n)$, but only by drawing on a metalanguage whose expressive power goes beyond that of $\mathcal{L}$. To define a truth predicate for the metalanguage would require a still higher metametalanguage, and so on.

(ii) However if formal theory $\mathbf{Th}_{\mathcal{L},st}$ is inconsistent this is not surpriising if we define a formula $\mathbf{True}(n) = \mathbf{True}(n;\mathbf{Th}_{\mathcal{L},st})$ by drawing only on a language $\mathcal{L}$.

(iii) Note that if under assumption $Con(\mathbf{Th}_{\mathcal{L},st})$ we define a formula $\mathbf{True}(n;\mathbf{Th}_{\mathcal{L},st})$ by drawing only on a language $\mathcal{L}$ by reductio ad absurdum it follows

$$\neg Con(\mathbf{Th}_{\mathcal{L},st}). \tag{22}$$

**Remark 2.1.3**. (i) Let $ZFC_{st}$ be a theory $ZFC_{st} \triangleq ZFC + \exists M_{st}^{ZFC}$. In this paper under assumption $Con(ZFC_{st})$ we define a formula $\mathbf{True}(n;ZFC_{st})$ by drawing only on a language $\mathcal{L}_{ZFC_{st}}$ by using Generalized Löb's theorem [4], [5]cite: Foulzon17cite: Foukzon15. Thus by reductio ad absurdum it follows

$$\neg Con(ZFC + \exists M_{st}^{ZFC}). \tag{23}$$

(ii) However note that in this case we obtain $\neg Con(ZFC_{st})$ by using approach that completely different in comparison with approach based on derivation of the countable Russell's set $\Re_{\omega}$ with conditions (15).

5mm . **2.2.Generalized Löb's theorem** 2mm

**Definition 2.2.1**. Let $\mathbf{Th}_{\mathcal{L}}^{\#}$ be first order theory and $Con(\mathbf{Th}_{\mathcal{L}}^{\#})$. A theory $\mathbf{Th}_{\mathcal{L}}^{\#}$ is complete if, for every formula $A$ in the theory's language $\mathcal{L}$, that formula $A$ or its negation $\neg A$ is provable in $\mathbf{Th}_{\mathcal{L}}^{\#}$, i.e., for any wff $A$, always $\mathbf{Th}_{\mathcal{L}}^{\#} \vdash A$ or $\mathbf{Th}_{\mathcal{L}}^{\#} \vdash \neg A$.

**Definition 2.2.2**.Let $\mathbf{Th}_{\mathcal{L}}$ be first order theory and $Con(\mathbf{Th}_{\mathcal{L}})$. We will say that a theory $\mathbf{Th}_{\mathcal{L}}^{\#}$ is completion of the theory $\mathbf{Th}_{\mathcal{L}}$ if (i) $\mathbf{Th}_{\mathcal{L}} \subset \mathbf{Th}_{\mathcal{L}}^{\#}$, (ii) a theory $\mathbf{Th}_{\mathcal{L}}^{\#}$ is complete.

**Theorem 2.2.1**. [4], [5]cite: Foulzon17cite: Foukzon15. Assume that: $Con(ZFC_{st})$, where $ZFC_{st} \triangleq ZFC + \exists M_{st}^{ZFC}$. Then there exists completion $ZFC_{st}^{\#}$ of the theory $ZFC_{st}$ such that the following condtions holds:

(i) For every formula $A$ in the language of $ZFC$ that formula $[A]_{M_{st}^{ZFC}}$ or formula $[\neg A]_{M_{st}^{ZFC}}$ is provable in $ZFC_{st}^{\#}$ i.e., for any wff $A$, always $ZFC_{st}^{\#} \vdash [A]_{M_{st}^{ZFC}}$ or $ZFC_{st}^{\#} \vdash [\neg A]_{M_{st}^{ZFC}}$.

(ii) $ZFC_{st}^{\#} = \cup_{m\in\mathbb{N}} \mathbf{Th}_m$, where for any $m$ a theory $\mathbf{Th}_{m+1}$ is finite extension of the theory $\mathbf{Th}_m$.

(iii) Let $\Pr_m^{st}(y,x)$ be recursive relation such that: $y$ is a Gödel number of a proof of the wff of the theory $\mathbf{Th}_m$ and $x$ is a Gödel number of this wff. Then the relation $\Pr_m^{st}(y,x)$ is expressible in the theory $\mathbf{Th}_m$ by canonical Gödel encoding and really asserts provability in $\mathbf{Th}_m$.

(iv) Let $\Pr_{st}^{\#}(y,x)$ be relation such that: $y$ is a Gödel number of a proof of the wff of the theory $ZFC_{st}^{\#}$ and $x$ is a Gödel number of this wff. Then the relation $\Pr_{st}^{\#}(y,x)$is expressible in the theory $ZFC_{st}^{\#}$ by the following formula

$$\Pr_{st}^{\#}(y,x) \Leftrightarrow \exists m(m \in \mathbb{N}) \Pr_m^{st}(y,x)) \tag{24}$$

(v) The predicate $\Pr_{st}^{\#}(y,x)$ really asserts provability in the set theory $ZFC_{st}^{\#}$.

**Remark 2.2.1**. Note that the relation $\Pr_m^{st}(y,x)$ is expressible in the theory $\mathbf{Th}_m$ since a theory $\mathbf{Th}_m$ is an finite extension of the recursively axiomatizable theory $ZFC$ and therefore the predicate $\Pr_m^{st}(y,x)$ exists since any theory $\mathbf{Th}_m$ is recursively axiomatizable.

**Remark 2.2.2**. Note that a theory $ZFC_{st}^{\#}$ obviously is not recursively axiomatizable nevertheless Gödel encoding holds by Remark 2.2.1.

**Theorem 2.2.2**. Assume that: $Con(ZFC_{st})$, where $ZFC_{st} \triangleq ZFC + \exists M_{st}^{ZFC}$. Then truth predicate $\mathbf{True}(n)$ is expressible by using only first order language by the following formula

$$\mathbf{True}(g(A)) \Leftrightarrow \exists y(y \in \mathbb{N})\exists m(m \in \mathbb{N}) \Pr_m^{st}(y,g(A))). \tag{25}$$

Proof. Assume that:

$$ZFC_{st}^{\#} \vdash [A]_{M_{st}^{ZFC}}. \tag{26}$$

It follows from (26) there exists $m^* = m^*(g(A))$ such that $\mathbf{Th}_{m^*} \vdash [A]_{M_{st}^{ZFC}}$ and therefore by (24) we obtain

$$\Pr_{st}^{\#}(y, g(A)) \Leftrightarrow \Pr_{m^*}^{st}(y, g(A))). \qquad (27)$$

From (24) immediately by definitions one obtains (25).

**Remark 2.2.3**. Note that Theorem 2.1.1 in tis case reads

$$Con(ZFC_{st}) \Rightarrow \neg\exists \mathbf{True}(x) \left(\mathbf{True}(g(A)) \Leftrightarrow [A]_{M_{st}^{\mathbf{ZFC}}}\right). \qquad (28)$$

**Theorem 2.2.3**. $\neg Con(ZFC_{st})$.

Proof. Assume that: $Con(ZFC_{st})$. From (25) and (28) one obtains a condradiction $Con(ZFC_{st}) \wedge \neg Con(ZFC_{st})$ (see Remark 2.1.3) and therefore by *reductio ad absurdum* it follows $\neg Con(ZFC_{st})$.

**Theorem 2.2.4**. [4], [5]cite: Foulzon17cite: Foukzon15. Let $M_{Nst}^{ZFC}$ be a nonstandard model of $ZFC$ and let $M_{st}^{PA}$ be a standard model of $PA$.

We assume now that $M_{st}^{PA} \subset M_{Nst}^{ZFC}$ and denote such nonstandard model of the set theory $ZFC$ by $M_{Nst}^{ZFC} = M_{Nst}^{ZFC}[PA]$. Let $ZFC_{Nst}$ be the theory $ZFC_{Nst} = ZFC + M_{Nst}^{ZFC}[PA]$. Assume that: $Con(ZFC_{Nst})$, where $ZFC_{st} \triangleq ZFC + \exists M_{Nst}^{ZFC}$. Then there exists completion $ZFC_{Nst}^{\#}$ of the theory $ZFC_{Nst}$ such that the following condtions holds:

(i) For every formula $A$ in the language of $ZFC$ that formula $[A]_{M_{Nst}^{ZFC}}$ or formula $[\neg A]_{M_{Nst}^{ZFC}}$ is provable in $ZFC_{Nst}^{\#}$ i.e., for any wff $A$, always $ZFC_{Nst}^{\#} \vdash [A]_{M_{Nst}^{ZFC}}$ or $ZFC_{Nst}^{\#} \vdash [\neg A]_{M_{Nst}^{ZFC}}$.

(ii) $ZFC_{Nst}^{\#} = \bigcup_{m \in \mathbb{N}} \mathbf{Th}_m$, where for any $m$ a theory $\mathbf{Th}_{m+1}$ is finite extension of the theory $\mathbf{Th}_m$.

(iii) Let $\Pr_m^{Nst}(y, x)$ be recursive relation such that: $y$ is a Gödel number of a proof of the wff of the theory $\mathbf{Th}_m$ and $x$ is a Gödel number of this wff. Then the relation $\Pr_m^{Nst}(y, x)$ is expressible in the theory $\mathbf{Th}_m$ by canonical Gödel encoding and really asserts provability in $\mathbf{Th}_m$.

(iv) Let $\Pr_{Nst}^{\#}(y, x)$ be relation such that: $y$ is a Gödel number of a proof of the wff of the theory $ZFC_{Nst}^{\#}$ and $x$ is a Gödel number of this wff. Then the relation $\Pr_{Nst}^{\#}(y, x)$ is expressible in the theory $ZFC_{Nst}^{\#}$ by the following formula

$$\Pr_{Nst}^{\#}(y, x) \Leftrightarrow \exists m(m \in M_{st}^{PA}) \Pr_m^{Nst}(y, x)) \qquad (29)$$

(v) The predicate $\Pr_{Nst}^{\#}(y, x)$ really asserts provability in the set theory $ZFC_{Nst}^{\#}$.

**Remark 2.2.4**. Note that the relation $\Pr_m^{Nst}(y, x)$ is expressible in the theory $\mathbf{Th}_m$ since a theory $\mathbf{Th}_m$ is an finite extension of the recursively axiomatizable theory $ZFC$ and therefore the predicate $\Pr_m^{Nst}(y, x)$ exists since any theory $\mathbf{Th}_m$ is recursively axiomatizable.

**Remark 2.2.5**.Note that a theory $ZFC_{Nst}^{\#}$ obviously is not recursively axiomatizable nevertheless Gödel encoding holds by Remark 2.2.1.

**Theorem 2.2.5**.Assume that: $Con(ZFC_{Nst})$, where $ZFC_{Nst} \triangleq ZFC + \exists M_{Nst}^{ZFC}, M_{st}^{PA} \subset M_{Nst}^{ZFC}$.

Then truth predicate $\mathbf{True}(n)$ is expressible by using first order language by the following formula

$$\mathbf{True}(g(A)) \Leftrightarrow \exists y(y \in M_{st}^{PA}) \exists m(m \in M_{st}^{PA}) \Pr_m^{Nst}(y, g(A)). \qquad (30)$$

Proof. Assume that:

$$ZFC_{Nst}^{\#} \vdash [A]_{M_{Nst}^{ZFC}}. \qquad (31)$$

It follows from (29) there exists $m^* = m^*(g(A))$ such that $\mathbf{Th}_{m^*} \vdash [A]_{M_{Nst}^{ZFC}}$ and therefore by (31) we obtain

$$\Pr_{Nst}^{\#}(y, g(A)) \Leftrightarrow \Pr_{m^*}^{Nst}(y, g(A))). \qquad (32)$$

From (32) immediately by definitions one obtains (30).

**Remark 2.2.6**. Note that Theorem 2.1.1 in this case reads

$$Con(ZFC_{Nst}) \Rightarrow \neg\exists \mathbf{True}(x) \left(\mathbf{True}(g(A)) \Leftrightarrow [A]_{M_{Nst}^{\mathbf{ZFC}}}\right). \qquad (33)$$

**Theorem 2.2.6**. $\neg Con(ZFC_{Nst})$.

Proof. Assume that: $Con(ZFC_{Nst})$. From (30) and (33) one obtains a condradiction $Con(ZFC_{Nst}) \wedge \neg Con(ZFC_{Nst})$ and therefore by *reductio ad absurdum* it follows

$\neg Con(ZFC_{Nst})$.

**Theorem 2.2.7**. Assume that:$Con\left(\overline{ZFC}_2^{Hs}\right)$, where $\overline{ZFC}_2^{Hs} \triangleq ZFC_2^{Hs} + \exists M_{st}^{ZFC_2^{Hs}}$. Then there exists completion $\overline{ZFC}_2^{Hs\#}$ of the theory $\overline{ZFC}_2^{Hs}$ such that the following condtions holds:

(i) For every first order wff formula $A$ (wff$_1$ $A$) in the language of $ZFC_2^{Hs}$ that formula $[A]_{M_{st}^{ZFC_2^{Hs}}}$ or formula $[\neg A]_{M_{st}^{ZFC_2^{Hs}}}$ is provable in $\overline{ZFC}_2^{Hs\#}$ i.e., for any wff$_1$ $A$, always $\overline{ZFC}_2^{Hs\#} \vdash [A]_{M_{st}^{ZFC_2^{Hs}}}$ or $\overline{ZFC}_2^{Hs\#} \vdash [\neg A]_{M_{st}^{ZFC_2^{Hs}}}$.

(ii) $\overline{ZFC}_2^{Hs\#} = \bigcup_{m\in\mathbb{N}} \mathbf{Th}_m$, where for any $m$ a theory $\mathbf{Th}_{m+1}$ is finite extension of the theory $\mathbf{Th}_m$.

(iii) Let $\Pr_m^{st}(y,x)$ be recursive relation such that: $y$ is a Gödel number of a proof of the wff$_1$of the theory $\mathbf{Th}_m$ and $x$ is a Gödel number of this wff$_1$. Then the relation $\Pr_m^{st}(y,x)$ is expressible in the theory $\mathbf{Th}_m$ by canonical Gödel encoding and really asserts provability in $\mathbf{Th}_m$.

(iv) Let $\Pr_{st}^{\#}(y,x)$ be relation such that: $y$ is a Gödel number of a proof of the wff of the set theory $\overline{ZFC}_2^{Hs\#}$ and $x$ is a Gödel number of this wff$_1$. Then the relation $\Pr_{st}^{\#}(y,x)$ is expressible in the set theory $\overline{ZFC}_2^{Hs\#}$ by the following formula

$$\Pr_{st}^{\#}(y,x) \Leftrightarrow \exists m(m \in \mathbb{N}) \Pr_m^{st}(y,x)) \tag{34}$$

(v) The predicate $\Pr_{st}^{\#}(y,x)$ really asserts provability in the set theory $\overline{ZFC}_2^{Hs\#}$.

**Remark 2.2.7**. Note that the relation $\Pr_m^{st}(y,x)$ is expressible in the theory $\mathbf{Th}_m$ since a theory $\mathbf{Th}_m$ is an finite extension of the finite axiomatizable theory $ZFC_2^{Hs}$ and therefore the predicate $\Pr_m^{Nst}(y,x)$ exists since any theory $\mathbf{Th}_m$ is recursively axiomatizable.

**Remark 2.2.8**. Note that a theory $ZFC_{Nst}^{\#}$ obviously is not recursively axiomatizable nevertheless Gödel encoding holds by Remark 2.2.1.

**Theorem 2.2.8**. Assume that: $Con\left(\overline{ZFC}_2^{Hs}\right)$, where $\overline{ZFC}_2^{Hs} \triangleq ZFC_2^{Hs} + \exists M_{st}^{ZFC_2^{Hs}}$. Then truth predicate $\mathbf{True}(n)$ is expressible by using first order language by the following formula

$$\mathbf{True}(g(A)) \Leftrightarrow \exists y(y \in \mathbb{N})\exists m(m \in \mathbb{N}) \Pr_m^{st}(y,g(A)), \tag{35}$$

where $A$ is wff$_1$.

Proof. Assume that:

$$\overline{ZFC}_2^{Hs\#} \vdash [A]_{M_{st}^{ZFC_2^{Hs}}}. \tag{36}$$

It follows from (34) there exists $m^* = m^*(g(A))$ such that $\mathbf{Th}_{m^*} \vdash [A]_{M_{st}^{ZFC_2^{Hs}}}$ and therefore by (36) we obtain

$$\Pr_{st}^{\#}(y,g(A)) \Leftrightarrow \Pr_{m^*}^{st}(y,g(A))). \tag{37}$$

From (37) immediately by definitions one obtains (35).

**Remark 2.2.9**. Note that in considerd case Tarski's undefinability theorem (2.1.1) reads

$$Con\left(\overline{ZFC}_2^{Hs\#}\right) \Rightarrow \neg\exists \mathbf{True}(x) \left(\mathbf{True}(g(A)) \Leftrightarrow [A]_{M_{st}^{ZFC_2^{Hs}}}\right), \tag{38}$$

where $A$ is wff$_1$.

**Theorem 2.2.9**. $\neg Con\left(\overline{ZFC}_2^{Hs\#}\right)$.

Proof. Assume that: $Con\left(\overline{ZFC}_2^{Hs\#}\right)$. From (35) and (38) one obtains a condradiction $Con\left(\overline{ZFC}_2^{Hs\#}\right) \wedge \neg Con\left(\overline{ZFC}_2^{Hs\#}\right)$ and therefore by *reductio ad absurdum* it follows $\neg Con\left(\overline{ZFC}_2^{Hs\#}\right)$.

## 6mm . 3. Derivation of the inconsistent provably definable set in set theory $\overline{ZFC}_2^{Hs}, ZFC_{st}$ and $ZFC_{Nst}$ 3mm

### 5mm . 3.1. Derivation of the inconsistent provably definable set in set theory $\overline{ZFC}_2^{Hs}$ 2mm

**Definition 3.1.1.** (i) Let $\Phi$ be a wff of $\overline{ZFC}_2^{Hs}$. We will say that $\Phi$ is a first order $n$-place open wff if $\Phi$ contains free occurrences of the first order individual variables $X_1,\ldots,X_n$ and quantifiers only over any first order individual variables $Y_1,\ldots,Y_m$.

(ii) Let $\widetilde{\mathfrak{I}}_2^{Hs}$ be the countable set of the all first order provable definable sets $X$, i.e. sets such that $\overline{ZFC}_2^{Hs} \vdash \exists! X\Psi(X)$, where $\Psi(X) = \Psi_{M_{st}}(X)$ is a first order 1-place open wff that contains only first order variables (we will denote such wff for short by wff$_1$), with all bound variables restricted to standard model $M_{st} = M_{st}^{ZFC_2^{Hs}}$, i.e.

$$\forall Y\Big\{Y \in \widetilde{\mathfrak{I}}_2^{Hs} \Leftrightarrow \overline{ZFC}_2^{Hs} \vdash \exists \Psi_{M_{st}}(X)[([\Psi_{M_{st}}(X)] \in \Gamma_{X,M_{st}}^{Hs}/\sim_X) \wedge [\exists! X[\Psi_{M_{st}}(X) \wedge Y = X]]]\Big\}, \quad (39)$$

or in a short notatons

$$\forall Y\Big\{Y \in \widetilde{\mathfrak{I}}_2^{Hs} \Leftrightarrow \overline{ZFC}_2^{Hs} \vdash \exists \Psi(X)[([\Psi(X)] \in \Gamma_X^{Hs}/\sim_X) \wedge [\exists! X[\Psi(X) \wedge Y = X]]]\Big\}. \quad (40)$$

**Notation 3.1.1.** In this subsection we often write for short $\Psi(X), \mathcal{F}_X^{Hs}, \Gamma_X^{Hs}$ instead $\Psi_{M_{st}}(X)$, $\mathcal{F}_{X,M_{st}}^{Hs}, \Gamma_{X,M_{st}}^{Hs}$ but this should not lead to a confusion.

**Assumption 3.1.1.** We assume now for simplicity but without loss of generality that

$$\mathcal{F}_{X,M_{st}}^{Hs} \in M_{st} \quad (41)$$

and therefore by definition of model $M_{st} = M_{st}^{ZFC_2^{Hs}}$ one obtains $\Gamma_{X,M_{st}}^{Hs} \in M_{st}^{ZFC_2^{Hs}}$.

Let $X \notin_{\vdash_{\overline{ZFC}_2^{Hs}}} Y$ be a predicate such that $X \notin_{\vdash_{\overline{ZFC}_2^{Hs}}} Y \leftrightarrow \overline{ZFC}_2^{Hs} \vdash X \notin Y$. Let $\widetilde{\mathfrak{R}}_2^{Hs}$ be the countable set of the all sets such that

$$\forall X\Big(X \in \widetilde{\mathfrak{I}}_2^{Hs}\Big)\Big[X \in \widetilde{\mathfrak{R}}_2^{Hs} \leftrightarrow X \notin_{\vdash_{\overline{ZFC}_2^{Hs}}} X\Big]. \quad (42)$$

From (42) one obtains

$$\widetilde{\mathfrak{R}}_2^{Hs} \in \widetilde{\mathfrak{R}}_2^{Hs} \leftrightarrow \widetilde{\mathfrak{R}}_2^{Hs} \notin_{\vdash_{\overline{ZFC}_2^{Hs}}} \widetilde{\mathfrak{R}}_2^{Hs}. \quad (43)$$

But obviously (43) immediately gives a contradiction

$$\Big(\widetilde{\mathfrak{R}}_2^{Hs} \in \widetilde{\mathfrak{R}}_2^{Hs}\Big) \wedge \Big(\widetilde{\mathfrak{R}}_2^{Hs} \notin_{\vdash_{\overline{ZFC}_2^{Hs}}} \widetilde{\mathfrak{R}}_2^{Hs}\Big). \quad (44)$$

**Remark 3.1.1.** Note that a contradiction (44) in fact is a contradiction inside $\overline{ZFC}_2^{Hs}$ for the reason that predicate $X \notin_{\vdash_{\overline{ZFC}_2^{Hs}}} Y$ is expressible by first order language as predicate of $\overline{ZFC}_2^{Hs}$ (see subsection 1.2, Theorem 1.2.8 (ii)-(iii) and therefore countable sets $\widetilde{\mathfrak{I}}_2^{Hs}$ and $\widetilde{\mathfrak{R}}_2^{Hs}$ are sets in the sense of the set theory $\overline{ZFC}_2^{Hs}$.

**Remark 3.1.2.** Note that by using Gödel encoding the above stated contradiction can be shipped in special completion $\overline{ZFC}_2^{Hs\#}$ of $\overline{ZFC}_2^{Hs}$, see subsection 1.2, Theorem 1.2.8.

**Remark 3.1.3.** (i) Note that Tarski's undefinability theorem cannot bloked the equivalence (43) since this theorem is no longer holds by Proposition 2.2.1. (Generalized Löbs Theorem).

(ii) In additional note that: since Tarski's undefinability theorem has been proved under the same assumption $\exists M_{st}^{ZFC_2^{Hs}}$ by <u>*reductio ad absurdum*</u> it follows again $\neg Con(ZFC_{Nst})$, see Theorem 1.2.10.

**Remark 3.1.4**. More formally we can to explain the gist of the contradictions derived in this paper (see section 4) as follows.

Let $M$ be Henkin model of $ZFC_2^{Hs}$. Let $\widetilde{\mathfrak{R}}_2^{Hs}$ be the set of the all sets of $M$ provably definable in $\overline{ZFC}_2^{Hs}$, and let $\widetilde{\mathfrak{R}}_2^{Hs} = \left\{x \in \widetilde{\mathfrak{I}}_2^{Hs} : \Box(x \notin x)\right\}$ where $\Box A$ means 'sentence $A$ derivable in $\overline{ZFC}_2^{Hs}$', or some appropriate modification thereof. We replace now formula (39) by the following formula

$$\forall Y\left\{Y \in \widetilde{\mathfrak{I}}_2^{Hs} \leftrightarrow \exists\Psi(X)[([\Psi(X)] \in \Gamma_X^{Hs}/\sim_X) \wedge \Box\exists!X[\Psi(X) \wedge Y = X]]\right\}. \quad (45)$$

and we replace formula (42) by the following formula

$$\forall X\left(X \in \widetilde{\mathfrak{I}}_2^{Hs}\right)\left[X \in \widetilde{\mathfrak{R}}_2^{Hs} \Leftrightarrow \Box(X \notin X)\right]. \quad (46)$$

**Definition 3.1.2**. We rewrite now (45) in the following equivalent form

$$\forall Y\left\{Y \in \widetilde{\mathfrak{I}}_2^{Hs} \Leftrightarrow \exists\Psi(X)[([\Psi(X)]_{Hs} \in \Gamma_X^{\star Hs}/\sim_X) \wedge (Y = X)]\right\}, \quad (47)$$

where the countable set $\Gamma_X^{\star Hs}/\sim_X$ is defined by the following formula

$$\forall\Psi(X)\{[\Psi(X)] \in \Gamma_X^{\star Hs}/\sim_X \Leftrightarrow [([\Psi(X)]_{Hs} \in \Gamma_X^{Hs}/\sim_X) \wedge \Box\exists!X\Psi(X)]\} \quad (48)$$

**Definition 3.1.3**. Let $\widetilde{\mathfrak{R}}_2^{Hs}$ be the countable set of the all sets such that

$$\forall X\left(X \in \widetilde{\mathfrak{I}}_2^{Hs}\right)\left[X \in \widetilde{\mathfrak{R}}_2^{Hs} \Leftrightarrow \Box X \notin X\right]. \quad (49)$$

**Remark 3.1.5**. Note that $\widetilde{\mathfrak{R}}_2^{Hs} \in \widetilde{\mathfrak{I}}_2^{Hs}$ since $\widetilde{\mathfrak{R}}_2^{Hs}$ is a set definable by the first order 1-place open wff$_1$ :

$$\Psi\left(Z, \widetilde{\mathfrak{R}}_2^{Hs}\right) \triangleq \forall X\left(X \in \widetilde{\mathfrak{I}}_2^{Hs}\right)[X \in Z \Leftrightarrow \Box(X \notin X)]. \quad (50)$$

From (49) and Remark 3.1.4 one obtains

$$\widetilde{\mathfrak{R}}_2^{Hs} \in \widetilde{\mathfrak{R}}_2^{Hs} \Leftrightarrow \Box\left(\widetilde{\mathfrak{R}}_2^{Hs} \notin \widetilde{\mathfrak{R}}_2^{Hs}\right). \quad (51)$$

But (51) immediately gives a contradiction

$$\overline{ZFC}_2^{Hs} \vdash \left(\widetilde{\mathfrak{R}}_2^{Hs} \in \widetilde{\mathfrak{R}}_2^{Hs}\right) \wedge \left(\widetilde{\mathfrak{R}}_2^{Hs} \notin \widetilde{\mathfrak{R}}_2^{Hs}\right). \quad (52)$$

**Remark 3.1.6**. Note that contradiction (52) is a contradiction inside $\overline{ZFC}_2^{Hs}$ for the reason that the countable set $\widetilde{\mathfrak{I}}_2^{Hs}$ is a set in the sense of the set theory $\overline{ZFC}_2^{Hs}$.

In order to obtain a contradiction inside $\overline{ZFC}_2^{Hs}$ without any reference to Assumption 3.1.1 we introduce the following definitions.

**Definition 3.1.4**. We define now the countable set $\Gamma_v^{\star Hs}/\sim_v$ by the following formula

$$\forall y\left\{[y]_{Hs} \in \Gamma_v^{\star Hs}/\sim_v \Leftrightarrow ([y]_{Hs} \in \Gamma_v^{Hs}/\sim_v) \wedge \widehat{\mathbf{Fr}}_2^{Hs}(y,v) \wedge [\Box\exists!X\Psi_{y,v}(X)]\right\}. \quad (53)$$

**Definition 3.1.5**. We choose now $\Box A$ in the following form

$$\Box A \triangleq Bew_{\overline{ZFC}_2^{Hs}}(\#A) \wedge \left[Bew_{\overline{ZFC}_2^{Hs}}(\#A) \Rightarrow A\right]. \quad (54)$$

Here $Bew_{\overline{ZFC}_2^{Hs}}(\#A)$ is a canonical Gödel formula which says to us that there exists proof in $\overline{ZFC}_2^{Hs}$ of the formula $A$ with Gödel number $\#A$.

**Remark 3.1.7**. Note that the Definition 3.1.5 holds as definition of predicate really asserting provability of the first order sentence $A$ in $\overline{ZFC}_2^{Hs}$.

**Definition 3.1.6**. Using Definition 3.1.5, we replace now formula (48) by the following formula

$$\forall\Psi(X)\{[\Psi(X)] \in \Gamma_X^{\star Hs}/\sim_X \iff \exists\Psi(X)([\Psi(X)] \in \Gamma_X^{Hs}/\sim_X)\wedge$$
$$\wedge\left[Bew_{\overline{ZFC}_2^{Hs}}(\#(\exists!X[\Psi(X)\wedge Y=X]))\right]\wedge \quad (55)$$
$$\wedge\left[Bew_{\overline{ZFC}_2^{Hs}}(\#(\exists!X[\Psi(X)\wedge Y=X])) \Rightarrow \exists!X[\Psi(X)\wedge Y=X]\right]\Big\}.$$

**Definition 3.1.7**. Using Definition 3.1.5, we replace now formula (49) by the following formula

$$\forall X\left(X \in \widetilde{\mathfrak{I}}_2^{Hs}\right)\left[X \in \widetilde{\mathfrak{R}}_2^{Hs} \iff \left[Bew_{\overline{ZFC}_2^{Hs}}(\#(X \notin X))\right]\wedge\right.$$
$$\left.\wedge\left[Bew_{\overline{ZFC}_2^{Hs}}(\#(X \notin X)) \Rightarrow X \notin X\right]\right]. \quad (56)$$

**Definition 3.1.8**. Using Proposition 2.1.1 and Remark 2.1.10 [6]cite: Foukzon19, we replace now formula (53) by the following formula

$$\forall y\{[y]_{Hs} \in \Gamma_v^{\star Hs}/\sim_v \iff$$
$$([y]_{Hs} \in \Gamma_v^{Hs}/\sim_v)\wedge\widehat{\mathbf{Fr}_2^{Hs}}(y,v)\wedge\left[Bew_{\overline{ZFC}_2^{Hs}}(\#\exists!X[\Psi_{y,v}(X)\wedge Y=X])\right]\wedge \quad (57)$$
$$\wedge\left[Bew_{\overline{ZFC}_2^{Hs}}(\#\exists!X[\Psi_{y,v}(X)\wedge Y=X]) \Rightarrow \exists!X[\Psi_{y,v}(X)\wedge Y=X]\right]\Big\}.$$

**Definition 3.1.9**. Using Definitions 3.1.4-3.1.6, we define now the countable set $\widetilde{\mathfrak{I}}_2^{\star Hs}$ by formula

$$\forall Y\left\{Y \in \widetilde{\mathfrak{I}}_2^{\star Hs} \iff \exists y\left[([y] \in \Gamma_v^{\star Hs}/\sim_v)\wedge\left(g_{\overline{ZFC}_2^{Hs}}(X)=v\right)\right]\right\}. \quad (58)$$

**Remark 3.1.8**. Note that from the second order axiom schema of replacement (12) it follows directly that $\widetilde{\mathfrak{I}}_2^{\star Hs}$ is a set in the sense of the set theory $\overline{ZFC}_2^{Hs}$.

**Definition 3.1.10**. Using Definition 3.1.8 we replace now formula (56) by the following formula

$$\forall X\left(X \in \widetilde{\mathfrak{I}}_2^{\star Hs}\right)$$
$$\left[X \in \widetilde{\mathfrak{R}}_2^{\star Hs} \iff \left[Bew_{\overline{ZFC}_2^{Hs}}(\#(X \notin X))\right]\wedge\left[Bew_{\overline{ZFC}_2^{Hs}}(\#(X \notin X)) \Rightarrow X \notin X\right]\right]. \quad (59)$$

**Remark 3.1.9**. Notice that the expression (60)

$$\left[Bew_{\overline{ZFC}_2^{Hs}}(\#(X \notin X))\right]\wedge\left[Bew_{\overline{ZFC}_2^{Hs}}(\#(X \notin X)) \Rightarrow X \notin X\right] \quad (60)$$

obviously is a well formed formula of $\overline{ZFC}_2^{Hs}$ and therefore a set $\widetilde{\mathfrak{R}}_2^{\star Hs}$ is a set in the sense of $\overline{ZFC}_2^{Hs}$.

**Remark 3.1.10**. Note that $\widetilde{\mathfrak{R}}_2^{\star Hs} \in \widetilde{\mathfrak{I}}_2^{\star Hs}$ since $\widetilde{\mathfrak{R}}_2^{\star Hs}$ is a set definable by 1-place open wff

$$\Psi\left(Z,\widetilde{\mathfrak{R}}_2^{\star Hs}\right) \triangleq$$
$$\forall X\left(X \in \widetilde{\mathfrak{I}}_2^{Hs}\right)[X \in Z \iff \quad (61)$$
$$\left[Bew_{\overline{ZFC}_2^{Hs}}(\#(X \notin X))\right]\wedge\left[Bew_{\overline{ZFC}_2^{Hs}}(\#(X \notin X)) \Rightarrow X \notin X\right]].$$

**Theorem 3.1.1**. Set theory $\overline{ZFC}_2^{Hs} \triangleq ZFC_2^{Hs} + \exists M_{st}^{ZFC_2^{Hs}}$ is inconsistent.

Proof. From (59) we obtain

$$\widetilde{\Re}_2^{\star Hs} \in \widetilde{\Re}_2^{\star Hs} \Leftrightarrow \left[ Bew_{\overline{ZFC}_2^{Hs}}\left(\#\left(\widetilde{\Re}_2^{\star Hs} \notin \widetilde{\Re}_2^{\star Hs}\right)\right)\right] \wedge$$
$$\wedge \left[ Bew_{\overline{ZFC}_2^{Hs}}\left(\#\left(\widetilde{\Re}_2^{\star Hs} \notin \widetilde{\Re}_2^{\star Hs}\right)\right) \Rightarrow \widetilde{\Re}_2^{\star Hs} \notin \widetilde{\Re}_2^{\star Hs}\right]. \quad (62)$$

(a) Assume now that:

$$\widetilde{\Re}_2^{\star Hs} \in \widetilde{\Re}_2^{\star Hs}. \quad (63)$$

Then from (62) we obtain $\vdash_{\overline{ZFC}_2^{Hs}} Bew_{\overline{ZFC}_2^{Hs}}\left(\#\left(\widetilde{\Re}_2^{\star Hs} \notin \widetilde{\Re}_2^{\star Hs}\right)\right)$ and $\vdash_{\overline{ZFC}_2^{Hs}} Bew_{\overline{ZFC}_2^{Hs}}\left(\#\left(\widetilde{\Re}_2^{\star Hs} \notin \widetilde{\Re}_2^{\star Hs}\right)\right) \Rightarrow \widetilde{\Re}_2^{\star Hs} \notin \widetilde{\Re}_2^{\star Hs}$, therefore $\vdash_{\overline{ZFC}_2^{Hs}} \widetilde{\Re}_2^{\star Hs} \notin \widetilde{\Re}_2^{\star Hs}$ and so

$$\vdash_{\overline{ZFC}_2^{Hs}} \widetilde{\Re}_2^{\star Hs} \in \widetilde{\Re}_2^{\star Hs} \Rightarrow \widetilde{\Re}_2^{\star Hs} \notin \widetilde{\Re}_2^{\star Hs}. \quad (64)$$

From (63)-(64) we obtain

$$\widetilde{\Re}_2^{\star Hs} \in \widetilde{\Re}_2^{\star Hs}, \widetilde{\Re}_2^{\star Hs} \in \widetilde{\Re}_2^{\star Hs} \Rightarrow \widetilde{\Re}_2^{\star Hs} \notin \widetilde{\Re}_2^{\star Hs} \vdash \widetilde{\Re}_2^{\star Hs} \notin \widetilde{\Re}_2^{\star Hs}$$

and thus $\vdash_{\overline{ZFC}_2^{Hs}} \left(\widetilde{\Re}_2^{Hs} \in \widetilde{\Re}_2^{Hs}\right) \wedge \left(\widetilde{\Re}_2^{Hs} \notin \widetilde{\Re}_2^{Hs}\right)$.

(b) Assume now that

$$\left[ Bew_{\overline{ZFC}_2^{Hs}}\left(\#\left(\widetilde{\Re}_2^{\star Hs} \notin \widetilde{\Re}_2^{\star Hs}\right)\right)\right] \wedge$$
$$\wedge \left[ Bew_{\overline{ZFC}_2^{Hs}}\left(\#\left(\widetilde{\Re}_2^{\star Hs} \notin \widetilde{\Re}_2^{\star Hs}\right)\right) \Rightarrow \widetilde{\Re}_2^{\star Hs} \notin \widetilde{\Re}_2^{\star Hs}\right]. \quad (65)$$

Then from (65) we obtain $\vdash \widetilde{\Re}_2^{\star Hs} \notin \widetilde{\Re}_2^{\star Hs}$. From (65) and (62) we obtain $\vdash_{\overline{ZFC}_2^{Hs}} \widetilde{\Re}_2^{\star Hs} \in \widetilde{\Re}_2^{\star Hs}$, so $\vdash_{\overline{ZFC}_2^{Hs}} \widetilde{\Re}_2^{\star Hs} \notin \widetilde{\Re}_2^{\star Hs}, \widetilde{\Re}_2^{\star Hs} \in \widetilde{\Re}_2^{\star Hs}$ which immediately gives us a contradiction $\vdash_{\overline{ZFC}_2^{Hs}} \left(\widetilde{\Re}_2^{\star Hs} \in \widetilde{\Re}_2^{\star Hs}\right) \wedge \left(\widetilde{\Re}_2^{\star Hs} \notin \widetilde{\Re}_2^{\star Hs}\right)$.

**Definition 3.1.11**. We choose now $\square A$ in the following form

$$\square A \triangleq \overline{Bew}_{\overline{ZFC}_2^{Hs}}(\#A), \quad (66)$$

or in the following equivalent form

$$\square A \triangleq \overline{Bew}_{\overline{ZFC}_2^{Hs}}(\#A) \wedge \left[\overline{Bew}_{\overline{ZFC}_2^{Hs}}(\#A) \Rightarrow A\right] \quad (67)$$

similar to (46). Here $\overline{Bew}_{\overline{ZFC}_2^{Hs}}(\#A)$ is a Gödel formula which really asserts provability in $\overline{ZFC}_2^{Hs}$ of the formula $A$ with Gödel number $\#A$.

**Remark 3.1.11**. Notice that the Definition 3.1.12 with formula (66) holds as definition of predicate really asserting provability in $\overline{ZFC}_2^{Hs}$.

**Definition 3.1.12**. Using Definition 3.1.11 with formula (66), we replace now formula (48) by the following formula

$$\forall\Psi(X)\Big\{[\Psi(X)] \in \overline{\Gamma}_X^{\star Hs}/\sim_X \Leftrightarrow \exists\Psi(X)([\Psi(X)] \in \Gamma_X^{Hs}/\sim_X) \wedge$$
$$\wedge\left[\overline{Bew}_{\overline{ZFC}_2^{Hs}}(\#(\exists!X[\Psi(X) \wedge Y = X]))\right]\Big\}. \quad (68)$$

**Definition 3.1.13**. Using Definition 3.1.11 with formula (66), we replace now formula (49) by the following formula

$$\forall X\left(X \in \widetilde{\Im}_2^{Hs}\right)\left[X \in \widetilde{\Re}_2^{Hs} \Leftrightarrow \left[\overline{Bew}_{\overline{ZFC}_2^{Hs}}(\#(X \notin X))\right]\right] \quad (69)$$

**Definition 3.1.14**. Using Definition 3.1.11 with formula (66),we replace now formula (53) by the following formula

$$\forall y\{[y]_{Hs} \in \Gamma_v^{\star Hs}/\sim_v \iff$$
$$([y]_{Hs} \in \Gamma_v^{Hs}/\sim_v) \wedge \widehat{\mathbf{Fr}}_2^{Hs}(y,v) \wedge \left[\overline{Bew}_{\overline{ZFC}_2^{Hs}}(\#\exists! X[\Psi_{y,v}(X) \wedge Y = X])\right]\Big\}. \tag{70}$$

**Definition 3.1.15**. Using Definitions 3.1.12-3.1.16, we define now the countable set $\widetilde{\mathfrak{I}}_2^{\star Hs}$ by formula

$$\forall Y\left\{Y \in \widetilde{\mathfrak{I}}_2^{\star Hs} \iff \exists y\left[([y] \in \Gamma_v^{\star Hs}/\sim_v) \wedge \left(g_{\overline{ZFC}_2^{Hs}}(X) = v\right)\right]\right\}. \tag{71}$$

**Remark 3.1.12**. Note that from the axiom schema of replacement (12) it follows directly that $\widetilde{\mathfrak{I}}_2^{\star Hs}$ is a set in the sense of the set theory $\overline{ZFC}_2^{Hs}$.

**Definition 3.1.16**. Using Definition 3.1.15 we replace now formula (69) by the following formula

$$\forall X\left(X \in \widetilde{\mathfrak{I}}_2^{\star Hs}\right)\left[X \in \widetilde{\mathfrak{R}}_2^{\star Hs} \iff \left[\overline{Bew}_{\overline{ZFC}_2^{Hs}}(\#(X \notin X))\right]\right]. \tag{72}$$

**Remark 3.1.13**. Notice that the expressions (73)

$$\left[\overline{Bew}_{\overline{ZFC}_2^{Hs}}(\#(X \notin X))\right]$$
$$\text{and} \tag{73}$$
$$\left[\overline{Bew}_{\overline{ZFC}_2^{Hs}}(\#(X \notin X))\right] \wedge \left[\overline{Bew}_{\overline{ZFC}_2^{Hs}}(\#(X \notin X)) \Rightarrow X \notin X\right]$$

obviously are a well formed formula of $\overline{ZFC}_2^{Hs}$ and therefore collection $\widetilde{\mathfrak{R}}_2^{\star Hs}$ is a set in the sense of $\overline{ZFC}_2^{Hs}$.

**Remark 3.1.14**. Note that $\widetilde{\mathfrak{R}}_2^{\star Hs} \in \widetilde{\mathfrak{I}}_2^{\star Hs}$ since $\widetilde{\mathfrak{R}}_2^{\star Hs}$ is a set definable by 1-place open wff$_1$

$$\Psi\left(Z, \widetilde{\mathfrak{R}}_2^{\star Hs}\right) \triangleq \forall X\left(X \in \widetilde{\mathfrak{I}}_2^{\star Hs}\right)\left[X \in Z \iff \overline{Bew}_{\overline{ZFC}_2^{Hs}}(\#(X \notin X))\right]. \tag{74}$$

**Theorem 3.1.2**. Set theory $\overline{ZFC}_2^{Hs} \triangleq ZFC_2^{Hs} + \exists M_{st}^{ZFC_2^{Hs}}$ is inconsistent.

Proof. From (72) we obtain

$$\widetilde{\mathfrak{R}}_2^{\star Hs} \in \widetilde{\mathfrak{R}}_2^{\star Hs} \iff \left[\overline{Bew}_{\overline{ZFC}_2^{Hs}}\left(\#\left(\widetilde{\mathfrak{R}}_2^{\star Hs} \notin \widetilde{\mathfrak{R}}_2^{\star Hs}\right)\right)\right]. \tag{75}$$

(a) Assume now that:

$$\widetilde{\mathfrak{R}}_2^{\star Hs} \in \widetilde{\mathfrak{R}}_2^{\star Hs}. \tag{76}$$

Then from (75) we obtain $\vdash_{\overline{ZFC}_2^{Hs}} \overline{Bew}_{\overline{ZFC}_2^{Hs}}\left(\#\left(\widetilde{\mathfrak{R}}_2^{\star Hs} \notin \widetilde{\mathfrak{R}}_2^{\star Hs}\right)\right)$ and therefore $\vdash_{\overline{ZFC}_2^{Hs}} \widetilde{\mathfrak{R}}_2^{\star Hs} \notin \widetilde{\mathfrak{R}}_2^{\star Hs}$,
thus we obtain

$$\vdash_{\overline{ZFC}_2^{Hs}} \widetilde{\mathfrak{R}}_2^{\star Hs} \in \widetilde{\mathfrak{R}}_2^{\star Hs} \Rightarrow \widetilde{\mathfrak{R}}_2^{\star Hs} \notin \widetilde{\mathfrak{R}}_2^{\star Hs}. \tag{77}$$

From (76)-(77) we obtain $\widetilde{\mathfrak{R}}_2^{\star Hs} \in \widetilde{\mathfrak{R}}_2^{\star Hs}$ and $\widetilde{\mathfrak{R}}_2^{\star Hs} \in \widetilde{\mathfrak{R}}_2^{\star Hs} \Rightarrow \widetilde{\mathfrak{R}}_2^{\star Hs} \notin \widetilde{\mathfrak{R}}_2^{\star Hs}$, thus $\vdash_{\overline{ZFC}_2^{Hs}} \widetilde{\mathfrak{R}}_2^{\star Hs} \notin \widetilde{\mathfrak{R}}_2^{\star Hs}$ and finally we obtain $\vdash_{\overline{ZFC}_2^{Hs}} \left(\widetilde{\mathfrak{R}}_2^{Hs} \in \widetilde{\mathfrak{R}}_2^{Hs}\right) \wedge \left(\widetilde{\mathfrak{R}}_2^{Hs} \notin \widetilde{\mathfrak{R}}_2^{Hs}\right)$.

(b) Assume now that

$$\left[Bew_{\overline{ZFC}_2^{Hs}}\left(\#\left(\widetilde{\mathfrak{R}}_2^{\star Hs} \notin \widetilde{\mathfrak{R}}_2^{\star Hs}\right)\right)\right]. \tag{78}$$

Then from (78) we obtain $\vdash_{\overline{ZFC}_2^{Hs}} \widetilde{\mathfrak{R}}_2^{\star Hs} \notin \widetilde{\mathfrak{R}}_2^{\star Hs}$. From (78) and (75) we obtain

$\vdash_{\overline{ZFC}_2^{Hs}} \widetilde{\mathfrak{R}}_2^{\star Hs} \in \widetilde{\mathfrak{R}}_2^{\star Hs}$, thus $\vdash_{\overline{ZFC}_2^{Hs}} \widetilde{\mathfrak{R}}_2^{\star Hs} \notin \widetilde{\mathfrak{R}}_2^{\star Hs}$ and $\vdash_{\overline{ZFC}_2^{Hs}} \widetilde{\mathfrak{R}}_2^{\star Hs} \in \widetilde{\mathfrak{R}}_2^{\star Hs}$ which immediately gives us a contradiction $\vdash_{\overline{ZFC}_2^{Hs}} \left(\widetilde{\mathfrak{R}}_2^{\star Hs} \in \widetilde{\mathfrak{R}}_2^{\star Hs}\right) \wedge \left(\widetilde{\mathfrak{R}}_2^{\star Hs} \notin \widetilde{\mathfrak{R}}_2^{\star Hs}\right)$.

5mm . **3.2. Derivation of the inconsistent provably definable set in set theory** $ZFC_{st}$ 2mm

Let $\Im_{st}$ be the countable set of all sets $X$ such that $ZFC_{st} \vdash \exists! X\Psi(X),$ where $\Psi(X)$ is a 1-place open wff of $ZFC$ i.e.,

$$\forall Y\{Y \in \Im_{st} \Leftrightarrow ZFC_{st} \vdash \exists\Psi(X)[([\Psi(X)] \in \Gamma_X^{st}/\sim_X) \wedge \exists! X[\Psi(X) \wedge Y = X]]\}. \quad (79)$$

Let $X \notin_{\vdash_{ZFC_{st}}} Y$ be a predicate such that $X \notin_{\vdash_{ZFC_{st}}} Y \Leftrightarrow ZFC_{st} \vdash X \notin Y.$ Let $\Re$ be the countable set of the all sets such that

$$\forall X\Big[X \in \Re_{st} \Leftrightarrow (X \in \Im_{st}) \wedge \Big(X \notin_{\vdash_{ZFC_{st}}} X\Big)\Big]. \quad (80)$$

From (80) one obtains

$$\Re_{st} \in \Re_{st} \Leftrightarrow \Re_{st} \notin_{\vdash_{ZFC_{st}}} \Re_{st}. \quad (81)$$

But (81) immediately gives a contradiction

$$(\Re_{st} \in \Re_{st}) \wedge (\Re_{st} \notin \Re_{st}). \quad (82)$$

**Remark 3.2.1**. Note that a contradiction (82) is a contradiction inside $ZFC_{st}$ for the reason that predicate $X \notin_{\vdash_{ZFC_{st}}} Y$ is expressible by using first order language as predicate of $ZFC_{st}$ (see subsection 4.1) and therefore countable sets $\Im_{st}$ and $\Re_{st}$ are sets in the sense of the set theory $ZFC_{st}$.

**Remark 3.2.2**. Note that by using Gödel encoding the above stated contradiction can be shipped in special completion $ZFC_{st}^{\#}$ of $ZFC_{st},$ see subsection 1.2, Theorem 1.2.2 (i).

**Designation 3.2.1** (i) Let $M_{st}^{ZFC}$ be a standard model of $ZFC$ and

(ii) let $ZFC_{st}$ be the theory $ZFC_{st} = ZFC + \exists M_{st}^{ZFC},$

(iii) let $\Im_{st}$ be the set of the all sets of $M_{st}^{ZFC}$ provably definable in $ZFC_{st},$ and let $\Re_{st} = \{X \in \Im_{st} : \Box_{st}(X \notin X)\},$ where $\Box_{st}A$ means: 'sentence $A$ derivable in $ZFC_{st}$', or some appropriate modification thereof.

We replace now (79) by formula

$$\forall Y\{Y \in \Im_{st} \leftrightarrow \Box_{st}[\exists\Psi(\cdot)\exists! X[\Psi(X) \wedge Y = X]]\}, \quad (83)$$

and we replace (80) by formula

$$\forall X\Big[X \in \Re_{st} \leftrightarrow (X \in \Im_{st}) \wedge \Box_{st}\Big(X \notin X\Big)\Big]. \quad (84)$$

Assume that $ZFC_{st} \vdash \Re_{st} \in \Im_{st}.$ Then, we have that: $\Re_{st} \in \Re_{st}$ iff $\Box_{st}(\Re_{st} \notin \Re_{st}),$ which immediately gives us $\Re_{st} \in \Re_{st}$ iff $\Re_{st} \notin \Re_{st}.$ But this is a contradiction, i.e., $ZFC_{st} \vdash (\Re_{st} \in \Re_{st}) \wedge (\Re_{st} \notin \Re_{st}).$ We choose now $\Box_{st}A$ in the following form

$$\Box_{st}A \triangleq Bew_{ZFC_{st}}(\#A) \wedge [Bew_{ZFC_{st}}(\#A) \Rightarrow A]. \quad (85)$$

Here $Bew_{ZFC_{st}}(\#A)$ is a canonycal Gödel formula which says to us that there exists proof in $ZFC_{st}$ of the formula $A$ with Gödel number $\#A \in M_{st}^{PA}$.

**Remark 3.2.3**. Notice that Definition 3.2.6 holds as definition of predicate really asserting provability in $ZFC_{st}$.

**Definition 3.2.1**. We rewrite now (83) in the following equivalent form

$$\forall Y\Big\{Y \in \widetilde{\Im}_{st} \Leftrightarrow \exists\Psi(X)[([\Psi(X)]_{st} \in \Gamma_X^{\star st}/\sim_X) \wedge (Y = X)]\Big\}, \quad (86)$$

where the countable collection $\Gamma_X^{\star Hs}/\sim_X$ is defined by the following formula

$$\forall\Psi(X)\{[\Psi(X)]_{st} \in \Gamma_X^{\star st}/\sim_X \Leftrightarrow [([\Psi(X)]_{st} \in \Gamma_X^{st}/\sim_X) \wedge \Box_{st}\exists! X\Psi(X)]\} \quad (87)$$

**Definition 3.2.2**. Let $\widetilde{\Re}_{st}$ be the countable collection of the all sets such that

$$\forall X\Big(X \in \widetilde{\Im}_{st}\Big)\Big[X \in \widetilde{\Re}_{st} \Leftrightarrow \Box_{st}(X \notin X)\Big]. \quad (88)$$

**Remark 3.2.4**. Note that $\widetilde{\Re}_2^{Hs} \in \widetilde{\Im}_2^{Hs}$ since $\widetilde{\Re}_2^{Hs}$ is a collection definable by 1-place open wff

$$\Psi\left(Z,\widetilde{\Re}_{st}\right) \triangleq \forall X\left(X \in \widetilde{\Im}_{st}\right)[X \in Z \Leftrightarrow \Box_{st}(X \notin X)]. \tag{89}$$

**Definition 3.2.3**. By using formula (85) we rewrite now (86) in the following equivalent form

$$\forall Y\left\{Y \in \widetilde{\Im}_{st} \Leftrightarrow \exists\Psi(X)[([\Psi(X)]_{st} \in \Gamma_X^{\star st}/\sim_X) \wedge (Y = X)]\right\}, \tag{90}$$

where the countable collection $\Gamma_X^{\star Hs}/\sim_X$ is defined by the following formula

$$\begin{gathered}\forall\Psi(X)\{[\Psi(X)]_{st} \in \Gamma_X^{\star st}/\sim_X \Leftrightarrow \\ [([\Psi(X)]_{st} \in \Gamma_X^{st}/\sim_X) \wedge Bew_{ZFC_{st}}(\#\exists!X\Psi(X))] \wedge \\ \wedge[Bew_{ZFC_{st}}(\#\exists!X\Psi(X)) \Rightarrow \exists!X\Psi(X)]\}\end{gathered} \tag{91}$$

**Definition 3.2.4**. Using formula (85), we replace now formula (88) by the following formula

$$\begin{gathered}\forall X\left(X \in \widetilde{\Im}_{st}\right)\left[X \in \widetilde{\Re}_{st} \Leftrightarrow [Bew_{ZFC_{st}}(\#(X \notin X))] \wedge \right. \\ \wedge[Bew_{ZFC_{st}}(\#(X \notin X))].\end{gathered} \tag{92}$$

**Definition 3.2.5**. Using Proposition 2.1.1. and Remark 2.2.2 [6]cite: Foukzon19, we replace now formula (89) by the following formula

$$\begin{gathered}\forall y\{[y]_{st} \in \Gamma_v^{\star st}/\sim_v \Leftrightarrow \\ ([y]_{st} \in \Gamma_v^{st}/\sim_v) \wedge \widehat{\mathbf{Fr}}_{st}(y,v) \wedge [Bew_{ZFC_{st}}(\#\exists!X[\Psi_{y,v}(X) \wedge Y = X])] \wedge \\ \wedge[Bew_{ZFC_{st}}(\#\exists!X[\Psi_{y,v}(X) \wedge Y = X]) \Rightarrow \exists!X[\Psi_{y,v}(X) \wedge Y = X]]\}.\end{gathered} \tag{93}$$

**Definition 3.2.6**. Using Definitions 3.2.3-3.2.5, we define now the countable set $\widetilde{\Im}_{st}^{\star}$ by formula

$$\forall Y\left\{Y \in \widetilde{\Im}_{st}^{\star} \Leftrightarrow \exists y[([y]_{st} \in \Gamma_v^{\star st}/\sim_v) \wedge (g_{ZFC_{st}}(X) = v)]\right\}. \tag{94}$$

**Remark 3.2.5**. Note that from the axiom schema of replacement it follows directly that $\widetilde{\Im}_{st}^{\star}$ is a set in the sense of the set theory $ZFC_{st}$.

**Definition 3.2.7**. Using Definition 3.2.6 we replace now formula (92) by the following formula

$$\begin{gathered}\forall X\left(X \in \widetilde{\Im}_{st}^{\star}\right) \\ \left[X \in \widetilde{\Re}_{st}^{\star} \Leftrightarrow [Bew_{ZFC_{st}}(\#(X \notin X))] \wedge [Bew_{\overline{ZFC}_{st}}(\#(X \notin X)) \Rightarrow X \notin X]\right].\end{gathered} \tag{95}$$

**Remark 3.2.6**. Notice that the expression (96)

$$[Bew_{ZFC_{st}}(\#(X \notin X))] \wedge [Bew_{ZFC_{st}}(\#(X \notin X)) \Rightarrow X \notin X] \tag{96}$$

obviously is a well formed formula of $ZFC_{st}$ and therefore collection $\widetilde{\Re}_{st}^{\star}$ is a set in the sense of $\overline{ZFC}_2^{Hs}$.

**Remark 3.2.7**.Note that $\widetilde{\Re}_{st}^{\star} \in \widetilde{\Im}_{st}^{\star}$ since $\widetilde{\Re}_{st}^{\star}$ is a collection definable by 1-place open wff

$$\begin{gathered}\Psi\left(Z,\widetilde{\Re}_{st}^{\star}\right) \triangleq \\ \forall X\left(X \in \widetilde{\Im}_{st}^{\star}\right)[X \in Z \Leftrightarrow \\ [Bew_{ZFC_{st}}(\#(X \notin X))] \wedge [Bew_{ZFC_{st}}(\#(X \notin X)) \Rightarrow X \notin X]].\end{gathered} \tag{97}$$

**Theorem 3.2.1**. Set theory $ZFC_{st} \triangleq ZFC + \exists M_{st}^{ZFC}$ is inconsistent.

Proof. From (95) we obtain

$$\widetilde{\Re}^{\star}_{st} \in \widetilde{\Re}^{\star}_{st} \Leftrightarrow \left[ Bew_{ZFC_{st}}\left(\#\left(\widetilde{\Re}^{\star}_{st} \notin \widetilde{\Re}^{\star}_{st}\right)\right)\right] \wedge \wedge\left[ Bew_{ZFC_{st}}\left(\#\left(\widetilde{\Re}^{\star}_{st} \notin \widetilde{\Re}^{\star}_{st}\right)\right) \Rightarrow \widetilde{\Re}^{\star}_{st} \notin \widetilde{\Re}^{\star}_{st}\right]. \quad (98)$$

(a) Assume now that:

$$\widetilde{\Re}^{\star}_{st} \in \widetilde{\Re}^{\star}_{st}. \quad (99)$$

Then from (98) we obtain $\vdash Bew_{ZFC_{st}}\left(\#\left(\widetilde{\Re}^{\star}_{st} \notin \widetilde{\Re}^{\star}_{st}\right)\right)$ and $\vdash Bew_{ZFC_{st}}\left(\#\left(\widetilde{\Re}^{\star}_{st} \notin \widetilde{\Re}^{\star}_{st}\right)\right) \Rightarrow \widetilde{\Re}^{\star}_{st} \notin \widetilde{\Re}^{\star}_{st}$, therefore $\vdash \widetilde{\Re}^{\star}_{st} \notin \widetilde{\Re}^{\star}_{st}$ and so

$$\vdash_{ZFC_{st}} \widetilde{\Re}^{\star}_{st} \in \widetilde{\Re}^{\star}_{st} \Rightarrow \widetilde{\Re}^{\star}_{st} \notin \widetilde{\Re}^{\star}_{st}. \quad (100)$$

From (99)-(100) we obtain $\widetilde{\Re}^{\star}_{st} \in \widetilde{\Re}^{\star}_{st}, \widetilde{\Re}^{\star}_{st} \in \widetilde{\Re}^{\star}_{st} \Rightarrow \widetilde{\Re}^{\star}_{st} \notin \widetilde{\Re}^{\star}_{st} \vdash \widetilde{\Re}^{\star}_{st} \notin \widetilde{\Re}^{\star}_{st}$ and therefore $\vdash_{ZFC_{st}} \left(\widetilde{\Re}^{\star}_{st} \in \widetilde{\Re}^{\star}_{st}\right) \wedge \left(\widetilde{\Re}^{\star}_{st} \notin \widetilde{\Re}^{\star}_{st}\right)$.

(b) Assume now that

$$\left[ Bew_{ZFC_{st}}\left(\#\left(\widetilde{\Re}^{\star}_{st} \notin \widetilde{\Re}^{\star}_{st}\right)\right)\right] \wedge \wedge\left[ Bew_{ZFC_{st}}\left(\#\left(\widetilde{\Re}^{\star}_{st} \notin \widetilde{\Re}^{\star}_{st}\right)\right) \Rightarrow \widetilde{\Re}^{\star}_{st} \notin \widetilde{\Re}^{\star}_{st}\right]. \quad (101)$$

Then from (101) we obtain $\vdash \widetilde{\Re}_2^{\star Hs} \notin \widetilde{\Re}_2^{\star Hs}$. From (101) and (98) we obtain $\vdash_{\overline{ZFC}_2^{Hs}} \widetilde{\Re}_2^{\star Hs} \in \widetilde{\Re}_2^{\star Hs}$, so $\vdash_{\overline{ZFC}_2^{Hs}} \widetilde{\Re}_2^{\star Hs} \notin \widetilde{\Re}_2^{\star Hs}, \widetilde{\Re}_2^{\star Hs} \in \widetilde{\Re}_2^{\star Hs}$ which immediately gives us a contradiction $\vdash_{\overline{ZFC}_2^{Hs}} \left(\widetilde{\Re}_2^{\star Hs} \in \widetilde{\Re}_2^{\star Hs}\right) \wedge \left(\widetilde{\Re}_2^{\star Hs} \notin \widetilde{\Re}_2^{\star Hs}\right)$.

5mm . **3.3. Derivation of the inconsistent provably definable set in** $ZFC_{Nst}$ 2mm

**Designation 3.3.1**. (i) Let $\overline{PA}$ be a first order theory which contain usual postulates of Peano arithmetic [8]cite: Cohen66 and recursive defining equations for every primitive recursive function
as desired.

(ii) Let $M^{ZFC}_{Nst}$ be a nonstandard model of $ZFC$ and let $M^{\overline{PA}}_{st}$ be a standard model of $\overline{PA}$. We assume now that $M^{\overline{PA}}_{st} \subset M^{ZFC}_{Nst}$ and denote such nonstandard model of $ZFC$ by $M^{ZFC}_{Nst}[\overline{PA}]$.

(iii) Let $ZFC_{Nst}$ be the theory $ZFC_{Nst} = ZFC + M^{ZFC}_{Nst}[\overline{PA}]$.

(iv) Let $\Im_{Nst}$ be the set of the all sets of $M^{ZFC}_{st}[\overline{PA}]$ provably definable in $ZFC_{Nst}$, and let $\Re_{Nst} = \{X \in \Im_{Nst} : \Box_{Nst}(X \notin X)\}$ where $\Box_{Nst}A$ means 'sentence $A$ derivable in $ZFC_{Nst}$', or some appropriate modification thereof. We replace now (45) by formula

$$\forall Y\{Y \in \Im_{Nst} \leftrightarrow \Box_{Nst}[\exists\Psi(\cdot)\exists!X[\Psi(X) \wedge Y = X]]\}, \quad (102)$$

and we replace (46) by formula

$$\forall X\left[X \in \Re_{Nst} \leftrightarrow (X \in \Im_{Nst}) \wedge \Box_{Nst}\left(X \notin X\right)\right]. \quad (103)$$

Assume that $ZFC_{Nst} \vdash \Re_{Nst} \in \Im_{Nst}$. Then, we have that: $\Re_{Nst} \in \Re_{Nst}$ iff $\Box_{Nst}(\Re_{Nst} \notin \Re_{Nst})$, which immediately gives us $\Re_{Nst} \in \Re_{Nst}$ iff $\Re_{Nst} \notin \Re_{Nst}$. But this is a contradiction, i.e., $ZFC_{Nst} \vdash (\Re_{Nst} \in \Re_{Nst}) \wedge (\Re_{Nst} \notin \Re_{Nst})$. We choose now $\Box_{Nst}A$ in the following form

$$\Box_{Nst}A \triangleq Bew_{ZFC_{Nst}}(\#A) \wedge [Bew_{ZFC_{Nst}}(\#A) \Rightarrow A]. \quad (104)$$

Here $Bew_{ZFC_{Nst}}(\#A)$ is a canonycal Gödel formula which says to us that there exists proof in $ZFC_{Nst}$ of the formula $A$ with Gödel number $\#A \in M^{PA}_{st}$.

**Remark 3.3.1**. Notice that definition (104) holds as definition of predicate really asserting provability in $ZFC_{Nst}$.

**Designation 3.3.2**. (i) Let $g_{ZFC_{Nst}}(u)$ be a Gödel number of given an expression $u$ of $ZFC_{Nst}$.

(ii) Let $\mathbf{Fr}_{Nst}(y,v)$ be the relation: $y$ is the Gödel number of a wff of $ZFC_{Nst}$ that contains free occurrences of the variable with Gödel number $v$ [6]cite: Foukzon19, [10]cite: Rossberg.

(iii) Let $\wp_{Nst}(y,v,v_1)$ be a Gödel number of the following wff: $\exists!X[\Psi(X) \wedge Y = X]$, where $g_{ZFC_{Nst}}(\Psi(X)) = y, g_{ZFC_{Nst}}(X) = v,\ g_{ZFC_{Nst}}(Y) = v_1$.

(iv) Let $\Pr_{ZFC_{Nst}}(z)$ be a predicate asserting provability in $ZFC_{Nst}$.

**Remark 3.3.2**. Let $\mathfrak{I}_{Nst}$ be the countable collection of all sets $X$ such that $ZFC_{Nst} \vdash \exists!X\Psi(X)$, where $\Psi(X)$ is a 1-place open wff i.e.,

$$\forall Y\{Y \in \mathfrak{I}_{Nst} \Leftrightarrow ZFC_{Nst} \vdash \exists\Psi(X)\exists!X[\Psi(X) \wedge Y = X]\}. \tag{105}$$

We rewrite now (105) in the following form

$$\begin{gathered}\forall Y\{Y \in \mathfrak{I}^{\star}_{Nst} \Leftrightarrow \\ (g_{ZFC_{Nst}}(Y) = v_1) \wedge \exists y\widehat{\mathbf{Fr}}_{Nst}(y,v) \wedge (g_{ZFC_{Nst}}(X) = v) \wedge [\Pr_{ZFC_{Nst}}(\wp_{Nst}(y,v,v_1)) \wedge \\ \wedge[\Pr_{ZFC_{Nst}}(\wp_{Nst}(y,v,v_1)) \Rightarrow \exists!X[\Psi(X) \wedge Y = X]]]\}\end{gathered} \tag{106}$$

**Designation 3.3.3**. Let $\wp_{Nst}(z)$ be a Gödel number of the following wff: $Z \notin Z$, where $g_{ZFC_{Nst}}(Z) = z$.

**Remark 3.3.3**. Let $\mathfrak{R}_{Nst}$ above by formula (103), i.e.,

$$\forall Z\big[Z \in \mathfrak{R}_{Nst} \leftrightarrow (Z \in \mathfrak{I}_{Nst}) \wedge \square_{Nst}\big(Z \notin Z\big)\big]. \tag{107}$$

We rewrite now (107) in the following form

$$\begin{gathered}\forall Z[Z \in \mathfrak{R}^{\star}_{Nst} \leftrightarrow (Z \in \mathfrak{I}^{\star}_{Nst}) \wedge g_{ZFC_{Nst}}(Z) = z \wedge \Pr_{ZFC_{Nst}}(\wp_{Nst}(z))] \wedge \\ \wedge\big[\Pr_{ZFC_{Nst}}(\wp_{Nst}(z)) \Rightarrow Z \notin Z\big].\end{gathered} \tag{108}$$

**Theorem 3.3.1**. $ZFC_{Nst} \vdash \mathfrak{R}^{\star}_{Nst} \in \mathfrak{R}^{\star}_{Nst} \wedge \mathfrak{R}^{\star}_{Nst} \notin \mathfrak{R}^{\star}_{Nst}$.

## 5mm . 3.4. Generalized Tarski's undefinability lemma 2mm

**Remark 3.4.1**. Remind that: (i) if **Th** is a theory, let $T_{\mathbf{Th}}$ be the set of Godel numbers of theorems of **Th**,[10], (ii) the property $x \in T_{\mathbf{Th}}$ is said to be is expressible in **Th** by wff $\mathbf{True}(x_1)$ if the following properties are satisfied [10]**cite: Rossberg**:

(a) if $n \in T_{\mathbf{Th}}$ then $\mathbf{Th} \vdash \mathbf{True}(\bar{n})$, (b) if $n \notin T_{\mathbf{Th}}$ then $\mathbf{Th} \vdash \neg\mathbf{True}(\bar{n})$.

**Remark 3.4.2**. Notice it follows from (a)$\wedge$(b) that $\neg[(\mathbf{Th} \nvdash \mathbf{True}(\bar{n})) \wedge (\mathbf{Th} \nvdash \neg\mathbf{True}(\bar{n}))]$.

**Theorem 3.4.1**. (Tarski's undefinability Lemma) [10]**cite: Rossberg**. Let **Th** be a consistent theory with equality in the language $\mathscr{L}$ in which the diagonal function $D$ is representable and let $g_{\mathbf{Th}}(u)$ be a Gödel number of given an expression $u$ of **Th.** Then the property $x \in T_{\mathbf{Th}}$ is not expressible in **Th**.

**Proof**. By the diagonalization lemma applied to $\neg\mathbf{True}(x_1)$ there is a sentence $\mathcal{F}$ such that: (c) $\mathbf{Th} \vdash \mathcal{F} \Leftrightarrow \neg\mathbf{True}(\bar{q})$, where $q$ is the Godel number of $\mathcal{F}$, i.e. $g_{\mathbf{Th}}(\mathcal{F}) = q$.

**Case 1**. Suppose that $\mathbf{Th} \vdash \mathcal{F}$, then $q \in T_{\mathbf{Th}}$. By (a), $\mathbf{Th} \vdash \mathbf{True}(\bar{q})$. But, from $\mathbf{Th} \vdash \mathcal{F}$and (c), by biconditional elimination, one obtains $\mathbf{Th} \vdash \neg\mathbf{True}(\bar{q})$.Hence **Th** is inconsistent, ontradicting our hypothesis.

**Case 2**. Suppose that $\mathbf{Th} \nvdash \mathcal{F}$, then $q \notin T_{\mathbf{Th}}$. By (b), $\mathbf{Th} \vdash \neg\mathbf{True}(\bar{q})$. Hence, by (c) and biconditional elimination, $\mathbf{Th} \vdash \mathcal{F}$.Thus, in either case a contradiction is reached.

**Definition 3.4.1**. If **Th** is a theory, let $T_{\mathbf{Th}}$ be the set of Godel numbers of theorems of **Th** and let $g_{\mathbf{Th}}(u)$ be a Gödel number of given an expression $u$ of **Th.** The property $x \in T_{\mathbf{Th}}$ is said to be is a strongly expressible in **Th** by wff $\mathbf{True}^*(x_1)$ if the following properties are satisfied:

(a) if $n \in T_{\mathbf{Th}}$ then $\mathbf{Th} \vdash \mathbf{True}^*(\bar{n}) \wedge (\mathbf{True}^*(\bar{n}) \Rightarrow g^{-1}_{\mathbf{Th}}(n))$,

(b) if $n \notin T_{\mathbf{Th}}$ then $\mathbf{Th} \vdash \neg\mathbf{True}^*(\bar{n})$.

**Theorem 3.4.2**. (Generalized Tarski's undefinability Lemma).Let **Th** be a consistent theory

with equality in the language $\mathcal{L}$ in which the diagonal function $D$ is representable and let $g_{\mathbf{Th}}(u)$ be a Gödel number of given an expression $u$ of **Th**. Then the property $x \in T_{\mathbf{Th}}$ is not strongly expressible in **Th**.

**Proof**. By the diagonalization lemma applied to $\neg\mathbf{True}^*(x_1)$ there is a sentence $\mathcal{F}^*$ such that: (c) $\mathbf{Th} \vdash \mathcal{F}^* \iff \neg\mathbf{True}^*(\overline{q})$, where $q$ is the Godel number of $\mathcal{F}^*$, i.e. $g_{\mathbf{Th}}(\mathcal{F}^*) = q$.

**Case 1**. Suppose that $\mathbf{Th} \vdash \mathcal{F}^*$, then $q \in T_{\mathbf{Th}}$. By (a), $\mathbf{Th} \vdash \mathbf{True}^*(\overline{q})$. But, from $\mathbf{Th} \vdash \mathcal{F}^*$and (c), by biconditional elimination, one obtains $\mathbf{Th} \vdash \neg\mathbf{True}^*(\overline{q})$. Hence **Th** is inconsistent, contradicting our hypothesis.

**Case 2**. Suppose that $\mathbf{Th} \nvdash \mathcal{F}^*$, then $q \notin T_{\mathbf{Th}}$. By (b), $\mathbf{Th} \vdash \neg\mathbf{True}^*(\overline{q})$. Hence, by (c) and biconditional elimination, $\mathbf{Th} \vdash \mathcal{F}^*$. Thus, in either case a contradiction is reached.

**Remark 3.4.3**. Notice that Tarski's undefinability theorem cannot blocking the biconditionals

$$\begin{gathered}\Re \in \Re \iff \Re \notin \Re\,, \Re_{st} \in \Re_{st} \iff \Re_{st} \notin \Re_{st}\,,\\ \Re_{Nst} \in \Re_{Nst} \iff \Re_{Nst} \notin \Re_{Nst},\end{gathered} \tag{109}$$

see subsection 2.2.

## 5mm . 3.5. Generalized Tarski's undefinability theorem 2mm

**Remark 3.5.1**. (**I**) Let $\mathbf{Th}_1^\#$ be the theory $\mathbf{Th}_1^\# \triangleq \overline{ZFC}_2^{Hs}$.

In addition under assumption $\widetilde{Con}(\mathbf{Th}_1^\#)$, we establish a countable sequence of the consistent extensions of the theory $\mathbf{Th}_1^\#$ such that:

(i)$\mathbf{Th}_1^\# \subsetneqq \ldots \subsetneqq \mathbf{Th}_i^\# \subsetneqq \mathbf{Th}_{i+1}^\# \subsetneqq \ldots \mathbf{Th}_\infty^\#$, where

(ii) $\mathbf{Th}_{i+1}^\#$is a finite consistent extension of $\mathbf{Th}_i^\#$,

(iii) $\mathbf{Th}_\infty^\# = \bigcup_{i\in\mathbb{N}} \mathbf{Th}_i^\#$,

(iv) $\mathbf{Th}_\infty^\#$ proves the all sentences of $\mathbf{Th}_1^\#$, which is valid in $M$, i.e.,$M \models A \Rightarrow \mathbf{Th}_\infty^\# \vdash A$, see see Subsection 4.1, Proposition 4.1.1.

(**II**) Let $\mathbf{Th}_{1,st}^\#$ be $\mathbf{Th}_{1,st}^\# \triangleq ZFC_{st}$.

In addition under assumption $\widetilde{Con}(\mathbf{Th}_{1,st}^\#)$, we establish a countable sequence of the consistent extensions of the theory $\mathbf{Th}_1^\#$ such that:

(i) $\mathbf{Th}_{1,st}^\# \subsetneqq \ldots \subsetneqq \mathbf{Th}_{i,st}^\# \subsetneqq \mathbf{Th}_{i+1,st}^\# \subsetneqq \ldots \mathbf{Th}_{\infty,st}^\#$, where

(ii) $\mathbf{Th}_{i+1,st}^\#$ is a finite consistent extension of $\mathbf{Th}_{i,st}^\#$,

(iii) $\mathbf{Th}_{\infty,st}^\# = \bigcup_{i\in\mathbb{N}} \mathbf{Th}_{i,st}^\#$,

(iv) $\mathbf{Th}_{\infty,st}^\#$ proves the all sentences of $\mathbf{Th}_{1,st}^\#$, which valid in $M_{st}^{ZFC}$, i.e., $M_{st}^{ZFC} \models A \Rightarrow \mathbf{Th}_{\infty,st}^\# \vdash A$, see Subsection 4.1, Proposition 4.1.1.

(**III**) Let $\mathbf{Th}_{1,Nst}^\#$ be $\mathbf{Th}_{1,Nst}^\# \triangleq ZFC_{Nst}$.

In addition under assumption $\widetilde{Con}(\mathbf{Th}_{1,Nst}^\#)$, we establish a countable sequence of the consistent extensions of the theory $\mathbf{Th}_1^\#$ such that:

(i)$\mathbf{Th}_{1,Nst}^\# \subsetneqq \ldots \subsetneqq \mathbf{Th}_{i,Nst}^\# \subsetneqq \mathbf{Th}_{i+1,st}^\# \subsetneqq \ldots \mathbf{Th}_{\infty,Nst}^\#$, where

(ii) $\mathbf{Th}_{i+1,Nst}^\#$ is a finite consistent extension of $\mathbf{Th}_{i,Nst}^\#$,

(iii) $\mathbf{Th}_{\infty,st}^\# = \bigcup_{i\in\mathbb{N}} \mathbf{Th}_{i,st}^\#$

(iv) $\mathbf{Th}_{\infty,st}^\#$ proves the all sentences of $\mathbf{Th}_{1,st}^\#$, which valid in $M_{Nst}^{ZFC}[PA]$, i.e., $M_{Nst}^{ZFC}[PA] \models A \Rightarrow \mathbf{Th}_{\infty,Nst}^\# \vdash A$, see see Subsection 4.1, Proposition 4.1.1.

**Remark 3.5.2**. (**I**) Let $\Im_i, i = 1,2,\ldots$ be the set of the all sets of $M$ provably definable in $\mathbf{Th}_i^\#$,

$$\forall Y\{Y \in \Im_i \leftrightarrow \Box_i \exists \Psi(\cdot) \exists! X[\Psi(X) \wedge Y = X]\}. \tag{110}$$

and let $\Re_i = \{x \in \Im_i : \Box_i(x \notin x)\}$ where $\Box_i A$ means sentence $A$ derivable in $\mathbf{Th}_i^\#$.Then we have that $\Re_i \in \Re_i$ iff $\Box_i(\Re_i \notin \Re_i)$, which immediately gives us $\Re_i \in \Re_i$ iff $\Re_i \notin \Re_i$. We choose now $\Box_i A, i = 1,2,\ldots$ in the following form

$$\Box_i A \triangleq Bew_i(\#A) \wedge [Bew_i(\#A) \Rightarrow A]. \tag{111}$$

Here $Bew_i(\#A), i = 1,2,\ldots$ is a canonycal Gödel formulae which says to us that there exists proof in $\mathbf{Th}_i^\#, i = 1,2,\ldots$of the formula $A$ with Gödel number $\#A$.

(**II**) Let $\Im_{i,st}, i = 1,2,\ldots$ be the set of the all sets of $M_{st}^{ZFC}$ provably definable in $\mathbf{Th}_{i,st}^{\#}$,

$$\forall Y\{Y \in \Im_{i,st} \leftrightarrow \Box_{i,st}\exists\Psi(\cdot)\exists!X[\Psi(X) \wedge Y = X]\}. \tag{112}$$

and let $\Re_{i,st} = \{x \in \Im_{i,st} : \Box_{i,st}(x \notin x)\}$ where $\Box_{i,st}A$ means sentence $A$ derivable in $\mathbf{Th}_{i,st}^{\#}$.
Then we have that $\Re_{i,st} \in \Re_{i,st}$ iff $\Box_{i,st}(\Re_{i,st} \notin \Re_{i,st})$, which immediately gives us $\Re_{i,st} \in \Re_{i,st}$ iff $\Re_{i,st} \notin \Re_{i,st}$. We choose now $\Box_{i,st}A, i = 1,2,\ldots$ in the following form

$$\Box_{i,st}A \triangleq Bew_{i,st}(\#A) \wedge [Bew_{i,st}(\#A) \Rightarrow A]. \tag{113}$$

Here $Bew_{i,st}(\#A), i = 1,2,\ldots$ is a canonycal Gödel formulae which says to us that there exists proof in $\mathbf{Th}_{i,st}^{\#}, i = 1,2,\ldots$ of the formula $A$ with Gödel number $\#A$.

(**III**) Let $\Im_{i,Nst}, i = 1,2,\ldots$ be the set of the all sets of $M_{Nst}^{ZFC}[PA]$ provably definable in $\mathbf{Th}_{i,Nst}^{\#}$,

$$\forall Y\{Y \in \Im_{i,Nst} \leftrightarrow \Box_{i,Nst}\exists\Psi(\cdot)\exists!X[\Psi(X) \wedge Y = X]\}. \tag{114}$$

and let $\Re_{i,Nst} = \{x \in \Im_{i,Nst} : \Box_{i,Nst}(x \notin x)\}$ where $\Box_{i,Nst}A$ means sentence $A$ derivable in $\mathbf{Th}_{i,Nst}^{\#}$. Then we have that $\Re_{i,Nst} \in \Re_{i,Nst}$ iff $\Box_{i,Nst}(\Re_{i,Nst} \notin \Re_{i,Nst})$, which immediately gives us $\Re_{i,Nst} \in \Re_{i,Nst}$ iff $\Re_{i,Nst} \notin \Re_{i,Nst}$.

We choose now $\Box_{i,Nst}A, i = 1,2,\ldots$ in the following form

$$\Box_{i,Nst}A \triangleq Bew_{i,Nst}(\#A) \wedge [Bew_{i,Nst}(\#A) \Rightarrow A]. \tag{115}$$

Here $Bew_{i,Nst}(\#A), i = 1,2,\ldots$ is a canonical Gödel formulae which says to us that there exists proof in $\mathbf{Th}_{i,Nst}^{\#}, i = 1,2,\ldots$ of the formula $A$ with Gödel number $\#A$.

**Remark 3.5.3** Notice that definitions (111), (113) and (115) hold as definitions of predicates really asserting provability in $\mathbf{Th}_{i}^{\#}, \mathbf{Th}_{i,st}^{\#}$ and $\mathbf{Th}_{i,Nst}^{\#}, i = 1,2,\ldots$ correspondingly.

**Remark 3.5.4**. Of course the all theories $\mathbf{Th}_{i}^{\#}, \mathbf{Th}_{i,st}^{\#}, \mathbf{Th}_{i,Nst}^{\#}, i = 1,2,\ldots$ are inconsistent, see subsection 4.1.

**Remark 3.5.5**. (**I**) Let $\Im_{\infty}$ be the set of the all sets of $M$ provably definable in $\mathbf{Th}_{\infty}^{\#}$,

$$\forall Y\{Y \in \Im_{\infty} \leftrightarrow \Box_{\infty}\exists\Psi(\cdot)\exists!X[\Psi(X) \wedge Y = X]\}. \tag{116}$$

and let $\Re_{\infty} = \{x \in \Im_{\infty} : \Box_{\infty}(x \notin x)\}$, where $\Box_{\infty}A$ means 'sentence $A$ derivable in $\mathbf{Th}_{\infty}^{\#}$. Then, we have that $\Re_{\infty} \in \Re_{\infty}$ iff $\Box_{\infty}(\Re_{\infty} \notin \Re_{\infty})$, which immediately gives us $\Re_{\infty} \in \Re_{\infty}$ iff $\Re_{\infty} \notin \Re_{\infty}$. We choose now $\Box_{\infty}A, i = 1,2,\ldots$ in the following form

$$\Box_{\infty}A \triangleq \exists i[Bew_{i}(\#A) \wedge [Bew_{i}(\#A) \Rightarrow A]]. \tag{117}$$

(**II**) Let $\Im_{\infty,st}$ be the set of the all sets of $M_{st}^{ZFC}$ provably definable in $\mathbf{Th}_{\infty,st}^{\#}$,

$$\forall Y\{Y \in \Im_{\infty,st} \leftrightarrow \Box_{\infty,st}\exists\Psi(\cdot)\exists!X[\Psi(X) \wedge Y = X]\}. \tag{118}$$

and let $\Re_{\infty,st}$ be the set $\Re_{\infty,st} = \{x \in \Im_{\infty,st} : \Box_{\infty,st}(x \notin x)\}$, where $\Box_{\infty,st}A$ means 'sentence $A$ derivable in $\mathbf{Th}_{\infty,st}^{\#}$. Then, we have that $\Re_{\infty,st} \in \Re_{\infty,st}$ iff $\Box_{\infty,st}(\Re_{\infty,st} \notin \Re_{\infty,st})$, which immediately gives us $\Re_{\infty,st} \in \Re_{\infty,st}$ iff $\Re_{\infty,st} \notin \Re_{\infty,st}$. We choose now $\Box_{\infty,st}A, i = 1,2,\ldots$ in the following form

$$\Box_{\infty,st}A \triangleq \exists i[Bew_{i,st}(\#A) \wedge [Bew_{i,st}(\#A) \Rightarrow A]]. \tag{119}$$

(**III**) Let $\Im_{\infty,Nst}$ be the set of the all sets of $M_{Nst}^{ZFC}[PA]$ provably definable in $\mathbf{Th}_{\infty,Nst}^{\#}$,

$$\forall Y\{Y \in \Im_{\infty,Nst} \leftrightarrow \Box_{\infty,Nst}\exists\Psi(\cdot)\exists!X[\Psi(X) \wedge Y = X]\}. \tag{120}$$

and let $\Re_{\infty,Nst}$ be the set $\Re_{\infty,Nst} = \{x \in \Im_{\infty,Nst} : \Box_{\infty,Nst}(x \notin x)\}$ where $\Box_{\infty,Nst}A$ means 'sentence $A$ derivable in $\mathbf{Th}_{\infty,Nst}^{\#}$. Then, we have that $\Re_{\infty,Nst} \in \Re_{\infty,Nst}$ iff $\Box_{\infty,Nst}(\Re_{\infty,Nst} \notin \Re_{\infty,Nst})$, which immediately gives us $\Re_{\infty,Nst} \in \Re_{\infty,Nst}$ iff $\Re_{\infty,Nst} \notin \Re_{\infty,Nst}$. We choose now $\Box_{\infty,Nst}A, i = 1,2,\ldots$ in the following form

$$\Box_{\infty,Nst}A \triangleq \exists i[Bew_{i,Nst}(\#A) \wedge [Bew_{i,Nst}(\#A) \Rightarrow A]]. \tag{121}$$

**Remark 3.5.6**. Notice that definitions (117), (119) and (121) hold as definitions of a predicate really asserting provability in $\mathbf{Th}_{\infty}^{\#}, \mathbf{Th}_{\infty,st}^{\#}$ and $\mathbf{Th}_{\infty,Nst}^{\#}$ correspondingly.

**Remark 3.5.7**. Of course all the theories $\mathbf{Th}_{\infty}^{\#}, \mathbf{Th}_{\infty,st}^{\#}$ and $\mathbf{Th}_{\infty,Nst}^{\#}$ are inconsistent, see subsection 4.1.

**Remark 3.5.8**. Notice that under naive consideration the set $\Im_{\infty}$ and $\Re_{\infty}$ can be defined directly using a truth predicate, which of couse is not available in the language of $ZFC_{2}^{Hs}$ (but iff $ZFC_2^{Hs}$ is consistent) by well-known Tarski's undefinability theorem [10]cite: Rossberg.

**Theorem 3.5.1. Tarski's undefinability theorem**: **(I)** Let $\mathbf{Th}_{\mathcal{L}}$ be first order theory with formal language $\mathcal{L}$, which includes negation and has a Gödel numbering $g(\circ)$ such that for every $\mathcal{L}$-formula $A(x)$ there is a formula $B$ such that $B \leftrightarrow A(g(B))$ holds. Assume that $\mathbf{Th}_{\mathcal{L}}$ has a standard model $M_{st}^{\mathbf{Th}_{\mathcal{L}}}$ and $Con(\mathbf{Th}_{\mathcal{L},st})$ where

$$\mathbf{Th}_{\mathcal{L},st} \triangleq \mathbf{Th}_{\mathcal{L}} + \exists M_{st}^{\mathbf{Th}_{\mathcal{L}}}. \quad (122)$$

Let $T^*$ be the set of Gödel numbers of $\mathcal{L}$-sentences true in $M_{st}^{\mathbf{Th}_{\mathcal{L}}}$. Then there is no $\mathcal{L}$-formula $\mathbf{True}(n)$ (truth predicate) which defines $T^*$. That is, there is no $\mathcal{L}$-formula $\mathbf{True}(n)$ such that for every $\mathcal{L}$-formula $A$,

$$\mathbf{True}(g(A)) \iff A \quad (123)$$

holds.

**(II)** Let $\mathbf{Th}_{\mathcal{L}}^{Hs}$ be second order theory with Henkin semantics and formal language $\mathcal{L}$, which includes negation and has a Gödel numbering $g(\circ)$ such that for every $\mathcal{L}$-formula $A(x)$ there is a formula $B$ such that $B \leftrightarrow A(g(B))$ holds.

Assume that $\mathbf{Th}_{\mathcal{L}}^{Hs}$ has a standard model $M_{st}^{\mathbf{Th}_{\mathcal{L}}^{Hs}}$ and $Con(\mathbf{Th}_{\mathcal{L},st}^{Hs})$, where

$$\mathbf{Th}_{\mathcal{L},st}^{Hs} \triangleq \mathbf{Th}_{\mathcal{L}}^{Hs} + \exists M_{st}^{\mathbf{Th}_{\mathcal{L}}^{Hs}} \quad (124)$$

Let $T^*$ be the set of Gödel numbers of the all $\mathcal{L}$-sentences true in $M$. Then there is no $\mathcal{L}$-formula $\mathbf{True}(n)$ (truth predicate) which defines $T^*$. That is, there is no $\mathcal{L}$-formula $\mathbf{True}(n)$ such that for every $\mathcal{L}$-formula $A$,

$$\mathbf{True}(g(A)) \iff A \quad (125)$$

holds.

**Remark 3.5.9**. Notice that the proof of Tarski's undefinability theorem in this form is again by simple reductio ad absurdum. Suppose that an $\mathcal{L}$- formula True(n) defines $T^*$. In particular, if $A$ is a sentence of $\mathbf{Th}_{\mathcal{L}}$ then $\mathbf{True}(g(A))$ holds in $\mathbb{N}$ iff $A$ is true in $M_{st}^{\mathbf{Th}_{\mathcal{L}}}$. Hence for all $A$, the Tarski $T$-sentence $\mathbf{True}(g(A)) \iff A$ is true in $M_{st}^{\mathbf{Th}_{\mathcal{L}}}$. But the diagonal lemma yields a counterexample to this equivalence, by giving a "Liar" sentence $S$ such that $S \iff \neg\mathbf{True}(g(S))$ holds in $M_{st}^{\mathbf{Th}_{\mathcal{L}}}$. Thus no $\mathcal{L}$-formula $\mathbf{True}(n)$ can define $T^*$.

**Remark 3.5.10**. Notice that the formal machinery of this proof is wholly elementary except for the diagonalization that the diagonal lemma requires. The proof of the diagonal lemma is likewise surprisingly simple; for example, it does not invoke recursive functions in any way. The proof does assume that every $\mathcal{L}$-formula has a Gödel number, but the specifics of a coding method are not required.

**Remark 3.5.11**. The undefinability theorem does not prevent truth in one consistent theory from being defined in a stronger theory. For example, the set of (codes for) formulas of first-order Peano arithmetic that are true in $\mathbb{N}$ is definable by a formula in second order arithmetic. Similarly, the set of true formulas of the standard model of second order arithmetic (or n-th order arithmetic for any n) can be defined by a formula in first-order $ZFC$.

**Remark1.3.5.12**. Notice that Tarski's undefinability theorem cannot blocking the biconditionals

$$\begin{gathered}\mathfrak{R}_i \in \mathfrak{R}_i \iff \mathfrak{R}_i \notin \mathfrak{R}_i, i \in \mathbb{N}, \\ \mathfrak{R}_\infty \in \mathfrak{R}_\infty \iff \mathfrak{R}_\infty \notin \mathfrak{R}_\infty, \text{etc.,}\end{gathered} \quad (126)$$

see Remark 3.5.14 below.

**Remark 3.5.13.(I)** We define again the set $\Im_\infty$ but now by using generalized truth predicate $\mathbf{True}_\infty^\#(g(A), A)$ such that

$$\mathbf{True}^{\#}_{\infty}(g(A),A) \Leftrightarrow \exists i[Bew_i(\#A) \wedge [Bew_i(\#A) \Rightarrow A]] \Leftrightarrow$$
$$\mathbf{True}_{\infty}(g(A)) \wedge [\mathbf{True}_{\infty}(g(A)) \Rightarrow A] \Leftrightarrow A, \quad (127)$$
$$\mathbf{True}_{\infty}(g(A)) \Leftrightarrow \exists i Bew_i(\#A).$$

holds.

(**II**) We define the set $\Im_{\infty,st}$ using generalized truth predicate $\mathbf{True}^{\#}_{\infty,st}(g(A),A)$ such that

$$\mathbf{True}^{\#}_{\infty,st}(g(A),A) \Leftrightarrow \exists i[Bew_{i,st}(\#A) \wedge [Bew_{i,st}(\#A) \Rightarrow A]] \Leftrightarrow$$
$$\mathbf{True}_{\infty,st}(g(A)) \wedge [\mathbf{True}_{\infty,st}(g(A)) \Rightarrow A] \Leftrightarrow A, \quad (128)$$
$$\mathbf{True}_{\infty,st}(g(A)) \Leftrightarrow \exists i Bew_{i,st}(\#A)$$

holds. Thus in contrast with naive definition of the sets $\Im_{\infty}$ and $\Re_{\infty}$ there is no any problem which arises from Tarski's undefinability theorem.

(**III**) We define a set $\Im_{\infty,Nst}$ using generalized truth predicate $\mathbf{True}^{\#}_{\infty,Nst}(g(A),A)$ such that

$$\mathbf{True}^{\#}_{\infty,Nst}(g(A),A) \Leftrightarrow \exists i[Bew_{i,Nst}(\#A) \wedge [Bew_{i,Nst}(\#A) \Rightarrow A]] \Leftrightarrow$$
$$\mathbf{True}_{\infty,Nst}(g(A)) \wedge [\mathbf{True}_{\infty,Nst}(g(A)) \Rightarrow A] \Leftrightarrow A, \quad (129)$$
$$\mathbf{True}_{\infty,Nst}(g(A)) \Leftrightarrow \exists i Bew_{i,Nst}(\#A).$$

holds. Thus in contrast with naive definition of the sets $\Im_{\infty,Nst}$ and $\Re_{\infty,Nst}$ there is no any problem which arises from Tarski's undefinability theorem.

**Remark 3.5.14**. In order to prove that set theory $ZFC_2^{Hs} + \exists M^{ZFC_2^{Hs}}$ is inconsistent without any reference to the set $\Im_{\infty}$, notice that by the properties of the extension $\mathbf{Th}^{\#}_{\infty}$ it follows that definition given by formula (127) is correct, i.e., for every $ZFC_2^{Hs}$-formula $\Phi$ such that $M^{ZFC_2^{Hs}} \models \Phi$ the following equivalence $\Phi \Leftrightarrow \mathbf{True}_{\infty}(g(\Phi),\Phi)$ holds.

**Theorem 3.5.2**. (**Generalized Tarski's undefinability theorem**) (see subsection 4.2, Proposition 4.2.1). Let $\mathbf{Th}_{\mathcal{L}}$ be a first order theory or the second order theory with Henkin semantics and with formal language $\mathcal{L}$, which includes negation and has a Gödel encoding $g(\bullet)$ such that for every $\mathcal{L}$-formula $A(x)$ there is a formula $B$ such that the equivalence $B \Leftrightarrow A(g(B))$ holds. Assume that $\mathbf{Th}_{\mathcal{L}}$ has a standard Model $M_{st}^{\mathbf{Th}}$. Then there is no $\mathcal{L}$-formula $\mathbf{True}(n), n \in \mathbb{N}$, such that for every $\mathcal{L}$-formula $A$ such that $M \models A$, the following equivalence holds

$$A \Leftrightarrow \mathbf{True}(g(A),A). \quad (130)$$

**Theorem 3.5.3**. (i) Set theory $\mathbf{Th}^{\#}_1 = ZFC_2^{Hs} + \exists M^{ZFC_2^{Hs}}$ is inconsistent;

(ii) Set theory $\mathbf{Th}^{\#}_{1,st} = ZFC + \exists M^{ZFC}_{st}$ is inconsistent; (iii) Set theory $\mathbf{Th}^{\#}_{1,Nst} = ZFC + \exists M^{ZFC}_{Nst}$ is inconsistent; (see subsection 4.2, Proposition 4.2.2).

**Proof**. (i) Notice that by the properties of the extension $\mathbf{Th}^{\#}_{\infty}$ of the theory $ZFC_2^{Hs} + \exists M^{ZFC_2^{Hs}} = \mathbf{T}h^{\#}_1$ it follows that

$$M^{ZFC_2^{Hs}} \models \Phi \Rightarrow \mathbf{Th}^{\#}_{\infty} \vdash \Phi. \quad (131)$$

Therefore formula (127) gives generalized "truth predicate" for the set theory $\mathbf{Th}^{\#}_{\mathbf{1}}$. By Theorem 3.5.2 one obtains a contradiction.

(ii) Notice that by the properties of the extension $\mathbf{Th}^{\#}_{\infty,Nst}$ of the theory$ZFC + \exists M^{ZFC}_{st} = \mathbf{Th}^{\#}_{1,st}$ it follows that

$$M^{ZFC}_{st} \models \Phi \Rightarrow \mathbf{Th}^{\#}_{\infty,st} \vdash \Phi. \quad (132)$$

Therefore formula (128) gives generalized "truth predicate" for the set theory $\mathbf{Th}^{\#}_{1,st}$. By Theorem 3.5.2 one obtains a contradiction.

(iii) Notice that by the properties of the extension $\mathbf{Th}^{\#}_{\infty,Nst}$ of the theory $ZFC + \exists M^{ZFC}_{Nst} = \mathbf{Th}^{\#}_{1,st}$ it follows that

$$M^{ZFC}_{Nst} \models \Phi \Rightarrow \mathbf{Th}^{\#}_{\infty,Nst} \vdash \Phi. \quad (133)$$

Therefore (129) gives generalized "truth predicate" for the set theory $\mathbf{Th}^{\#}_{1,Nst}$. By Theorem 3.5.2 one obtains a contradiction.

## 5mm . 3.6. Avoiding the contradictions from set theory $\overline{ZFC}_2^{Hs}$, $ZFC_{st}$ and set theory $ZFC_{Nst}$ using Quinean approach 2mm

In order to avoid difficulties mentioned above we use well known Quinean approach [19]cite: Quine37.

### 5mm . 3.6.1 Quinean set theory $NF$ 2mm

Remind that the primitive predicates of Russellian unramified typed set theory (TST), a streamlined version of the theory of types, are equality = and membership $\in$. TST has a linear hierarchy of types: type 0 consists of individuals otherwise undescribed. For each (meta-) natural number $n$, type $n+1$ objects are sets of type $n$ objects; sets of type $n$ have members of type $n-1$. Objects connected by identity must have the same type. The following two atomic formulas succinctly describe the typing rules: $x^n = y^n$ and $x^n \in y^{n+1}$.

The axioms of TST are:

**Extensionality**: sets of the same (positive) type with the same members are equal;

**Axiom schema of comprehension**:

If $\Phi(x^n)$ is a formula, then the set $\{x^n \mid \Phi(x^n)\}^{n+1}$ exists i.e., given any formula $\Phi(x^n)$, the formula

$$\exists A^{n+1} \forall x^n [x^n \in A^{n+1} \leftrightarrow \Phi(x^n)] \quad (134)$$

is an axiom where $A^{n+1}$ represents the set $\{x^n \mid \Phi(x^n)\}^{n+1}$ and is not free in $\Phi(x^n)$.

Quinean set theory. (New Foundations) seeks to eliminate the need for such superscripts.

New Foundations has a universal set, so it is a non-well founded set theory. That is to say, it is a logical theory that allows infinite descending chains of membership such as $x_n \in x_{n-1} \in \ldots x_3 \in x_2 \in x_1$. It avoids Russell's paradox by only allowing stratifiable formulas in the axiom of comprehension. For instance $x \in y$ is a stratifiable formula, but $x \in x$ is not (for details of how this works see below).

**Definition 3.6.1**. In New Foundations (*NF*) and related set theories, a formula $\Phi$ in the language of first-order logic with equality and membership is said to be stratified iff there is a function σ which sends each variable appearing in $\Phi$ [considered as an item of syntax] to a natural number (this works equally well if all integers are used) in such a way that any atomic formula $x \in y$ appearing in $\Phi$ satisfies $\sigma(x) + 1 = \sigma(y)$ and any atomic formula $x = y$ appearing in $\Phi$ satisfies $\sigma(x) = \sigma(y)$.

## Quinean Set Theory $NF$

**Axioms and stratification are**:

The well-formed formulas of New Foundations (*NF*) are the same as the well-formed formulas of TST, but with the type annotations erased. The axioms of *NF* are [19]cite: Quine37:

**Extensionality**: Two objects with the same elements are the same object.

A comprehension schema: All instances of TST Comprehension but with type indices dropped (and without introducing new identifications between variables).

By convention, NF's Comprehension schema is stated using the concept of stratified formula and making no direct reference to types. Comprehension then becomes.

**Stratified Axiom schema of comprehension**:

$\{x \mid \Phi^s\}$ exists for each stratified formula $\Phi^s$.

Even the indirect reference to types implicit in the notion of stratification can be eliminated. Theodore Hailperin showed in 1944 that Comprehension is equivalent to a finite conjunction of its instances, so that *NF* can be finitely axiomatized without any reference to the notion of type [20]cite: Hailperin44. Comprehension may seem to run afoul of problems similar to those in naive set theory, but this is not the case. For example, the existence of the impossible Russell class $\{x \mid x \notin x\}$ is not an axiom of $NF$, because $x \notin x$ cannot be stratified.

### 5mm . 3.6.2 Set theory $\overline{ZFC}_2^{Hs}, ZFC_{st}$ and set theory $ZFC_{Nst}$ with stratified axiom schema of replacement 2mm

The stratified axiom schema of replacement asserts that the image of a set under any function definable by stratified formula of the theory $ZFC_{st}$ will also fall inside a set.

**Stratified Axiom schema of replacement**:

Let $\Phi^s(x,y,w_1,w_2,\ldots,w_n)$ be any stratified formula in the language of $ZFC_{st}$ whose free variables are among $x,y,A,w_1,w_2,\ldots,w_n$, so that in particular $B$ is not free in $\Phi^s$. Then

$$\begin{gathered}\forall A\forall w_1\forall w_2\ldots\forall w_n[\forall x(x\in A\Rightarrow\exists!y\Phi^s(x,y,w_1,w_2,\ldots,w_n))\Rightarrow\\ \Rightarrow\exists B\forall x(x\in A\Rightarrow\exists y(y\in B\wedge\Phi^s(x,y,w_1,w_2,\ldots,w_n)))],\end{gathered}\tag{135}$$

i.e., if the relation $\Phi^s(x,y,\ldots)$ represents a definable function $f,A$ represents its domain, and $f(x)$ is a set for every $x\in A$, then the range of $f$ is a subset of some set $B$.

**Stratified Axiom schema of separation**:

Let $\Phi^s(x,w_1,w_2,\ldots,w_n)$ be any stratified formula in the language of $ZFC_{st}$ whose free variables are among $x,A,w_1,w_2,\ldots,w_n$, so that in particular $B$ is not free in $\Phi^s$. Then

$$\forall w_1\forall w_2\ldots\forall w_n\forall A\exists B\forall x[x\in B\Leftrightarrow(x\in A\wedge\Phi^s(x,w_1,w_2,\ldots,w_n))],\tag{136}$$

**Remark 3.6.1**. Notice that the stratified axiom schema of separation follows from the stratified axiom schema of replacement together with the axiom of empty set.

**Remark 3.6.2**. Notice that the stratified axiom schema of replacement (separation) obviously violeted any contradictions (82), (126), etc. mentioned above. The existence of the countable Russell sets $\Re_2^{*Hs}, \Re_{st}^*$ and $\Re_{Nst}^*$ is impossible, because $x\notin x$ cannot be stratified.

## 6mm . 4.Generalized Löbs Theorem 3mm

### 5mm . 4.1. Generalized Löbs Theorem. Second-Order theories with Henkin semantics 2mm

**Remark 4.1.1**. In this section we use second-order arithmetic $Z_2^{Hs}$ with Henkin semantics. Notice that any standard model $M_{st}^{Z_2^{Hs}}$ of second-order arithmetic $Z_2^{Hs}$ consisting of a set $\mathbb{N}$ of uusual natural numbers (which forms the range of individual variables) together with a constant 0 (an element of $\mathbb{N}$), a function S from $\mathbb{N}$ to $\mathbb{N}$, two binary operations + and $\cdot$ on $\mathbb{N}$, a binary relation < on $\mathbb{N}$, and a collection $D\subseteq 2^{\mathbb{N}}$ of subsets of $\mathbb{N}$, which is the range of the set variables. Omitting $D$ produces a model of the first order Peano arithmetic.

When $D = 2^{\mathbb{N}}$ is the full powerset of $\mathbb{N}$, the model $M_{st}^{Z_2}$ is called a full model. The use of full second-order semantics is equivalent to limiting the models of second-order arithmetic to the full models. In fact, the axioms of second-order arithmetic $Z_2^{fss}$ have only one full model. This follows from the fact that the axioms of Peano arithmetic with the ssecond-order induction axiom have inly one model under second-order semantics, see section 3.

Let **Th** be some fixed, but unspecified, consistent formal theory. For later convenience, we assume that the encoding is done in some fixed formal second order theory **S** and that **Th** contains **S**. We assume throughout this paper that formal second order theory **S** has an $\omega$-model $M_\omega^{\mathbf{S}}$. The sense in which **S** is contained in **Th** is better exemplified than explained: if **S** is a formal system of a second order arithmetic $Z_2^{Hs}$ and **Th** is, say, $ZFC_2^{Hs}$, then **Th** contains **S** in the sense that there is a well-known embedding, or interpretation, of **S** in **Th**. Since encoding is to take place in $M_\omega^{\mathbf{S}}$, it will have to have a large supply of constants and closed terms to be used as codes.

(e.g. in formal arithmetic, one has $\bar{0}, \bar{1}, \ldots$ .) **S** will also have certain function symbols to be described shortly. To each formula, $\Phi$, of the language of **Th** is assigned a closed term, $[\Phi]^c$, called the code of $\Phi$ [19]cite: Quine37. We note that if $\Phi(x)$ is a formula with free variable $x$, then $[\Phi(x)]^c$ is a closed term encoding the formula $\Phi(x)$ with $x$ viewed as a syntactic object and not as a parameter. Corresponding to the logical connectives and quantifiers are the function symbols, $neg(\bullet)$, $imp(\bullet)$, etc., such that for all first order formulae

$\Phi, \Psi$ : $\mathbf{S} \vdash neg([\Phi]^c) = [\neg\Phi]^c$, $\mathbf{S} \vdash imp([\Phi]^c, [\Psi]^c) = [\Phi \to \Psi]^c$ etc. Of particular importance is the substitution operator, represented by the function symbol $sub(\bullet, \bullet)$. For formulae $\Phi(x)$, terms $t$ with codes $[t]^c$ :

$$\mathbf{S} \vdash sub([\Phi(x)]^c, [t]^c) = [\Phi(t)]^c. \tag{137}$$

It well known that one can also encode derivations and have a binary relation $\mathbf{Prov_{Th}}(x,y)$ (read "$x$ proves $y$ " or "$x$ is a proof of $y$") such that for closed $t_1, t_2$ : $\mathbf{S} \vdash \mathbf{Prov_{Th}}(t_1, t_2)$ iff $t_1$ is the code of a derivation in **Th** of the formula with code $t_2$ . It follows that

$$\mathbf{Th} \vdash \Phi \text{ iff } \mathbf{S} \vdash \mathbf{Prov_{Th}}(t, [\Phi]^c) \tag{138}$$

for some closed term $t$. Thus we can define

$$\mathbf{Pr_{Th}}(y) \leftrightarrow \exists x \mathbf{Prov_{Th}}(x,y), \tag{139}$$

and therefore we obtain a predicate asserting provability.

**Remark 4.1.2**. (**I**) We note that it is not always the case that:

$$\mathbf{Th} \vdash \Phi \text{ iff } \mathbf{S} \vdash \mathbf{Pr_{Th}}([\Phi]^c), \tag{140}$$

unless **S** is fairly sound, e.g. this is the case when **S** and **Th** replaced by $\mathbf{S}_\omega = \mathbf{S} \restriction M_\omega^{\mathbf{Th}}$ and $\mathbf{Th}_\omega = \mathbf{Th} \restriction M_\omega^{\mathbf{Th}}$ correspondingly (see Designation 4.1.1 below).

(**II**) Notice that it is always the case that:

$$\mathbf{Th}_\omega \vdash \Phi_\omega \text{ iff } \mathbf{S}_\omega \vdash \mathbf{Pr_{Th_\omega}}([\Phi_\omega]^c), \tag{141}$$

i.e. that is the case when predicate $\mathbf{Pr_{Th_\omega}}(y), y \in M_\omega^{\mathbf{Th}}$ :

$$\mathbf{Pr_{Th_\omega}}(y) \leftrightarrow \exists x (x \in M_\omega^{\mathbf{Th}}) \mathbf{Prov_{Th_\omega}}(x,y) \tag{142}$$

really asserting provability.

It well known that the above encoding can be carried out in such a way that the following important conditions $\mathbf{D}1, \mathbf{D}2$ and $\mathbf{D}3$ are meet for all sentences:

$$\begin{gathered}\mathbf{D}1.\, \mathbf{Th} \vdash \Phi \text{ implies } \mathbf{S} \vdash \ \mathbf{Th} \vdash \mathbf{Pr_{Th}}([\Phi]^c), \\ \mathbf{D}2.\, \mathbf{S} \vdash \mathbf{Pr_{Th}}([\Phi]^c) \to \mathbf{Pr_{Th}}([\mathbf{Pr_{Th}}([\Phi]^c)]^c), \\ \mathbf{D}3.\, \mathbf{S} \vdash \mathbf{Pr_{Th}}([\Phi]^c) \wedge \mathbf{Pr_{Th}}([\Phi \to \Psi]^c) \to \mathbf{Pr_{Th}}([\Psi]^c).\end{gathered} \tag{143}$$

Conditions $\mathbf{D}1, \mathbf{D}2$ and $\mathbf{D}3$ are called the Derivability Conditions.

**Remark 4.1.3.**From (141)-(142) it follows that

$$\begin{gathered}\mathbf{D}4.\, \mathbf{Th}_\omega \vdash \Phi \text{ iff } \mathbf{S}_\omega \vdash \mathbf{Pr_{Th_\omega}}([\Phi_\omega]^c), \\ \mathbf{D}5.\, \mathbf{S}_\omega \vdash \mathbf{Pr_{Th_\omega}}([\Phi_\omega]^c) \leftrightarrow \mathbf{Pr_{Th_\omega}}([\mathbf{Pr_{Th_\omega}}([\Phi_\omega]^c)]^c), \\ \mathbf{D}6.\, \mathbf{S}_\omega \vdash \mathbf{Pr_{Th_\omega}}([\Phi_\omega]^c) \wedge \mathbf{Pr_{Th_\omega}}([\Phi_\omega \to \Psi]^c) \to \mathbf{Pr_{Th_\omega}}([\Psi_\omega]^c).\end{gathered} \tag{144}$$

Conditions $\mathbf{D}4, \mathbf{D}5$ and $\mathbf{D}6$ are called a Strong Derivability Conditions.

**Definition 4.1.1**. Let $\Phi$ be well formed formula (wff) of **Th**. Then wff $\Phi$ is called **Th**-sentence iff it has no free variables.

**Designation 4.1.1** (i) Assume that a theory **Th** has an $\omega$-model $M_\omega^{\mathbf{Th}}$ and $\Phi$ is a **Th**-sentence, then: $\Phi_{M_\omega^{\mathbf{Th}}} \triangleq \Phi \restriction M_\omega^{\mathbf{Th}}$ (we will write $\Phi_\omega$ instead $\Phi_{M_\omega^{\mathbf{Th}}}$) is a **Th**-sentence $\Phi$ with all quantifiers relativized to $\omega$-model $M_\omega^{\mathbf{Th}}$ [11] and $\mathbf{Th}_\omega \triangleq \mathbf{Th} \restriction M_\omega^{\mathbf{Th}}$ is a theory **Th** relativized to model $M_\omega^{\mathbf{Th}}$, i.e., any $\mathbf{Th}_\omega$-sentence has the form $\Phi_\omega$ for some $\mathbf{Th}$-sentence $\Phi$.

(ii) Assume that a theory $\mathbf{Th}$ has a standard model $M_{st}^{\mathbf{Th}}$ and $\Phi$ is a $\mathbf{Th}$-sentence, then:

$\Phi_{M_{st}^{\mathbf{Th}}} \triangleq \Phi \upharpoonright M_{st}^{\mathbf{Th}}$ (we will write $\Phi_{st}$ instead $\Phi_{M_{st}^{\mathbf{Th}}}$) is a $\mathbf{Th}$-sentence with all quantifiers relativized to a standard model $M_{st}^{\mathbf{Th}}$, and $\mathbf{Th}_{st} \triangleq \mathbf{Th} \upharpoonright M_{st}^{\mathbf{Th}}$ is a theory $\mathbf{Th}$ relativized to model $M_{st}^{\mathbf{Th}}$, i.e., any $\mathbf{Th}_{st}$-sentence has a form $\Phi_{st}$ for some $\mathbf{Th}$-sentence $\Phi$.

(iii) Assume that a theory $\mathbf{Th}$ has a non-standard model $M_{Nst}^{\mathbf{Th}}$ and $\Phi$ is a $\mathbf{Th}$-sentence, then:

$\Phi_{M_{Nst}^{\mathbf{Th}}} \triangleq \Phi \upharpoonright M_{Nst}^{\mathbf{Th}}$ (we will write $\Phi_{Nst}$ instead $\Phi_{M_{Nst}^{\mathbf{Th}}}$) is a $\mathbf{Th}$-sentence with all quantifiers relativized to non-standard model $M_{Nst}^{\mathbf{Th}}$, and $\mathbf{Th}_{Nst} \triangleq \mathbf{Th} \upharpoonright M_{Nst}^{\mathbf{Th}}$ is a theory $\mathbf{Th}$ relativized to model $M_{Nst}^{\mathbf{Th}}$, i.e., any $\mathbf{Th}_{Nst}$-sentence has a form $\Phi_{Nst}$ for some $\mathbf{Th}$-sentence $\Phi$.

(iv) Assume that a theory $\mathbf{Th}$ has a model $M = M^{\mathbf{Th}}$ and $\Phi$ is a $\mathbf{Th}$-sentence, then: $\Phi_{M^{\mathbf{Th}}}$ is a $\mathbf{Th}$-sentence with all quantifiers relativized to model $M^{\mathbf{Th}}$, and $\mathbf{Th}_M$ is a theory $\mathbf{Th}$ relativized to model $M^{\mathbf{Th}}$, i.e. any $\mathbf{Th}_M$-sentence has a form $\Phi_M$ for some $\mathbf{Th}$-sentence $\Phi$.

**Designation 4.1.2**. (i) Assume that a theory $\mathbf{Th}$ with a language $\mathcal{L}$ has an $\omega$-model $M_\omega^{\mathbf{Th}}$ and there exists $\mathbf{Th}$-sentence $S_{\mathcal{L}}$ such that: (a) $S_{\mathcal{L}}$ expressible by language $\mathcal{L}$ and

(b) $S_{\mathcal{L}}$ asserts that $\mathbf{Th}$ has a model $M_\omega^{\mathbf{Th}}$; we denote such $\mathbf{Th}$-sentence $S_{\mathcal{L}}$ by $Con(\mathbf{Th}; M_\omega^{\mathbf{Th}})$.

(ii) Assume that a theory $\mathbf{Th}$ with a language $\mathcal{L}$ has a non-standard model $M_{Nst}^{\mathbf{Th}}$ and there exists $\mathbf{Th}$-sentence $S_{\mathcal{L}}$ such that: (a) $S_{\mathcal{L}}$ expressible by language $\mathcal{L}$ and

(b) $S_{\mathcal{L}}$ asserts that $\mathbf{Th}$ has a non-standard model $M_{Nst}^{\mathbf{Th}}$; we denote such $\mathbf{Th}$-sentence $S_{\mathcal{L}}$ by $Con(\mathbf{Th}; M_{Nst}^{\mathbf{Th}})$.

(iii) Assume that a theory $\mathbf{Th}$ with a language $\mathcal{L}$ has an model $M^{\mathbf{Th}}$ and there exists $\mathbf{Th}$-sentence $S_{\mathcal{L}}$ such that: (a) $S_{\mathcal{L}}$ expressible by language $\mathcal{L}$ and

(b) $S_{\mathcal{L}}$ asserts that $\mathbf{Th}$ has a model $M^{\mathbf{Th}}$; we denote such $\mathbf{Th}$-sentence $S_{\mathcal{L}}$ by $Con(\mathbf{Th}; M^{\mathbf{Th}})$.

**Remark 4.1.4**. We emphasize that: (i) it is well known that there exist a $ZFC$-sentence $Con(ZFC; M^{ZFC})$ [8]cite: Cohen66,

(ii) obviously there exists a $ZFC_2^{Hs}$-sentence $Con\left(ZFC_2^{Hs}; M^{ZFC_2^{Hs}}\right)$ and there exists a $Z_2^{Hs}$-sentence $Con\left(Z_2^{Hs}; M^{Z_2^{Hs}}\right)$.

**Designation 4.1.3**. Assume that $Con(\mathbf{Th}; M^{\mathbf{Th}})$. Let $\widetilde{Con}(\mathbf{Th}; M^{\mathbf{Th}})$ be the formula:

$$\begin{gathered}
\widetilde{Con}(\mathbf{Th}; M^{\mathbf{Th}}) \triangleq \\
\forall t_1(t_1 \in M_\omega^{\mathbf{Th}}) \forall t_1'(t_1' \in M_\omega^{\mathbf{Th}}) \forall t_2(t_2 \in M_\omega^{\mathbf{Th}}) \forall t_2'(t_2' \in M_\omega^{\mathbf{Th}}) \\
\neg[\mathbf{Prov_{Th}}(t_1, [\Phi]^c) \wedge \mathbf{Prov_{Th}}(t_2, neg([\Phi]^c))], \\
\text{where} \\
t_1' = [\Phi]^c, t_2' = neg([\Phi]^c) \\
\text{or} \\
\widetilde{Con}(\mathbf{Th}; M_\omega^{\mathbf{Th}}) \triangleq \\
\forall \Phi \forall t_1(t_1 \in M_\omega^{\mathbf{Th}}) \forall t_2(t_2 \in M_\omega^{\mathbf{Th}}) \neg[\mathbf{Prov_{Th}}(t_1, [\Phi]^c) \wedge \mathbf{Prov_{Th}}(t_2, neg([\Phi]^c))]
\end{gathered} \tag{145}$$

and where $t_1, t_1', t_2, t_2'$ is a closed term.

**Lemma 4.1.1**. (**I**) Assume that: (i) a theory $\mathbf{Th}$ is recursively axiomatizable.

(ii) $Con(\mathbf{Th}; M^{\mathbf{Th}})$,

(iii) $M^{\mathbf{Th}} \models \widetilde{Con}(\mathbf{Th}; M^{\mathbf{Th}})$ and

(iv) $\mathbf{Th} \vdash \mathbf{Pr_{Th}}([\Phi]^c)$, where $\Phi$ is a closed formula.

Then $\mathbf{Th} \nvdash \mathbf{Pr_{Th}}([\neg\Phi]^c)$.

(**II**) Assume that: (i) a theory $\mathbf{Th}$ is recursively axiomatizable.

(ii) $Con(\mathbf{Th}; M_\omega^{\mathbf{Th}})$

(iii) $M_\omega^{\mathbf{Th}} \models \widetilde{Con}(\mathbf{Th}; M^{\mathbf{Th}})$ and

(iv) $\mathbf{Th}_\omega \vdash \mathbf{Pr}_{\mathbf{Th}_\omega}([\Phi_\omega]^c),$ where $\Phi_\omega$ is a closed formula.
Then $\mathbf{Th}_\omega \nvdash \mathbf{Pr}_{\mathbf{Th}_\omega}([\neg\Phi_\omega]^c)$.

**Proof**. **(I)** Let $\widetilde{Con}_{\mathbf{Th}}(\Phi; M^{\mathbf{Th}})$ be the formula :

$$\widetilde{Con}(\Phi; M^{\mathbf{Th}}) \triangleq$$
$$\forall t_1(t_1 \in M^{\mathbf{Th}})\forall t_2(t_2 \in M^{\mathbf{Th}})\neg[\mathbf{Prov}_{\mathbf{Th}}(t_1,[\Phi]^c) \wedge \mathbf{Prov}_{\mathbf{Th}}(t_2, neg([\Phi]^c))],$$
$$\text{i.e.} \qquad (146)$$
$$\forall t_1(t_1 \in M^{\mathbf{Th}})\forall t_2(t_2 \in M^{\mathbf{Th}})\neg[\mathbf{Prov}_{\mathbf{Th}}(t_1,[\Phi]^c) \wedge \mathbf{Prov}_{\mathbf{Th}}(t_2, neg([\Phi]^c))]\leftrightarrow$$
$$\{\neg\exists t_1(t_1 \in M^{\mathbf{Th}})\neg\exists t_2(t_2 \in M^{\mathbf{Th}})[\mathbf{Prov}_{\mathbf{Th}}(t_1,[\Phi]^c) \wedge \mathbf{Prov}_{\mathbf{Th}}(t_2, neg([\Phi]^c))]\}.$$

where $t_1, t_2$ is a closed term. From (i)-(ii) it follows that theory $\mathbf{Th} + \widetilde{Con}(\mathbf{Th}; M^{\mathbf{Th}})$ is consistent. We note that $\mathbf{Th} + \widetilde{Con}(\mathbf{Th}; M^{\mathbf{Th}}) \vdash \widetilde{Con}_{\mathbf{Th}}(\Phi; M^{\mathbf{Th}})$ for any closed $\Phi$. Suppose that $\mathbf{Th} \vdash \mathbf{Pr}_{\mathbf{Th}}([\neg\Phi]^c)$, then (iii) gives

$$\mathbf{Th} \vdash \mathbf{Pr}_{\mathbf{Th}}([\Phi]^c) \wedge \mathbf{Pr}_{\mathbf{Th}}([\neg\Phi]^c). \qquad (147)$$

From (139) and (147) we obtain

$$\exists t_1 \exists t_2[\mathbf{Prov}_{\mathbf{Th}}(t_1,[\Phi]^c) \wedge \mathbf{Prov}_{\mathbf{Th}}(t_2, neg([\Phi]^c))]. \qquad (148)$$

But the formula (146) contradicts the formula (148). Therefore $\mathbf{Th} \nvdash \mathbf{Pr}_{\mathbf{Th}}([\neg\Phi]^c)$.

**Remark 4.1.5**. In additional note that under the following conditions:
(i) a theory **Th** is recursively axiomatizable,
(ii) $Con(\mathbf{Th}; M_{st}^{\mathbf{Th}})$, and
(iii) $M_{st}^{\mathbf{Th}} \models \widetilde{Con}(\mathbf{Th}; M_{st}^{\mathbf{Th}})$ predicate $\mathbf{Pr}_{\mathbf{Th}}([\Psi]^c)$ really asserts provability, one obtains

$$\mathbf{Th} \vdash \Phi \wedge \neg\Phi \qquad (149)$$

and therefore by reductio ad absurdum again one obtains $\mathbf{Th} \nvdash \mathbf{Pr}_{\mathbf{Th}}([\neg\Phi]^c)$.
**(II)** Let $\widetilde{Con}_{\mathbf{Th}}(\Phi; M_\omega^{\mathbf{Th}})$ be the formula:

$$\widetilde{Con}_{\mathbf{Th}}(\Phi; M_\omega^{\mathbf{Th}}) \triangleq$$
$$\forall t_1(t_1 \in M_\omega^{\mathbf{Th}})\forall t_2(t_2 \in M_\omega^{\mathbf{Th}})\neg[\mathbf{Prov}_{\mathbf{Th}}(t_1,[\Phi]^c) \wedge \mathbf{Prov}_{\mathbf{Th}}(t_2, neg([\Phi]^c))],$$
$$\text{i.e.} \qquad (150)$$
$$\forall t_1(t_1 \in M_\omega^{\mathbf{Th}})\forall t_2(t_2 \in M_\omega^{\mathbf{Th}})\neg[\mathbf{Prov}_{\mathbf{Th}}(t_1,[\Phi]^c) \wedge \mathbf{Prov}_{\mathbf{Th}}(t_2, neg([\Phi]^c))]\leftrightarrow$$
$$\{\neg\exists t_1(t_1 \in M_\omega^{\mathbf{Th}})\neg\exists t_2(t_2 \in M_\omega^{\mathbf{Th}})[\mathbf{Prov}_{\mathbf{Th}}(t_1,[\Phi]^c) \wedge \mathbf{Prov}_{\mathbf{Th}}(t_2, neg([\Phi]^c))]\}.$$

This case is trivial because formula $\mathbf{Pr}_{\mathbf{Th}_\omega}([\neg\Phi]^c)$ by the Strong Derivability Condition **D**4, see formulae (144), really asserts provability of the $\mathbf{Th}_\omega$-sentence $\neg\Phi_\omega$. But this is a contradiction.

**Lemma 4.1.2**. **(I)** Assume that: (i) a theory **Th** is recursively axiomatizable.
(ii) $Con(\mathbf{Th}; M^{\mathbf{Th}})$,
(iii) $M^{\mathbf{Th}} \models \widetilde{Con}(\mathbf{Th})$ and
(iv) $\mathbf{Th} \vdash \mathbf{Pr}_{\mathbf{Th}}([\neg\Phi]^c)$, where $\Phi$ is a closed formula. Then $\mathbf{Th} \nvdash \mathbf{Pr}_{\mathbf{Th}}([\Phi]^c)$,
**(II)** Assume that: (i) a theory **Th** is recursively axiomatizable.
(ii) $Con(\mathbf{Th}; M_\omega^{\mathbf{Th}})$
(iii) $M_\omega^{\mathbf{Th}} \models \widetilde{Con}(\mathbf{Th})$ and
(iv) $\mathbf{Th}_\omega \vdash \mathbf{Pr}_{\mathbf{Th}_\omega}([\neg\Phi_\omega]^c)$,
where $\Phi_\omega$ is a closed formula. Then $\mathbf{Th}_\omega \nvdash \mathbf{Pr}_{\mathbf{Th}_\omega}([\Phi_\omega]^c)$.
**Proof**. Similarly as Lemma 4.1.1 above.
**Example 4.1.1**. (i) Let $\mathbf{Th} = \mathbf{PA}$ be Peano arithmetic and $\Phi \Leftrightarrow 0 = 1$.

Assume that: (i) $Con(\mathbf{PA}; M^{\mathbf{PA}})$

(ii) $M^{\mathbf{PA}} \models \widetilde{Con}(\mathbf{PA}; M^{\mathbf{PA}})$ where $M^{\mathbf{PA}}$ is a model of $\mathbf{PA}$.

Then obviously $\mathbf{PA} \vdash \mathbf{Pr_{PA}}(0 \neq 1)$ since $\mathbf{PA} \vdash 0 \neq 1$ and therefore by Lemma 4.1.1 $\mathbf{PA} \nvdash \mathbf{Pr_{PA}}(0 = 1)$.

(ii) Let $Con(\mathbf{PA}; M^{\mathbf{PA}})$ , $M^{\mathbf{PA}} \models \neg\widetilde{Con}(\mathbf{PA}; M^{\mathbf{PA}})$ and let $\mathbf{PA}^{\clubsuit}$ be a theory $\mathbf{PA}^{\clubsuit} = \mathbf{PA} + \neg\widetilde{Con}(\mathbf{PA}; M^{\mathbf{PA}})$ and $\Phi \Leftrightarrow 0 = 1$. Then obviously

$$\mathbf{PA}^{\clubsuit} \vdash [\mathbf{Pr_{PA}}(0 \neq 1)] \wedge [\mathbf{Pr_{PA}}(0 = 1)]. \tag{151}$$

and therefore

$$\mathbf{PA}^{\clubsuit} \vdash \mathbf{Pr_{PA}}(0 \neq 1), \tag{152}$$

and

$$\mathbf{PA}^{\clubsuit} \vdash \mathbf{Pr_{PA}}(0 = 1). \tag{153}$$

However by Löbs theorem

$$\mathbf{PA}^{\clubsuit} \nvdash 0 = 1. \tag{154}$$

(iii) Let $Con(\mathbf{PA}^{\clubsuit}; M^{\mathbf{PA}^{\clubsuit}})$ , $M^{\mathbf{PA}^{\clubsuit}} \models \widetilde{Con}(\mathbf{PA}^{\clubsuit}; M^{\mathbf{PA}^{\clubsuit}})$ and $\Phi \Leftrightarrow 0 = 1$. Then obviously $\mathbf{PA}^{\clubsuit} \vdash \mathbf{Pr}_{\mathbf{PA}^{\clubsuit}}(0 \neq 1)$ since $\mathbf{PA}^{\clubsuit} \vdash 0 \neq 1$ and therefore by Lemma 4.1.1 we obtain. $\mathbf{PA}^{\clubsuit} \nvdash \mathbf{Pr}_{\mathbf{PA}^{\clubsuit}}(0 = 1)$.

**Remark 4.1.6**. Notice that there is no standard model of $\mathbf{PA}^{\clubsuit}$.

**Assumption 4.1.1**. Let **Th** be a second order theory with Henkin semantics. We assume now that:

(i) the language of **Th** consists of:

numerals $\overline{0}, \overline{1}, \ldots$

countable set of the numerical variables: $\{v_0, v_1, \ldots\}$

countable set $\mathcal{F}_1$ of the first order variables, i.e.

a set of variables: $\mathcal{F}_1 = \{x, y, z, X, Y, Z, \mathfrak{I}, \mathfrak{R}, \ldots\}$

countable set $\mathcal{F}_2$ of the first order variables, i.e.

a set of variables: $\mathcal{F}_2 = \{f_0^n, R_0^n, f_1^n, R_1^n, \ldots\}$

countable set of the $n$-ary function symbols: $f_0^n, f_1^n, \ldots$

countable set of the $n$-ary relation symbols: $R_0^n, R_1^n, \ldots$

connectives: $\neg, \to$

quantifier: $\forall$.

(ii) A theory **Th** is recursively axiomatizable.

(iii) **Th** contains $ZFC_2^{Hs}$ or $ZFC$ or $NF$ and $\widetilde{Con}(\mathbf{Th}; M^{\mathbf{Th}})$ is expressible in **Th** by a single statement of **Th**;

(iv) **Th** has an $\omega$-model $M_\omega^{\mathbf{Th}}$ and $M_\omega^{\mathbf{Th}} \models \widetilde{Con}(\mathbf{Th}; M_\omega^{\mathbf{Th}})$; or

(v) **Th** has an nonstandard model $M_{Nst}^{\mathbf{Th}} = M_{Nst}^{\mathbf{Th}}[PA] \supset M_{st}^{\mathbf{PA}}$ and $M_{Nst}^{\mathbf{Th}} \models \widetilde{Con}(\mathbf{Th}; M_{Nst}^{\mathbf{Th}})$.

**Definition 4.1.2**. A **Th**-wff $\Phi$ (well-formed formula $\Phi$) is closed, i.e. $\Phi$ is a sentence, i.e. if it has no free variables; a wff is open if it has free variables. We'll use the slang '$k$-place open wff' to mean a wff with $k$ distinct free variables.

**Definition 4.1.3**. We will say that $\mathbf{Th}_\infty^\#$ is a nice theory or a nice extension of the **Th** iff the following properties holds:

(i) $\mathbf{Th}_\infty^\#$ contains **Th**;

(ii) Let $\Phi$ be any first order closed formula of **Th**, then $\mathbf{Th} \vdash \mathbf{Pr_{Th}}([\Phi]^c)$ implies $\mathbf{Th}_\infty^\# \vdash \Phi$;

(iii) Let $\Phi_\infty$ be any first order closed formula of $\mathbf{Th}_\infty^\#$, then $M_\omega^{\mathbf{Th}} \models \Phi_\infty$ implies $\mathbf{Th}_\infty^\# \vdash \Phi_\infty$, i.e. $Con(\mathbf{Th} + \Phi_\infty; M_\omega^{\mathbf{Th}})$ implies $\mathbf{Th}_\infty^\# \vdash \Phi_\infty$.

(iv) Let $\Phi_\infty$ be any first order closed formula of $\mathbf{Th}_\infty^\#$, then formulas $Con(\mathbf{Th} + \Phi_\infty; M_\omega^{\mathbf{Th}})$and $\widetilde{Con}(\mathbf{Th}_\infty^\# + \Phi_\infty; M_\omega^{\mathbf{Th}})$ are expressible in $\mathbf{Th}_\infty^\#$.

**Definition 4.1.4**. Let $L$ be a classical propositional logic $L$. Recall that a set $\Delta$ of $L$-wff's is said to be $L$-consistent, or consistent for short, if $\Delta \nvdash \bot$ and there are other equivalent formulations of consistency: (1) $\Delta$ is consistent, (2) $\mathbf{Ded}(\Delta) := \{A \mid \Delta \vdash A\}$ is not the set of all wff's, (3) there is a formula such that $\Delta \nvdash A$, (4) there are no formula $A$ such that $\Delta \vdash A$ and $\Delta \vdash \neg A$.

**Definition 4.1.5**. We will say that, $\mathbf{Th}_{\infty}^{\#}$ is a maximally nice theory or a maximally nice extension of the $\mathbf{Th}$ iff $\mathbf{Th}_{\infty}^{\#}$ is consistent and for any consistent nice extension $\mathbf{Th}_{\infty}^{\#\prime}$ of the $\mathbf{Th}$ : $\mathbf{Ded}(\mathbf{Th}_{\infty}^{\#}) \subseteqq \mathbf{Ded}(\mathbf{Th}_{\infty}^{\#\prime})$ implies $\mathbf{Ded}(\mathbf{Th}_{\infty}^{\#}) = \mathbf{Ded}(\mathbf{Th}_{\infty}^{\#\prime})$.

**Remark 4.1.7**. We note that a theory $\mathbf{Th}_{\infty}^{\#}$ depend on model $M_{\omega}^{\mathbf{Th}}$ or $M_{Nst}^{\mathbf{Th}}$, i.e. $\mathbf{Th}_{\infty}^{\#} = \mathbf{Th}_{\infty}^{\#}[M_{\omega}^{\mathbf{Th}}]$ or $\mathbf{Th}_{\infty}^{\#} = \mathbf{Th}_{\infty}^{\#}[M_{Nst}^{\mathbf{Th}}]$ correspondingly. We will consider now the case $\mathbf{Th}_{\infty}^{\#} \triangleq \mathbf{Th}_{\infty}^{\#}[M_{\omega}^{\mathbf{Th}}]$ without loss of generality.

**Remark 4.1.8**. Notice that in order to prove the statements: (i) $\neg Con(NF_{2}^{Hs}; M_{\omega}^{\mathbf{Th}})$, (ii) $\neg Con(NF; M_{\omega}^{\mathbf{Th}})$ the following Proposition 4.1.1 is necessary.

**Proposition 4.1.1**. (**Generalized Löbs Theorem**).

(**I**) Assume that:

(i) A theory $\mathbf{Th}$ is recursively axiomatizable.

(ii) $\mathbf{Th}$ is a second order theory with Henkin semantics.

(iii) $\mathbf{Th}$ contains $ZFC_{2}^{Hs}$.

(iv) $\mathbf{Th}$ has an $\omega$-model $M_{\omega}^{\mathbf{Th}}$, and

(v) the statement $\exists M_{\omega}^{\mathbf{Th}}$ is expressible by language of $\mathbf{Th}$ as a single sentence of $\mathbf{Th}$.

(vi) $M_{\omega}^{\mathbf{Th}} \models \widetilde{Con}(\mathbf{Th}; M_{\omega}^{\mathbf{Th}})$, where predicate $\widetilde{Con}(\mathbf{Th}; M_{\omega}^{\mathbf{Th}})$ is defined by formula 4.1.9.

Then theory $\mathbf{Th}$ can be extended to a maximally consistent nice theory $\mathbf{Th}_{\infty,st}^{\#} = \mathbf{Th}_{\infty,st}^{\#}[M_{\omega}^{\mathbf{Th}}]$. Below we write for short $\mathbf{Th}_{\infty,st}^{\#} \triangleq \mathbf{Th}_{\infty}^{\#} = \mathbf{Th}_{\infty}^{\#}[M_{\omega}^{\mathbf{Th}}]$.

**Remark 4.1.9**. We emphasize that (v) is valid for $ZFC$ despite the fact that the axioms of $ZFC$ are infinite, see [8]cite: Cohen66 Chapter II, section 7, p. 78.

(**II**) Assume that:

(i) A theory $\mathbf{Th}$ is recursively axiomatizable.

(ii) $\mathbf{Th}$ is a first order theory.

(iii) $\mathbf{Th}$ contains $ZFC$.

(iv) $\mathbf{Th}$ has an $\omega$-model $M_{\omega}^{\mathbf{Th}}$ and

(v) the statement $\exists M_{\omega}^{\mathbf{Th}}$ is expressible by language of $\mathbf{Th}$ as a single sentence of $\mathbf{Th}$.

(vi) $M_{\omega}^{\mathbf{Th}} \models \widetilde{Con}(\mathbf{Th}; M_{\omega}^{\mathbf{Th}})$, where predicate $\widetilde{Con}(\mathbf{Th}; M_{\omega}^{\mathbf{Th}})$ defined by formula 4.1.9.

Then theory $\mathbf{Th}_{\omega} \triangleq \mathbf{Th} \upharpoonright M_{\omega}^{\mathbf{Th}}$ can be extended to a maximally consistent nice theory $\mathbf{Th}_{\boldsymbol{\omega}}^{\#}$.

(**III**) Assume that:

(i) A theory $\mathbf{Th}$ is recursively axiomatizable.

(ii) $\mathbf{Th}$ is a first order theory.

(iii) $\mathbf{Th}$ contains $ZFC$.

(iv) $\mathbf{Th}$ has a nonstandard model $M_{Nst}^{\mathbf{Th}} = M_{Nst}^{\mathbf{Th}}[PA]$ and

(v) the statement $\exists M_{Nst}^{\mathbf{Th}}[PA]$ is expressible by language of $\mathbf{Th}$ as a single sentence of $\mathbf{Th}$.

(vi) $M_{Nst}^{\mathbf{Th}} \models \widetilde{Con}(\mathbf{Th}; M_{Nst}^{\mathbf{Th}})$, where predicate $\widetilde{Con}(\mathbf{Th}; M_{Nst}^{\mathbf{Th}})$ defined by formula (146).

Then theory $\mathbf{Th}$ can be extended to a maximally consistent nice theory $\mathbf{Th}_{\infty,Nst}^{\#} = \mathbf{Th}_{\infty,Nst}^{\#}[M_{Nst}^{\mathbf{Th}}]$.

**Remark 4.1.10**. We emphasize that (v) is valid for $ZFC$ despite the fact that the axioms of $ZFC$ are infinite, see [8]cite: Cohen66 Ch. II, section 7, p.78.

**Proof**. (**I**) Let $\Phi_1 \ldots \Phi_i \ldots$ be an enumeration of the all first order closed wff's of the theory $\mathbf{Th}$ (this can be achieved if the set of propositional variables, etc. can be enumerated).

Define a chain $\wp = \{\mathbf{Th}_{i,st}^{\#} | i \in \mathbb{N}\}, \mathbf{Th}_{1,st}^{\#} = \mathbf{Th}$ of consistent theories inductively as follows: assume that theory $\mathbf{Th}_{i,st}^{\#}$ is defined. Notice that below we write for short $\mathbf{Th}_{i,st}^{\#} \triangleq \mathbf{Th}_{i}^{\#}$.

(i) Suppose that the following statement (155) is satisfied

$$\left[\mathbf{Th}_{i}^{\#} \nvdash \mathbf{Pr}_{\mathbf{Th}_{i}^{\#}}([\Phi_i]^c)\right] \wedge \left[\mathbf{Th}_{i}^{\#} \nvdash \mathbf{Pr}_{\mathbf{Th}_{i}^{\#}}([\neg\Phi_i]^c)\right] \wedge M_{\omega}^{\mathbf{Th}} \models \Phi_i. \tag{155}$$

Note that

$$\mathbf{Th}_i^{\#} \nvdash \mathbf{Pr}_{\mathbf{Th}_i^{\#}}([\Phi_i]^c) \Leftrightarrow \mathbf{Th}_i^{\#} \nvdash \Phi_i,$$
$$\mathbf{Th}_i^{\#} \nvdash \mathbf{Pr}_{\mathbf{Th}_i^{\#}}([\neg\Phi_i]^c) \Leftrightarrow \mathbf{Th}_i^{\#} \nvdash \neg\Phi_i, \tag{156}$$

since predicate $\mathbf{Pr}_{\mathbf{Th}_i^{\#}}([\Phi_i]^c)$ really aserts provability in $\mathbf{Th}_i^{\#}$. Then we define a theory $\mathbf{Th}_{i+1}^{\#}$ as follows

$$\mathbf{Th}_{i+1}^{\#} \triangleq \mathbf{Th}_i^{\#} \cup \{\Phi_i\}. \tag{157}$$

**Remark 4.1.11**.Note that the predicate $\mathbf{Pr}_{\mathbf{Th}_{i+1}^{\#}}([\Phi_i]^c)$ is expressible in $\mathbf{Th}_{i+1}^{\#}$ since a theory $\mathbf{Th}_{i+1}^{\#}$ is a finite extension of the recursively axiomatizable theory $\mathbf{Th}$.

We will rewrite the conditions (155)-(157) using predicate $\mathbf{Pr}_{\mathbf{Th}_{i+1}^{\#}}^{\#}(\bullet)$ symbolically as follows:

$$\begin{gathered}
\mathbf{Th}_{i+1}^{\#} \vdash \mathbf{Pr}_{\mathbf{Th}_{i+1}^{\#}}^{\#}([\Phi_i]^c), \\
\mathbf{Pr}_{\mathbf{Th}_{i+1}^{\#}}^{\#}([\Phi_i]^c) \Leftrightarrow \\
\left[\neg\mathbf{Pr}_{\mathbf{Th}_i^{\#}}([\Phi_i]^c)\right] \wedge \left[\neg\mathbf{Pr}_{\mathbf{Th}_i^{\#}}([\neg\Phi_i]^c)\right] \wedge [M_\omega^{\mathbf{Th}} \vDash \Phi_i], \\
M_\omega^{\mathbf{Th}} \vDash \Phi_i \Leftrightarrow Con(\mathbf{Th}_\mathbf{i}^{\#} + \Phi_i; M_\omega^{\mathbf{Th}}), \\
\text{i.e.} \\
\mathbf{Pr}_{\mathbf{Th}_{i+1}^{\#}}^{\#}([\Phi_i]^c) \Leftrightarrow \\
\left[\neg\mathbf{Pr}_{\mathbf{Th}_i^{\#}}([\Phi_i]^c)\right] \wedge \left[\neg\mathbf{Pr}_{\mathbf{Th}_i^{\#}}([\neg\Phi_i]^c)\right] \wedge Con(\mathbf{Th}_i + \Phi_i; M_\omega^{\mathbf{Th}}), \\
\mathbf{Pr}_{\mathbf{Th}_{i+1}^{\#}}^{\#}([\Phi_i]^c) \Leftrightarrow \left[\neg\mathbf{Pr}_{\mathbf{Th}_i^{\#}}([\Phi_i]^c)\right] \wedge \left[\neg\mathbf{Pr}_{\mathbf{Th}_i^{\#}}([\neg\Phi_i]^c)\right] \\
\mathbf{Pr}_{\mathbf{Th}_{i+1}^{\#}}([\Phi_i]^c) \Rightarrow \mathbf{Th}_{i+1}^{\#} \vdash \Phi_i, \\
\mathbf{Th}_{i+1}^{\#} \vdash \mathbf{Pr}_{\mathbf{Th}_{i+1}^{\#}}^{\#}([\Phi_i]^c) \Rightarrow \Phi_i.
\end{gathered} \tag{158}$$

(ii) Suppose that the following statement (159) is satisfied

$$\left[\mathbf{Th}_i^{\#} \nvdash \mathbf{Pr}_{\mathbf{Th}_i^{\#}}([\Phi_i]^c)\right] \wedge \left[\mathbf{Th}_i^{\#} \nvdash \mathbf{Pr}_{\mathbf{Th}_i^{\#}}([\neg\Phi_i]^c)\right] \wedge M_\omega^{\mathbf{Th}} \vDash \neg\Phi_i. \tag{159}$$

Note that

$$\mathbf{Th}_i^{\#} \nvdash \mathbf{Pr}_{\mathbf{Th}_i^{\#}}([\Phi_i]^c) \Leftrightarrow \mathbf{Th}_i^{\#} \nvdash \Phi_i,$$
$$\mathbf{Th}_i^{\#} \nvdash \mathbf{Pr}_{\mathbf{Th}_i^{\#}}([\neg\Phi_i]^c) \Leftrightarrow \mathbf{Th}_i^{\#} \nvdash \neg\Phi_i, \tag{160}$$

since predicate $\mathbf{Pr}_{\mathbf{Th}_i^{\#}}([\neg\Phi_i]^c)$ really aserts provability in $\mathbf{Th}_i^{\#}$. Then we define a theory $\mathbf{Th}_{i+1}^{\#}$ as follows

$$\mathbf{Th}_{i+1}^{\#} \triangleq \mathbf{Th}_i^{\#} \cup \{\neg\Phi_i\}. \tag{161}$$

We will rewrite the conditions (159)-(161) using predicate $\mathbf{Pr}_{\mathbf{Th}_{i+1}^{\#}}^{\#}(\bullet)$, symbolically as follows:

$$\mathbf{Th}^{\#}_{i+1} \vdash \mathbf{Pr}^{\#}_{\mathbf{Th}^{\#}_{i+1}}([\neg\Phi_i]^c),$$
$$\mathbf{Pr}^{\#}_{\mathbf{Th}^{\#}_{i+1}}([\neg\Phi_i]^c) \Leftrightarrow \neg\mathbf{Pr}_{\mathbf{Th}^{\#}_{i}}([\neg\Phi_i]^c) \wedge [M^{\mathbf{Th}}_{\omega} \models \neg\Phi_i],$$
$$M^{\mathbf{Th}}_{\omega} \models \neg\Phi_i \Leftrightarrow Con(\mathbf{Th}^{\#}_{i}+(\neg\Phi_i);M^{\mathbf{Th}}_{\omega}),$$
$$\text{i.e.}$$
$$\mathbf{Pr}^{\#}_{\mathbf{Th}^{\#}_{i+1}}([\neg\Phi_i]^c) \Leftrightarrow \neg\mathbf{Pr}_{\mathbf{Th}^{\#}_{i}}([\neg\Phi_i]^c) \wedge Con(\mathbf{Th}_{i}+(\neg\Phi_i);M^{\mathbf{Th}}_{\omega}), \tag{162}$$
$$\mathbf{Pr}^{\#}_{\mathbf{Th}^{\#}_{i+1}}([\neg\Phi_i]^c) \Leftrightarrow \mathbf{Pr}_{\mathbf{Th}^{\#}_{i+1}}([\neg\Phi_i]^c),$$
$$\mathbf{Pr}_{\mathbf{Th}^{\#}_{i+1}}([\neg\Phi_i]^c) \Rightarrow \neg\Phi_i,$$
$$\mathbf{Th}^{\#}_{i+1} \vdash \mathbf{Pr}^{\#}_{\mathbf{Th}^{\#}_{i+1}}([\neg\Phi_i]^c) \Rightarrow \neg\Phi_i.$$

(iii) Suppose that the following statement (163) is satisfied

$$\mathbf{Th}^{\#}_{i} \vdash \mathbf{Pr}_{\mathbf{Th}^{\#}_{i}}([\Phi_i]^c) \tag{163}$$

and therefore $[\mathbf{Th}^{\#}_{i} \vdash \Phi_i] \wedge [M^{\mathbf{Th}}_{\omega} \models \Phi_i]$. Then we define a theory $\mathbf{Th}^{\#}_{i+1}$ as follows

$$\mathbf{Th}^{\#}_{i+1} \triangleq \mathbf{Th}^{\#}_{i}. \tag{164}$$

**Remark 4.1.12.** Note that predicate $\mathbf{Pr}^{\#}_{\mathbf{Th}^{\#}_{i+1}}([\Phi_i]^c)$ is expressible in $\mathbf{Th}^{\#}_{i}$ because $\mathbf{Th}^{\#}_{i}$ is a finite extension of the recursive theory **Th** and $Con(\mathbf{Th}^{\#}_{i}+\Phi_i;M^{\mathbf{Th}}_{\omega}) \in \mathbf{Th}^{\#}_{i+1}$.

(iv) Suppose that the following statement (165) is satisfied

$$\mathbf{Th}^{\#}_{i} \vdash \mathbf{Pr}_{\mathbf{Th}^{\#}_{i}}([\neg\Phi_i]^c) \tag{165}$$

and therefore $[\mathbf{Th}^{\#}_{i} \nvdash \neg\Phi_i] \wedge [M^{\mathbf{Th}}_{\omega} \models \neg\Phi_i]$.

Then we define theory $\mathbf{Th}^{\#}_{i+1}$ as follows:

$$\mathbf{Th}^{\#}_{i+1} \triangleq \mathbf{Th}^{\#}_{i}. \tag{166}$$

We define now a theory $\mathbf{Th}^{\#}_{\infty}$ as follows:

$$\mathbf{Th}^{\#}_{\infty} \triangleq \bigcup_{i\in\mathbb{N}} \mathbf{Th}^{\#}_{i}. \tag{167}$$

(**1**) First, notice that each $\mathbf{Th}^{\#}_{i}$ is consistent. This is done by induction on $i$ and by Lemmas 4.1.1-4.1.2. By assumption, the case is true when $i = 1$. Now, suppose $\mathbf{Th}^{\#}_{i}$ is consistent.

Then its deductive closure $\mathbf{Ded}(\mathbf{Th}^{\#}_{i})$ is also consistent.

(**2**) If statements (155)-(157) are satisfied, i.e. $\mathbf{Th}^{\#}_{i+1} \vdash \mathbf{Pr}^{\#}_{\mathbf{Th}^{\#}_{i+1}}([\Phi_i]^c)$ and $\mathbf{Th}^{\#}_{i+1} \vdash \Phi_i$, then clearly a theory $\mathbf{Th}^{\#}_{i+1} \triangleq \mathbf{Th}^{\#}_{i} \cup \{\Phi_i\}$ is consistent since it is a subset of closure $\mathbf{Ded}(\mathbf{Th}^{\#}_{i+1})$.

(**3**) If statements (159)-(161) are satisfied, i.e. $\mathbf{Th}^{\#}_{i+1} \vdash \mathbf{Pr}^{\#}_{\mathbf{Th}^{\#}_{i+1}}([\neg\Phi_i]^c)$ and $\mathbf{Th}^{\#}_{i+1} \vdash \neg\Phi_i$, then clearly $\mathbf{Th}^{\#}_{i+1} \triangleq \mathbf{Th}^{\#}_{i} \cup \{\neg\Phi_i\}$ is consistent since it is a subset of closure $\mathbf{Ded}(\mathbf{Th}^{\#}_{i+1})$.

(**4**) If the statement (163) is satisfied, i.e. $\mathbf{Th}^{\#}_{i} \vdash \mathbf{Pr}_{\mathbf{Th}^{\#}_{i}}([\Phi_i]^c)$ then clearly $\mathbf{Th}^{\#}_{i+1} \triangleq \mathbf{Th}^{\#}_{i}$ is consistent.

(**5**) If the statement (165) is satisfied, i.e. $\mathbf{Th}^{\#}_{i} \vdash \mathbf{Pr}_{\mathbf{Th}^{\#}_{i}}([\neg\Phi_i]^c)$ then clearly $\mathbf{Th}^{\#}_{i+1} \triangleq \mathbf{Th}^{\#}_{i}$ is consistent.

(**6**) Next, notice $\mathbf{Ded}(\mathbf{Th}^{\#}_{\infty})$ is maximally consistent nice extension of the $\mathbf{Ded}(\mathbf{Th})$. $\mathbf{Ded}(\mathbf{Th}^{\#}_{\infty})$ is consistent because, by the standard Lemma 4.1.3 below, it is the union of a chain of consistent sets. To see that $\mathbf{Ded}(\mathbf{Th}^{\#}_{\infty})$ is maximal, pick any wff $\Phi$. Then $\Phi$ is some $\Phi_i$ in the enumerated list of all wff's. Therefore for any $\Phi$ such that

$\mathbf{Th}_{i} \vdash \mathbf{Pr}_{\mathbf{Th}_{i}}([\Phi]^c)$ or $\mathbf{Th}^{\#}_{i} \vdash \mathbf{Pr}_{\mathbf{Th}^{\#}_{i}}([\neg\Phi]^c)$, either $\Phi \in \mathbf{Th}^{\#}_{\infty}$ or $\neg\Phi \in \mathbf{Th}^{\#}_{\infty}$. Since

$\mathbf{Ded}(\mathbf{Th}^{\#}_{i+1}) \subseteq \mathbf{Ded}(\mathbf{Th}^{\#}_{\infty})$, we have $\Phi \in \mathbf{Ded}(\mathbf{Th}^{\#}_{\infty})$ or $\neg\Phi \in \mathbf{Ded}(\mathbf{Th}^{\#}_{\infty})$, which implies that $\mathbf{Ded}(\mathbf{Th}^{\#}_{\infty})$ is maximally consistent nice extension of the $\mathbf{Ded}(\mathbf{Th})$.

**Definition 4.1.6.** We define now predicate $\mathbf{Pr}_{\mathbf{Th}^{\#}_{\infty}}([\Phi]^c)$ really asserting provability in $\mathbf{Th}^{\#}_{\infty}$ by the following formula

$$\mathbf{Pr}_{\mathbf{Th}^{\#}_{\infty}}([\Phi]^c) \Leftrightarrow \exists i(\Phi \in \mathbf{Th}^{\#}_i)\left[\mathbf{Pr}^{\#}_{\mathbf{Th}^{\#}_i}([\Phi]^c)\right]. \quad (168)$$

**Proof**. (**II**) **and** (**III**) similarly to (**I**).

**Lemma 4.1.3**. The union of a chain $\wp = \{\mathbf{\Gamma}_i | i \in \mathbb{N}\}$ of consistent sets $\mathbf{\Gamma}_i$, ordered by $\subseteq$ is consistent.

**Definition 4.1.7**. Let $\Psi = \Psi(x)$ be one-place open **Th**-wff such that the following condition:

$$\mathbf{Th} \triangleq \mathbf{Th}^{\#}_1 \vdash \exists! x_{\Psi}[\Psi(x_{\Psi})] \quad (169)$$

is satisfied.

**Remark 4.1.13**. We rewrite now the condition (168) using only the language of the theory $\mathbf{Th}^{\#}_1$ :

$$\begin{gathered}\{\mathbf{Th}^{\#}_1 \vdash \exists! x_{\Psi}[\Psi(x_{\Psi})]\} \Leftrightarrow \mathbf{Pr}_{\mathbf{Th}^{\#}_1}([\exists! x_{\Psi}[\Psi(x_{\Psi})]]^c) \wedge \\ \wedge\left\{\mathbf{Pr}_{\mathbf{Th}^{\#}_1}([\exists! x_{\Psi}[\Psi(x_{\Psi})]]^c) \Rightarrow \exists! x_{\Psi}[\Psi(x_{\Psi})]\right\}.\end{gathered} \quad (170)$$

**Definition 4.1.8**. We will say that, a set $y$ is a $\mathbf{Th}^{\#}_1$-set if there exist one-place open wff $\Psi(x)$ such that $y = x_{\Psi}$. We will write $y[\mathbf{Th}^{\#}_1]$ iff $y$ is a $\mathbf{Th}^{\#}_1$-set.

**Remark 4.1.14**. Note that

$$\begin{gathered}y[\mathbf{Th}^{\#}_1] \Leftrightarrow \exists\Psi(y = x_{\Psi}) \wedge \mathbf{Pr}_{\mathbf{Th}^{\#}_1}([\exists! x_{\Psi}[\Psi(x_{\Psi})]]^c) \\ \left\{\mathbf{Pr}_{\mathbf{Th}^{\#}_1}([\exists! x_{\Psi}[\Psi(x_{\Psi})]]^c) \Rightarrow \exists! x_{\Psi}[\Psi(x_{\Psi})]\right\}.\end{gathered} \quad (171)$$

**Definition 4.1.9**. Let $\Im_1$ be a set such that :

$$\forall x\left[x \in \Im_1 \leftrightarrow x \text{ is a } \mathbf{Th}^{\#}_1\text{-set}\right]. \quad (172)$$

**Proposition 4.1.2**. $\Im_1$ is a $\mathbf{Th}^{\#}_1$-set.

**Proof**. Let us consider an one-place open wff $\Psi(x)$ such that condition (169) is satisfied, i.e. $\mathbf{Th}^{\#}_1 \vdash \exists! x_{\Psi}[\Psi(x_{\Psi})]$. We note that there exists countable set $\mathcal{F}_{\Psi}$ of the one-place open wff's $\mathcal{F}_{\Psi} = \{\Psi_n(x)\}_{n \in \mathbb{N}}$ such that: (i) $\Psi(x) \in \mathcal{F}_{\Psi}$ and (ii)

$$\mathbf{Th} \triangleq \mathbf{Th}^{\#}_1 \vdash \exists! x_{\Psi}[[\Psi(x_{\Psi})] \wedge \{\forall n(n \in \mathbb{N})[\Psi(x_{\Psi}) \leftrightarrow \Psi_n(x_{\Psi})]\}]$$

or in the equivalent form

$$\begin{gathered}\mathbf{Th} \triangleq \mathbf{Th}^{\#}_1 \vdash \\ \mathbf{Pr}_{\mathbf{Th}^{\#}_1}([\exists! x_{\Psi}[\Psi(x_{\Psi})]]^c) \wedge \\ \left\{\mathbf{Pr}_{\mathbf{Th}^{\#}_1}([\exists! x_{\Psi}[\Psi(x_{\Psi})]]^c) \Rightarrow \exists! x_{\Psi}[\Psi(x_{\Psi})]\right\} \wedge \\ \left[\mathbf{Pr}_{\mathbf{Th}^{\#}_1}([\forall n(n \in \mathbb{N})[\Psi(x_{\Psi}) \leftrightarrow \Psi_n(x_{\Psi})]]^c)\right] \wedge \\ \mathbf{Pr}_{\mathbf{Th}^{\#}_1}([\forall n(n \in \mathbb{N})[\Psi(x_{\Psi}) \leftrightarrow \Psi_n(x_{\Psi})]]^c) \Rightarrow \forall n(n \in \mathbb{N})[\Psi(x_{\Psi}) \leftrightarrow \Psi_n(x_{\Psi})]\end{gathered} \quad (173)$$

or in the following equivalent form

$$\mathbf{Th}^{\#}_1 \vdash \exists! x_1[[\Psi_1(x_1)] \wedge \{\forall n(n \in \mathbb{N})[\Psi_1(x_1) \leftrightarrow \Psi_{n,1}(x_1)]\}]$$

or

$$\begin{gathered}\mathbf{Th}^{\#}_1 \vdash \\ \mathbf{Pr}_{\mathbf{Th}^{\#}_1}([\exists! x_1\Psi(x_1)]^c) \wedge \\ \left\{\mathbf{Pr}_{\mathbf{Th}^{\#}_1}([\exists! x_1\Psi(x_1)]^c) \Rightarrow \exists! x_1\Psi(x_1)\right\} \wedge \\ \left[\mathbf{Pr}_{\mathbf{Th}^{\#}_1}([\forall n(n \in \mathbb{N})[\Psi(x_1) \leftrightarrow \Psi_n(x_1)]]^c)\right] \wedge \\ \mathbf{Pr}_{\mathbf{Th}^{\#}_1}([\forall n(n \in \mathbb{N})[\Psi(x_1) \leftrightarrow \Psi_n(x_1)]]^c) \Rightarrow \forall n(n \in \mathbb{N})[\Psi(x_1) \leftrightarrow \Psi_n(x_1)],\end{gathered} \quad (174)$$

where we have set $\Psi(x) = \Psi_1(x_1), \Psi_n(x_1) = \Psi_{n,1}(x_1)$ and $x_\Psi = x_1$.

We note that any set $\mathcal{F}_{\Psi_k} = \{\Psi_{n,k}(x)\}_{n\in\mathbb{N}}, k = 1,2,\ldots$ such as mentioned above, defines an unique set $x_{\Psi_k}$, i.e. $\mathcal{F}_{\Psi_{k_1}} \cap \mathcal{F}_{\Psi_{k_2}} = \varnothing$ iff $x_{\Psi_{k_1}} \neq x_{\Psi_{k_2}}$. We note that a sets $\mathcal{F}_{\Psi_k}, k = 1,2,..$ are the part of the $ZFC_2^{Hs}$ or $ZFC$, i.e. a set $\mathcal{F}_{\Psi_k}$ is a set in the sense of $ZFC_2^{Hs}$ or $ZFC$.

Note that by using Gödel numbering one can replace any set $\mathcal{F}_{\Psi_k}, k = 1,2,..$ by a set $\Theta_k = g(\mathcal{F}_{\Psi_k})$ of the corresponding Gödel numbers such that

$$\Theta_k = g(\mathcal{F}_{\Psi_k}) = \{g(\Psi_{n,k}(x_k))\}_{n\in\mathbb{N}}, k = 1,2,\ldots . \tag{175}$$

It is easy to prove that any set $\Theta_k = g(\mathcal{F}_{\Psi_k}), k = 1,2,..$ is a $\mathbf{Th}_1^\#$-set. This is done by Gödel encoding (175), by the statament (173) and by axiom schemata of separation.

Let $g_{n,k} = g(\Psi_{n,k}(x_k)), k = 1,2,..$ be a Gödel number of the wff $\Psi_{n,k}(x_k)$. Therefore $g(\mathcal{F}_k) = \{g_{n,k}\}_{n\in\mathbb{N}}$, where we have set $\mathcal{F}_k = \mathcal{F}_{\Psi_k}$, $k = 1,2,..$ and

$$\forall k_1 \forall k_2 [\{g_{n,k_1}\}_{n\in\mathbb{N}} \cap \{g_{n,k_2}\}_{n\in\mathbb{N}} = \varnothing \leftrightarrow x_{k_1} \neq x_{k_2}]. \tag{176}$$

Let $\{\{g_{n,k}\}_{n\in\mathbb{N}}\}_{k\in\mathbb{N}}$ be a family of the sets $\{g_{n,k}\}_{n\in\mathbb{N}}, k = 1,2,\ldots$.By the axiom of choice one obtains unique set $\Im_1' = \{g_k\}_{k\in\mathbb{N}}$ such that $\forall k[g_k \in \{g_{n,k}\}_{n\in\mathbb{N}}]$. Finally one obtains a set $\Im_1$ from the set $\Im_1'$ by the axiom schema of replacement.

**Proposition 4.1.3**. Any set $\Theta_k = g(\mathcal{F}_{\Psi_k}), k = 1,2,..$ is a $\mathbf{Th}_1^\#$-set.

**Proof**. We define $g_{n,k} = g(\Psi_{n,k}(x_k)) = [\Psi_{n,k}(x_k)]^c, v_k = [x_k]^c$. Therefore $g_{n,k} = g(\Psi_{n,k}(x_k)) \leftrightarrow \mathbf{Fr}(g_{n,k}, v_k)$. Let us define now predicate $\Pi(g_{n,k}, v_k)$

$$\begin{gathered}\Pi(g_{n,k}, v_k) \leftrightarrow \mathbf{Pr}_{\mathbf{Th}_1^\#}([\exists! x_k[\Psi_{1,k}(x_1)]]^c) \wedge \\ \wedge \exists! x_k(v_k = [x_k]^c)\Big[\forall n(n \in \mathbb{N})\Big[\mathbf{Pr}_{\mathbf{Th}_1^\#}([[\Psi_{1,k}(x_k)]]^c) \leftrightarrow \mathbf{Pr}_{\mathbf{Th}_1^\#}(\mathbf{Fr}(g_{n,k}, v_k))\Big]\Big].\end{gathered} \tag{177}$$

We define now a set $\Theta_k$ such that

$$\begin{cases} \Theta_k = \Theta_k' \cup \{g_k\}, \\ \forall n(n \in \mathbb{N})[g_{n,k} \in \Theta_k' \leftrightarrow \Pi(g_{n,k}, v_k)] \end{cases} \tag{178}$$

Obviously definitions (177) and (178) are equivalent.

**Definition 4.1.10**. We define now the following $\mathbf{Th}_1^\#$-set $\Re_1 \subsetneqq \Im_1$ :

$$\forall x\big[x \in \Re_1 \Leftrightarrow (x \in \Im_1) \wedge \mathbf{Pr}_{\mathbf{Th}_1^\#}([x \notin x]^c) \wedge\big]. \tag{179}$$

**Proposition 4.1.4**. (i) $\mathbf{Th}_1^\# \vdash \exists \Re_1$,(ii) $\Re_1$ is a countable $\mathbf{Th}_1^\#$-set.

**Proof**. (i) Statement $\mathbf{Th}_1^\# \vdash \exists \Re_1$ follows immediately from the statement $\exists \Im_1$ and the axiom schema of separation, (ii) follows immediately from countability of a set $\Im_1$.Notice that $\Re_1$ is nonempty countable set such that $\mathbb{N} \subset \Re_1$, because for any $n \in \mathbb{N} : \mathbf{Th}_1^\# \vdash n \notin n$.

**Proposition 4.1.5**. A set $\Re_1$ is inconsistent.

**Proof**. From formula (179) we obtain

$$\mathbf{Th}_1^\# \vdash \Re_1 \in \Re_1 \Leftrightarrow \mathbf{Pr}_{\mathbf{Th}_1^\#}([\Re_1 \notin \Re_1]^c). \tag{180}$$

From (180) we obtain

$$\mathbf{Th}_1^\# \vdash \Re_1 \in \Re_1 \Leftrightarrow \Re_1 \notin \Re_1 \tag{181}$$

and therefore

$$\mathbf{Th}_1^\# \vdash (\Re_1 \in \Re_1) \wedge (\Re_1 \notin \Re_1). \tag{182}$$

But this is a contradiction.

**Definition 4.1.11**. Let $\Psi = \Psi(x)$ be one-place open **Th**-wff such that the following condition is satisfied:

$$\mathbf{Th}_i^\# \vdash \exists! x_\Psi[\Psi(x_\Psi)] \tag{183}$$

**Remark 4.1.15**. We rewrite now the condition (183) in the following equivalent form using only the language of the theory $\mathbf{Th}_i^\#$ :

$$\{\mathbf{Th}_i^\# \vdash \exists! x_\Psi[\Psi(x_\Psi)]\} \Leftrightarrow \mathbf{Pr}_{\mathbf{Th}_i^\#}([\exists! x_\Psi[\Psi(x_\Psi)]]^c) \tag{184}$$

**Definition 4.1.12**. We will say that, a set $y$ is a $\mathbf{Th}_i^\#$-set if there exist one-place open wff $\Psi(x)$ such that $y = x_\Psi$. We will write for short $y[\mathbf{Th}_i^\#]$ iff $y$ is a $\mathbf{Th}_i^\#$-set.

**Remark 4.1.16**. Note that

$$y[\mathbf{Th}_i^\#] \Leftrightarrow \exists\Psi\left\{(y = x_\Psi) \wedge \mathbf{Pr}_{\mathbf{Th}_i^\#}([\exists!x_\Psi[\Psi(x_\Psi)]]^c)\right\}. \quad (185)$$

**Definition 4.1.13**. Let $\Im_i$ be a set such that :

$$\forall x\left[x \in \Im_i \leftrightarrow x \text{ is a } \mathbf{Th}_i^\#\text{-set}\right]. \quad (186)$$

**Proposition 4.1.6**. $\Im_i$ is a $\mathbf{Th}_i^\#$-set.

**Proof**. Let us consider an one-place open wff $\Psi(x)$ such that conditions (183) are satisfied, i.e. $\mathbf{Th}_i^\# \vdash \exists!x_\Psi[\Psi(x_\Psi)]$. We note that there exists countable set $\mathcal{F}_\Psi$ of the one-place open wff's $\mathcal{F}_\Psi = \{\Psi_n(x)\}_{n\in\mathbb{N}}$ such that: (i) $\Psi(x) \in \mathcal{F}_\Psi$ and (ii)

$$\mathbf{Th}_i^\# \vdash \exists!x_\Psi[[\Psi(x_\Psi)] \wedge \{\forall n(n \in \mathbb{N})[\Psi(x_\Psi) \leftrightarrow \Psi_n(x_\Psi)]\}]$$

or in the equivalent form

$$\begin{gathered}\mathbf{Th}_i^\# \vdash \mathbf{Pr}_{\mathbf{Th}_i^\#}([\exists!x_\Psi[\Psi(x_\Psi)]]^c) \wedge \\ \left\{\mathbf{Pr}_{\mathbf{Th}_i^\#}([\exists!x_\Psi[\Psi(x_\Psi)]]^c) \Rightarrow \exists!x_\Psi[\Psi(x_\Psi)]\right\} \wedge \\ \left[\mathbf{Pr}_{\mathbf{Th}_i^\#}([\forall n(n \in \mathbb{N})[\Psi(x_\Psi) \leftrightarrow \Psi_n(x_\Psi)]]^c)\right] \wedge \\ \mathbf{Pr}_{\mathbf{Th}_i^\#}([\forall n(n \in \mathbb{N})[\Psi(x_\Psi) \leftrightarrow \Psi_n(x_\Psi)]]^c) \Rightarrow \forall n(n \in \mathbb{N})[\Psi(x_\Psi) \leftrightarrow \Psi_n(x_\Psi)]\end{gathered} \quad (187)$$

or in the following equivalent form

$$\mathbf{Th}_i^\# \vdash \exists!x_1[[\Psi_1(x_1)] \wedge \{\forall n(n \in \mathbb{N})[\Psi_1(x_1) \leftrightarrow \Psi_{n,1}(x_1)]\}]$$

or

$$\begin{gathered}\mathbf{Th}_i^\# \vdash \\ \mathbf{Pr}_{\mathbf{Th}_i^\#}([\exists!x_1\Psi(x_1)]^c) \wedge \\ \left\{\mathbf{Pr}_{\mathbf{Th}_i^\#}([\exists!x_1\Psi(x_1)]^c) \Rightarrow \exists!x_1\Psi(x_1)\right\} \wedge \\ \left[\mathbf{Pr}_{\mathbf{Th}_i^\#}([\forall n(n \in \mathbb{N})[\Psi(x_1) \leftrightarrow \Psi_n(x_1)]]^c)\right] \wedge \\ \mathbf{Pr}_{\mathbf{Th}_i^\#}([\forall n(n \in \mathbb{N})[\Psi(x_1) \leftrightarrow \Psi_n(x_1)]]^c) \Rightarrow \forall n(n \in \mathbb{N})[\Psi(x_1) \leftrightarrow \Psi_n(x_1)],\end{gathered} \quad (188)$$

where we have set $\Psi(x) \triangleq \Psi_1(x_1), \Psi_n(x_1) \triangleq \Psi_{n,1}(x_1)$ and $x_\Psi \triangleq x_1$. We note that any set $\mathcal{F}_{\Psi_k} = \{\Psi_{n,k}(x)\}_{n\in\mathbb{N}}, k = 1,2,\ldots$ such as mentioned above, defines an unique set $x_{\Psi_k}$, i.e. $\mathcal{F}_{\Psi_{k_1}} \cap \mathcal{F}_{\Psi_{k_2}} = \varnothing$ iff $x_{\Psi_{k_1}} \neq x_{\Psi_{k_2}}$. We note that a sets $\mathcal{F}_{\Psi_k}, k = 1,2,..$ are part of the $ZFC_2^{Hs}$, i.e. any set $\mathcal{F}_{\Psi_k}$ is a set in the sense of $ZFC_2^{Hs}$. Note that by using Gödel numbering one can to replace any set $\mathcal{F}_{\Psi_k}, k = 1,2,..$ by set $\Theta_k = g(\mathcal{F}_{\Psi_k})$ of the corresponding Gödel numbers such that

$$\Theta_k = g(\mathcal{F}_{\Psi_k}) = \{g(\Psi_{n,k}(x_k))\}_{n\in\mathbb{N}}, k = 1,2,\ldots . \quad (189)$$

It is easy to prove that any set $\Theta_k = g(\mathcal{F}_{\Psi_k}), k = 1,2,..$ is a $\mathbf{Th}_i^\#$-set. This is done by Gödel encoding, by the statament (183) and by the axiom schema of separation.

Let $g_{n,k} = g(\Psi_{n,k}(x_k)), k = 1,2,..$ be a Gödel number of the wff $\Psi_{n,k}(x_k)$. Therefore $g(\mathcal{F}_k) = \{g_{n,k}\}_{n\in\mathbb{N}}$, where we have set $\mathcal{F}_k = \mathcal{F}_{\Psi_k}$, $k = 1,2,..$ and

$$\forall k_1\forall k_2[\{g_{n,k_1}\}_{n\in\mathbb{N}} \cap \{g_{n,k_2}\}_{n\in\mathbb{N}} = \varnothing \leftrightarrow x_{k_1} \neq x_{k_2}]. \quad (190)$$

Let $\{\{g_{n,k}\}_{n\in\mathbb{N}}\}_{k\in\mathbb{N}}$ be a family of the all sets $\{g_{n,k}\}_{n\in\mathbb{N}}$. By axiom of choice one obtains a unique set $\Im_i' = \{g_k\}_{k\in\mathbb{N}}$ such that $\forall k[g_k \in \{g_{n,k}\}_{n\in\mathbb{N}}]$. Finally for any $i \in \mathbb{N}$ one obtains a set $\Im_i$ from the set $\Im_i'$ by the axiom schema of replacement.

**Proposition 4.1.7**. Any collection $\Theta_k = g(\mathcal{F}_{\Psi_k}), k = 1,2,..$ is a $\mathbf{Th}_i^\#$-set.

**Proof**. We define $g_{n,k} = g(\Psi_{n,k}(x_k)) = [\Psi_{n,k}(x_k)]^c, v_k = [x_k]^c$. Therefore $g_{n,k} = g(\Psi_{n,k}(x_k)) \leftrightarrow \mathbf{Fr}(g_{n,k}, v_k)$. Let us define now predicate $\Pi_i(g_{n,k}, v_k)$

$$\Pi_i(g_{n,k}, v_k) \Leftrightarrow \mathbf{Pr}_{\mathbf{Th}_i^\#}([\exists! x_k[\Psi_{1,k}(x_1)]]^c) \wedge \\ \wedge \exists! x_k(v_k = [x_k]^c)\Big[\forall n(n \in \mathbb{N})\Big[\mathbf{Pr}_{\mathbf{Th}_i^\#}([[\Psi_{1,k}(x_k)]]^c) \Leftrightarrow \mathbf{Pr}_{\mathbf{Th}_i^\#}(\mathbf{Fr}(g_{n,k}, v_k))\Big]\Big]. \tag{191}$$

We define now a set $\Theta_k$ such that

$$\begin{cases} \Theta_k = \Theta'_k \cup \{g_k\}, \\ \forall n(n \in \mathbb{N})[g_{n,k} \in \Theta'_k \leftrightarrow \Pi_i(g_{n,k}, v_k)]. \end{cases} \tag{192}$$

Obviously definitions (191) and (192) are equivalent.

**Definition 4.1.14**. We define now the following $\mathbf{Th}_i^\#$-set $\Re_i \subsetneqq \Im_i$ :

$$\forall x\Big[x \in \Re_i \Leftrightarrow (x \in \Im_i) \wedge \mathbf{Pr}_{\mathbf{Th}_i^\#}([x \notin x]^c)\Big]. \tag{193}$$

**Proposition 4.1.8**. (i) $\mathbf{Th}_i^\# \vdash \exists \Re_i$, (ii) $\Re_i$ is a countable $\mathbf{Th}_i^\#$-set, $i \in \mathbb{N}$.

**Proof**. (i) Statement $\mathbf{Th}_i^\# \vdash \exists \Re_i$ follows immediately by using statement $\exists \Im_i$ and axiom schema of separation. (ii) follows immediately from countability of a set $\Im_i$.

**Proposition 4.1.9**. Any set $\Re_i, i \in \mathbb{N}$ is inconsistent.

**Proof**. From the formula (193) we obtain

$$\mathbf{Th}_i^\# \vdash \Re_i \in \Re_i \Leftrightarrow \mathbf{Pr}_{\mathbf{Th}_i^\#}([\Re_i \notin \Re_i]^c). \tag{194}$$

From the formula (194) we obtain

$$\mathbf{Th}_i^\# \vdash \Re_i \in \Re_i \Leftrightarrow \Re_i \notin \Re_i \tag{195}$$

and therefore

$$\mathbf{Th}_i^\# \vdash (\Re_i \in \Re_i) \wedge (\Re_i \notin \Re_i). \tag{196}$$

But this is a contradiction.

**Definition 4.1.15**. A $\mathbf{Th}_\infty^\#$-wff $\Phi_\infty$ that is: (i) $\mathbf{Th}$-wff $\Phi$ or (ii) well-formed formula $\Phi_\infty$ which contains predicate $\mathbf{Pr}_{\mathbf{Th}_\infty^\#}([\Phi]^c)$ given by formula (4.1.28). An $\mathbf{Th}_\infty^\#$-wff $\Phi_\infty$ (well-formed formula $\Phi_\infty$) is closed, i.e. $\Phi_\infty$ is a sentence if $\Phi_\infty$ has no free variables; a wff is open if it has free variables.

**Definition 4.1.16**. Let $\Psi = \Psi(x)$ be one-place open $\mathbf{Th}_\infty^\#$-wff such that the following condition:

$$\mathbf{Th}_\infty^\# \vdash \exists! x_\Psi[\Psi(x_\Psi)] \tag{197}$$

is satisfied.

**Remark 4.1.17**. We rewrite now the condition (197) in the following equivalent form using only the language of the theory $\mathbf{Th}_\infty^\#$ :

$$\{\mathbf{Th}_\infty^\# \vdash \exists! x_\Psi[\Psi(x_\Psi)]\} \Leftrightarrow \mathbf{Pr}_{\mathbf{Th}_\infty^\#}([\exists! x_\Psi[\Psi(x_\Psi)]]^c) \tag{198}$$

**Definition 4.1.17**. We will say that, a set $y$ is a $\mathbf{Th}_\infty^\#$-set if there exists one-place open wff $\Psi(x)$ such that $y = x_\Psi$. We write $y[\mathbf{Th}_\infty^\#]$ iff $y$ is a $\mathbf{Th}_\infty^\#$-set.

**Definition 4.1.18**. Let $\Im_\infty$ be a set such that : $\forall x\big[x \in \Im_\infty \leftrightarrow x$ is a $\mathbf{Th}_\infty^\#$-set$\big]$.

**Proposition 4.1.10**. A set $\Im_\infty$ is a $\mathbf{Th}_\infty^\#$-set.

**Proof**. Let us consider an one-place open wff $\Psi(x)$ such that condition (197) is satisfied, i.e. $\mathbf{Th}_\infty^\# \vdash \exists! x_\Psi[\Psi(x_\Psi)]$. We note that there exists countable set $\mathcal{F}_\Psi$ of the one-place open wff's $\mathcal{F}_\Psi = \{\Psi_n(x)\}_{n \in \mathbb{N}}$ such that: (i) $\Psi(x) \in \mathcal{F}_\Psi$ and (ii)

$$\mathbf{Th}_{\infty}^{\#} \vdash \exists! x_{\Psi}[[\Psi(x_{\Psi})] \wedge \{\forall n(n \in \mathbb{N})[\Psi(x_{\Psi}) \leftrightarrow \Psi_n(x_{\Psi})]\}]$$

or in the equivalent form

$$\begin{gathered}\mathbf{Th}_{\infty}^{\#} \vdash \mathbf{Pr}_{\mathbf{Th}_{\infty}^{\#}}([\exists! x_{\Psi}[\Psi(x_{\Psi})]]^c) \wedge \\ \{\mathbf{Pr}_{\mathbf{Th}_{\infty}^{\#}}([\exists! x_{\Psi}[\Psi(x_{\Psi})]]^c) \Rightarrow \exists! x_{\Psi}[\Psi(x_{\Psi})]\} \wedge \\ [\mathbf{Pr}_{\mathbf{Th}_{\infty}^{\#}}([\forall n(n \in \mathbb{N})[\Psi(x_{\Psi}) \leftrightarrow \Psi_n(x_{\Psi})]]^c)] \wedge \\ \mathbf{Pr}_{\mathbf{Th}_{\infty}^{\#}}([\forall n(n \in \mathbb{N})[\Psi(x_{\Psi}) \leftrightarrow \Psi_n(x_{\Psi})]]^c) \Rightarrow \forall n(n \in \mathbb{N})[\Psi(x_{\Psi}) \leftrightarrow \Psi_n(x_{\Psi})]\end{gathered} \tag{199}$$

or in the following equivalent form

$$\mathbf{Th}_{\infty}^{\#} \vdash \exists! x_1[[\Psi_1(x_1)] \wedge \{\forall n(n \in \mathbb{N})[\Psi_1(x_1) \leftrightarrow \Psi_{n,1}(x_1)]\}]$$

or

$$\begin{gathered}\mathbf{Th}_{\infty}^{\#} \vdash \mathbf{Pr}_{\mathbf{Th}_i^{\#}}([\exists! x_1\Psi(x_1)]^c) \wedge \\ \{\mathbf{Pr}_{\mathbf{Th}_{\infty}^{\#}}([\exists! x_1\Psi(x_1)]^c) \Rightarrow \exists! x_1\Psi(x_1)\} \wedge \\ \left[\mathbf{Pr}_{\mathbf{Th}_i^{\#}}([\forall n(n \in \mathbb{N})[\Psi(x_1) \leftrightarrow \Psi_n(x_1)]]^c)\right] \wedge \\ \mathbf{Pr}_{\mathbf{Th}_i^{\#}}([\forall n(n \in \mathbb{N})[\Psi(x_1) \leftrightarrow \Psi_n(x_1)]]^c) \Rightarrow \forall n(n \in \mathbb{N})[\Psi(x_1) \leftrightarrow \Psi_n(x_1)].\end{gathered} \tag{200}$$

where we set $\Psi(x) = \Psi_1(x_1), \Psi_n(x_1) = \Psi_{n,1}(x_1)$ and $x_{\Psi} = x_1$. We note that any set $\mathcal{F}_{\Psi_k} = \{\Psi_{n,k}(x)\}_{n\in\mathbb{N}}, k = 1,2,\ldots$ such as mentioned above defines a unique set $x_{\Psi_k}$, i.e. $\mathcal{F}_{\Psi_{k_1}} \bigcap \mathcal{F}_{\Psi_{k_2}} = \varnothing$ iff $x_{\Psi_{k_1}} \neq x_{\Psi_{k_2}}$. We note that sets $\mathcal{F}_{\Psi_k}, k = 1,2,..$ are the part of the $ZFC_2^{Hs}$, i.e. a set $\mathcal{F}_{\Psi_k}$ is a set in the sense of $ZFC_2^{Hs}$. Note that by using Gödel numbering one can replace any set $\mathcal{F}_{\Psi_k}, k = 1,2,..$ by the set $\Theta_k = g(\mathcal{F}_{\Psi_k})$ of the corresponding Gödel numbers such that

$$\Theta_k = g(\mathcal{F}_{\Psi_k}) = \{g(\Psi_{n,k}(x_k))\}_{n\in\mathbb{N}}, k = 1,2,\ldots . \tag{201}$$

It is easy to prove that any set $\Theta_k = g(\mathcal{F}_{\Psi_k}), k = 1,2,..$ is a $\mathbf{Th}^{\#}$-set. This is done by Gödel encoding and by axiom schema of separation. Let $g_{n,k} = g(\Psi_{n,k}(x_k)), k = 1,2,..$ be a Gödel number of the wff $\Psi_{n,k}(x_k)$. Therefore $g(\mathcal{F}_k) = \{g_{n,k}\}_{n\in\mathbb{N}}$, where we have set $\mathcal{F}_k \triangleq \mathcal{F}_{\Psi_k}$, $k = 1,2,..$ and

$$\forall k_1 \forall k_2[\{g_{n,k_1}\}_{n\in\mathbb{N}} \bigcap \{g_{n,k_2}\}_{n\in\mathbb{N}} = \varnothing \leftrightarrow x_{k_1} \neq x_{k_2}]. \tag{202}$$

Let $\{\{g_{n,k}\}_{n\in\mathbb{N}}\}_{k\in\mathbb{N}}$ be a family of the sets $\{g_{n,k}\}_{n\in\mathbb{N}}, k = 1,2,\ldots$ . By axiom of choice one obtains an unique set $\mathfrak{I}' = \{g_k\}_{k\in\mathbb{N}}$ such that $\forall k[g_k \in \{g_{n,k}\}_{n\in\mathbb{N}}]$. Finally one obtains a set $\mathfrak{I}_\infty$ from the set $\mathfrak{I}'_\infty$ by axiom schema of replacement. Thus we can define $\mathbf{Th}_{\infty}^{\#}$-set $\mathfrak{R}_\infty \subsetneqq \mathfrak{I}_\infty$ :

$$\forall x[x \in \mathfrak{R}_\infty \leftrightarrow (x \in \mathfrak{I}_\infty) \wedge [\mathbf{Pr}_{\mathbf{Th}_{\infty}^{\#}}([x \notin x]^c)]]. \tag{203}$$

**Proposition 4.1.11**. Any set $\Theta_k = g(\mathcal{F}_{\Psi_k}), k = 1,2,..$ is a $\mathbf{Th}_{\infty}^{\#}$-set.

**Proof**. We define $g_{n,k} = g(\Psi_{n,k}(x_k)) = [\Psi_{n,k}(x_k)]^c, v_k = [x_k]^c$. Therefore $g_{n,k} = g(\Psi_{n,k}(x_k)) \leftrightarrow \mathbf{Fr}(g_{n,k}, v_k)$. Let us define now predicate $\Pi_\infty(g_{n,k}, v_k)$

$$\begin{gathered}\Pi_\infty(g_{n,k}, v_k) \Leftrightarrow \\ \mathbf{Pr}_{\mathbf{Th}_{\infty}^{\#}}([\exists! x_k[\Psi_{1,k}(x_1)]]^c) \wedge [\mathbf{Pr}_{\mathbf{Th}_{\infty}^{\#}}([\exists! x_k[\Psi_{1,k}(x_1)]]^c) \Rightarrow \exists! x_1\Psi(x_1)] \\ \wedge \exists! x_k(v_k = [x_k]^c)[\forall n(n \in \mathbb{N})[\mathbf{Pr}_{\mathbf{Th}_{\infty}^{\#}}([[\Psi_{1,k}(x_k)]]^c) \Leftrightarrow \mathbf{Pr}_{\mathbf{Th}_{\infty}^{\#}}(\mathbf{Fr}(g_{n,k}, v_k))]].\end{gathered} \tag{204}$$

We define now a set $\Theta_k$ such that

$$\begin{gathered}\Theta_k = \Theta'_k \cup \{g_k\}, \\ \forall n(n \in \mathbb{N})[g_{n,k} \in \Theta'_k \Leftrightarrow \Pi(g_{n,k}, v_k)]\end{gathered} \tag{205}$$

Obviously definitions (204) and (205) are equivalent by Proposition 4.1.1.

**Proposition 4.1.12**. (i) $\mathbf{Th}_{\infty}^{\#} \vdash \exists \mathfrak{R}_\infty$, (ii) $\mathfrak{R}_\infty$ is a countable $\mathbf{Th}_{\infty}^{\#}$-set.

**Proof**. (i) Statement $\mathbf{Th}_{\infty}^{\#} \vdash \exists \mathfrak{R}_\infty$ follows immediately from the statement $\exists \mathfrak{I}_\infty$ and axiom schema of separation [9], (ii) follows immediately from countability of the set $\Im_\infty$.

**Proposition 4.1.13**. Set $\Re_\infty$ is inconsistent.

**Proof**. From the formula (203) we obtain

$$\mathbf{Th}_\infty^\# \vdash \Re_\infty \in \Re_\infty \Leftrightarrow \mathbf{Pr}_{\mathbf{Th}_\infty^\#}([\Re_\infty \notin \Re_\infty]^c). \quad (206)$$

From (206) one obtains

$$\mathbf{Th}_\infty^\# \vdash \Re_\infty \in \Re_\infty \Leftrightarrow \Re_\infty \notin \Re_\infty \quad (207)$$

and therefore

$$\mathbf{Th}_\infty^\# \vdash (\Re_\infty \in \Re_\infty) \wedge (\Re_\infty \notin \Re_\infty). \quad (208)$$

But this is a contradiction.

**Remark 4.1.18**. Note that a contradictions mentioned above can be again avoid using canonical Quinean approach, see subsection 3.6.

## 5mm . 4.2. Proof of the inconsistency of the set theory $ZFC_2^{Hs} + \exists M_{st}^{ZFC_2^{Hs}}$ using generalized Tarski's undefinability theorem 2mm

In this section we will prove that a set theory $ZFC_2^{Hs} + \exists M^{ZFC_2^{Hs}}$ is inconsistent, without any reference to the sets $\Im_1, \Im_2, \ldots, \Im_\infty$ and corresponding inconsistent sets $\Re_1, \Re_2, \ldots, \Re_\infty$.

**Remark 4.2.1**. Note that a contradiction mentioned above is a strictly stronger then contradictions derived in subsection 4.1, and these contradictions are impossible to avoid by using Quinean approach, see subsection 3.6.

**Proposition 4.2.1**. (Generalized Tarski's undefinability theorem). Let $\mathbf{Th}_{\mathcal{L}}^{Hs}$ be second order theory with Henkin semantics and with formal language $\mathcal{L}$, which includes negation and has a Gödel encoding $g(\cdot)$ such that for every $\mathcal{L}$-formula $A(x)$ there is a formula $B$ such that $B \Leftrightarrow A(g(B))$ holds. Assume that $\mathbf{Th}_{\mathcal{L}}^{Hs}$ has a standard Model $M^{ZFC_2^{Hs}}$.

Then there is no $\mathcal{L}$-formula $\mathbf{True}(n)$ such that for every $\mathcal{L}$-formula $A$ such that $M^{ZFC_2^{Hs}} \models A$, the following equivalence holds

$$\left(M^{ZFC_2^{Hs}} \models A\right) \Leftrightarrow \mathbf{True}(g(A)). \quad (209)$$

**Proof**. The diagonal lemma yields a counterexample to this equivalence, by giving a "Liar" sentence $S$ such that $S \Leftrightarrow \neg\mathbf{True}(g(S))$ holds.

**Remark 4.2.2**. Above we has been defined the set $\Im_\infty$(see Definition 4.1.16) in fact using generalized truth predicate $\mathbf{True}_\infty^\#([\Phi]^c)$ such that

$$\mathbf{True}_\infty^\#([\Phi]^c) \Leftrightarrow \mathbf{Pr}_{\mathbf{Th}_\infty^\#}([\Phi]^c). \quad (210)$$

In order to prove that set theory $ZFC_2^{Hs} + \exists M^{ZFC_2^{Hs}}$ is inconsistent without any reference to the set $\Im_\infty$, notice that by the properties of the nice extension $\mathbf{Th}_\infty^\#$ it follows that definition given by biconditional (211) is correct, i.e., for every first order $ZFC_2^{Hs}$-formula $\Phi$ such that $M^{ZFC_2^{Hs}} \models \Phi$ and the following equivalence holds

$$\left(M^{ZFC_2^{Hs}} \models \Phi\right) \Leftrightarrow \mathbf{Pr}_{\mathbf{Th}_\infty^\#}([\Phi]^c), \quad (211)$$

where $\mathbf{Pr}_{\mathbf{Th}_\infty^\#}([\Phi]^c) \Rightarrow \Phi$.

**Proposition 4.2.2**. Set theory $\mathbf{Th}_1^\# = ZFC_2^{Hs} + \exists M^{ZFC_2^{Hs}}$ is inconsistent.

**Proof**. Notice that by the properties of the nice extension $\mathbf{Th}_\infty^\#$ of the $\mathbf{Th}_1^\#$ it follows that

$$\left(M^{ZFC_2^{Hs}} \models \Phi\right) \Leftrightarrow \mathbf{Th}_\infty^\# \vdash \Phi. \quad (212)$$

Therefore (210) gives generalized "truth predicate" for set theory $\mathbf{Th}_\infty^\#$. By Proposition 4.2.1 one obtains a contradiction.

**Remark 4.2.3**. A cardinal $\kappa$ is inaccessible iff $\kappa$ has the following reflection property: for all subsets $U \subset V_\kappa$, there exists $\alpha < \kappa$ such that $(V_\alpha, \in, U \cap V_\alpha)$ is an elementary substructure of

$(V_\kappa, \in, U)$. (In fact, the set of such $\alpha$ is closed unbounded in $\kappa$.)

Equivalently, $\kappa$ is $\Pi^0_n$ -indescribable for all $n \geq 0$.

**Remark 4.2.4**. Under *ZFC* it can be shown that $\kappa$ is inaccessible iff $(V_\kappa, \in)$ is a model of second order *ZFC* [5]cite: Foukzon15.

**Remark 4.2.5**. By the reflection property, there exists $\alpha < \kappa$ such that $(V_\alpha, \in)$ is a standard model of (first order) *ZFC*. Hence, the existence of an inaccessible cardinal is a stronger hypothesis than the existence of the standard model of $ZFC_2^{Hs}$.

## 5mm . 4.3. Derivation inconsistent countable set in set theory $ZFC_2$ with the full semantics 2mm

Let $\mathbf{Th} = \mathbf{Th}^{fss}$ be a second order theory with the full second order semantics. We assume now that $\mathbf{Th}$ contains $ZFC_2^{fss}$. We will write for short $\mathbf{Th}$, instead $\mathbf{Th}^{fss}$.

**Remark 4.3.1**. Notice that $M$ is a model of $ZFC_2^{fss}$ iff it is isomorphic to a model of the form $V_\kappa, \in \cap (V_\kappa \times V_\kappa)$, for $\kappa$ a strongly inaccessible ordinal.

**Remark 4.3.2**. Notice that a standard model for the language of first-order set theory is an ordered pair $\{D, I\}$ .Its domain, $D$, is a nonempty set and its interpretation function, $I$, assigns a set of ordered pairs to the two-place predicate "$\in$". A sentence is true in $\{D, I\}$ just in case it is satisfied by all assignments of first-order variables to members of $D$ and second-order variables to subsets of $D$; a sentence is satisfiable just in case it is true in some standard model; finally, a sentence is valid just in case it is true in all standard models.

**Remark 4.3.3**. Notice that:

(**I**) The assumption that $D$ and $I$ be sets is not without consequence. An immediate effect of this stipulation is that no standard model provides the language of set theory with its intended interpretation. In other words, there is no standard model $\{D, I\}$ in which $D$ consists of all sets and $I$ assigns the standard element-set relation to "$\in$". For it is a theorem of *ZFC* that there is no set of all sets and that there is no set of ordered-pairs $\{x, y\}$ for $x$ an element of $y$.

(**II**) Thus, on the standard definition of model:

(1) it is not at all obvious that the validity of a sentence is a guarantee of its truth;

(2) similarly, it is far from evident that the truth of a sentence is a guarantee of its satisfiability in some standard model;

(3) if there is a connection between satisfiability, truth, and validity, it is not one that can be “read off” standard model theory.

(**III**) Nevertheless this is not a problem in the first-order case since set theory provides us with two reassuring results for the language of first-order set theory. One result is the first order completeness theorem according to which first-order sentences are provable, if true in all models. Granted the truth of the axioms of the first-order predicate calculus and the truth preserving character of its rules of inference, we know that a sentence of the first-order language of set theory is true, if it is provable. Thus, since valid sentences are provable and provable sentences are true, we know that valid sentences are true. The connection between truth and satisfiability immediately follows: if $\phi$ is unsatisfiable, then $\neg\phi$, its negation, is true in all models and hence valid. Therefore, $\neg\phi$ is true and $\phi$ is false.

**Definition 4.3.1**. The language of second order arithmetic $Z_2$ is a two-sorted language: there are two kinds of terms, numeric terms and set terms.

**0** is a numeric term.

1. There are in nitely many numeric variables, $x_0, x_1, \ldots, x_n, \ldots$ each of which is a numeric term.

2. If $s$ is a numeric term then $\mathbf{S}s$ is a numeric term.

3. If $s, t$ are numeric terms then $+st$ and $\cdot st$ are numeric terms (abbreviated $s + t$ and $s \cdot t$).

4. There are infinitely many set variables, $X_0, X_1, \ldots, X_n \ldots$ each of which is a set term;

5. If $t$ is a numeric term and $S$ then $\in tS$ is an atomic formula (abbreviated $t \in S$).

6. If s and t are numeric terms then $= st$ and $< st$ are atomic formulas (abbreviated $s = t$ and $s < t$ correspondingly).

The formulas are built from the atomic formulas in the usual way.

As the examples in the definition suggest, we use upper case letters for set variables and lower case letters for numeric terms. (Note that the only set terms are the variables.) It will be more convenient to work with functions instead of sets, but within arithmetic, these are equivalent: one can use the pairing operation, and say that $X$ represents a function if for each $n$ there is exactly one $m$ such that the pair $(n,m)$ belongs to $X$.

We have to consider what we intend the semantics of this language to be. One possibility is the semantics of full second order logic: a model consists of a set $M$, representing the numeric objects, and interpretations of the various functions and relations (probably with the requirement that equality be the genuine equality relation), and a statement $\forall X \Phi(X)$ is satisfied by the model if for every possible subset of $M$, the corresponding statement holds.

**Remark 4.3.4**. Full second order logic has no corresponding proof system. An easy way to see this is to observe that it has no compactness theorem. For example, the only model (up to isomorphism) of Peano arithmetic together with the second order induction axiom: $\forall X(\mathbf{0} \in X \wedge \forall x(x \in X \Rightarrow \mathbf{S}x \in X) \Rightarrow \forall x(x \in X))$ is the standard model $\mathbb{N}$. This is easily seen: any model of Peano arithmetic has an initial segment isomorphic to $\mathbb{N}$; applying the induction axiom to this set, we see that it must be the whole of the model.

**Remark 4.3.5**.There is no completeness theorem for second-order logic. Nor do the axioms of second-order $ZFC$ imply a reflection principle which ensures that if a sentence of second-order set theory is true, then it is true in some standard model. Thus there may be sentences of the language of second-order set theory that are true but unsatisfiable, or sentences that are valid, but false. To make this possibility vivid, let $Z$ be the conjunction of all the axioms of second-order $ZFC$. $Z$ is "surely" true. But the existence of a model for $Z$ requires the existence of strongly inaccessible cardinals.

The axioms of second-order $ZFC$ don't entail the existence of strongly inaccessible cardinals, and hence the satisfiability of $Z$ is independent of second-order $ZFC$. Thus, $Z$ is true but its unsatisfiability is consistent with second-order $ZFC$ [5]cite: Foukzon15.

**Remark 4.3.6**. We remind that urlogic is the aspect of mathematicians' activity that consists of just writing down finite strings of symbols-sentences-according to some fixed rules. Those sentences are sentences of urlogic. Whether a string of symbols is a sentence of urlogic should be totally unproblematic. In summary, urlogic has the following characteristics [13]cite: Vaananen01:

(i) Sentences of urlogic are finite strings of symbols. That a string of symbols is a sentence of urlogic, is a non-mathematical judgement.

(ii) Some sentences are accepted as axioms. That a sentence is an axiom, is a non-mathematical judgement.

(iii) Derivations are made from axioms. The derivations obey certain rules of proof. That a derivation obeys the rules of proof, is a non-mathematical judgement.

(iv) Derived sentences can be asserted as facts.

If we take first-order set theory as the urlogic, the sentences of urlogic are sentences of first-order predicate logic with identity, with the binary predicate symbol $e$ as the only non-logical symbol. The axioms are the usual rules of first-order logic augmented with the Zermelo-Fraenkel axioms $ZFC$ of set theory. On the informal level we interpret the sentences of this urlogic as propositions about mathematical objects construed as sets.

**Remark 4.3.7**. In the case of second-order logic the sentences of urlogic are the sentences of second-order predicate logic. Depending on the context, the non-logical vocabulary may consist of symbols for the arithmetic of natural numbers, arithmetic of real numbers, and so forth.

Montague [21]cite: Montague65 gives second-order Peano axioms $Z_2$ for number theory, and second-order axioms $RCF_2$ for real closed fields. For full second-order logic there is a notion of "semantical" derivation:

We can derive $\Psi$ from $\Phi$ if every model of $\Phi$ is a model of $\Psi$.

Of course scanning through all models of $\Phi$ is a highly mathematical act.

Thus with respect to $ZFC_2^{fss}$, this is a semantically defined system and thus it is not standard to speak about it being contradictory if anything, one might attempt to prove that it has no models, which to be what is being done in section 3 and section 4 for $ZFC_2^{Hs}$.

**Remark 4.3.8**. Note that in order to avoid difficulties with "semantical" derivation mentioned above one considers first order theory $\mathbf{Th} \subset ZFC_2^{fss}$ which contains only first order wff of $ZFC_2^{Hs}$. Thus in order to prove that $ZFC_2^{fss}$ has no models or it being contradictory, one might use the same approach, which is done in section 3 and section 4 for $ZFC_2^{Hs}$.

**Definition 4.3.2**. Let $\Phi$ be a wff of $\overline{ZFC}_2^{Hs}$. We will say that $\Phi$ is a first order $n$ –place open wff if $\Phi$ contains free occurrences of the first order individual variables $X_1,\ldots,X_n$ and quantifiers only over any first order individual variables $Y_1,\ldots,Y_m$.

**Definition 4.3.3**. Let **Th** be a first order theory which contains only first order wff of $ZFC_2^{Hs}$. Using formula (141) one can define predicate $\mathbf{Pr}^{\#}_{\mathbf{Th}}(y)$ really asserting provability of the first order sentences in $\mathbf{Th} \subset ZFC_2^{fss}$:

$$\begin{gathered}\mathbf{Pr}^{\#}_{\mathbf{Th}}(y) \Leftrightarrow \mathbf{Pr}_{\mathbf{Th}}(y) \wedge [\mathbf{Pr}_{\mathbf{Th}}(y) \Rightarrow \Phi], \\ \mathbf{Pr}_{\mathbf{Th}}(y) \Leftrightarrow \exists x\left(x \in M_{\omega}^{Z_2^{fss}}\right)\mathbf{Prov}_{\mathbf{Th}}(x,y), \\ y = [\Phi]^c.\end{gathered} \tag{213}$$

**Theorem 4.3.1**. [12]**cite: RayoUzquiano99** . (Löb's Theorem for $ZFC_2^{fss}$.) Let $\Phi$ be any first order closed formula with code $y = [\Phi]^c \in M_{\omega}^{Z_2}$, then $\mathbf{Th} \vdash \mathbf{Pr}_{\mathbf{Th}}([\Phi]^c)$ implies $\mathbf{Th} \vdash \Phi$.

**Proof**. Assume that

(#) $\mathbf{Th} \vdash \mathbf{Pr}_{\mathbf{Th}}([\Phi]^c)$.

Note that

(1) $\mathbf{Th} \nvdash \neg\Phi$. Otherwise one obtains $\mathbf{Th} \vdash \mathbf{Pr}_{\mathbf{Th}}([\neg\Phi]^c) \wedge \mathbf{Pr}_{\mathbf{Th}}([\Phi]^c)$, but this is a contradiction.

(2) Assume now that (2.i) $\mathbf{Th} \vdash \mathbf{Pr}_{\mathbf{Th}}([\Phi]^c)$ and (2.ii) $\mathbf{Th} \nvdash \Phi$.

From (1) and (2.ii) it follows that

(3) $\mathbf{Th} \nvdash \neg\Phi$ and $\mathbf{Th} \nvdash \Phi$.

Let $\mathbf{Th}_{\neg\Phi}$ be a theory

(4)$\mathbf{Th}_{\neg\Phi} \triangleq \mathbf{Th} \cup \{\neg\Phi\}$. From (3) it follows that

(5) $Con(\mathbf{Th}_{\neg\Phi})$.

From (4) and (5) it follows that

(6) $\mathbf{Th}_{\neg\Phi} \vdash \mathbf{Pr}_{\mathbf{Th}_{\neg\Phi}}([\neg\Phi]^c)$.

From (4) and (#) it follows that

(7) $\mathbf{Th}_{\neg\Phi} \vdash \mathbf{Pr}_{\mathbf{Th}_{\neg\Phi}}([\Phi]^c)$.

From (6) and (7) follows that

(8) $\mathbf{Th}_{\neg\Phi} \vdash \mathbf{Pr}_{\mathbf{Th}_{\neg\Phi}}([\Phi]^c) \wedge \mathbf{Pr}_{\mathbf{Th}_{\neg\Phi}}([\neg\Phi]^c)$, but this is a contradiction.

**Definition 4.3.4**. Let $\Psi = \Psi(x)$ be first order one-place open wff such that:

$$\mathbf{Th} \vdash \exists! x_{\Psi}[\Psi(x_{\Psi})]. \tag{214}$$

Then we will say that, a set $y$ is a **Th**-set iff there is exist first order one-place open wff $\Psi(x)$ such that $y = x_{\Psi}$. We write $y[\mathbf{Th}]$ iff $y$ is a **Th**-set.

**Remark 4.3.9**. Note that

$$\begin{gathered} y[\mathbf{Th}] \Leftrightarrow \\ \exists\Psi[(y = x_{\Psi}) \wedge \mathbf{Pr}_{\mathbf{Th}}([\exists! x_{\Psi}[\Psi(x_{\Psi})]]^c)] \end{gathered} \tag{215}$$

**Definition 4.3.5**. Let $\Im$ be a collection such that : $\forall x\left[x \in \Im \leftrightarrow x \text{ is a } \mathbf{Th}\text{-set}\right]$.

**Proposition 4.3.1**. A set $\Im$ is a **Th**-set.

**Definition 4.3.6**. We define now a **Th**-set $\Re_c \subsetneqq \Im$ :

$$\forall x[x \in \Re_c \leftrightarrow (x \in \Im) \wedge \mathbf{Pr_{Th}}([x \notin x]^c)]. \quad (216)$$

**Proposition 4.3.2**. (i) **Th** $\vdash \exists \Re_c$, (ii) $\Re_c$ is a countable **Th**-set.

**Proof**. (i) Statement **Th** $\vdash \exists \Re_c$ follows immediately by using statement $\exists \Im$ and axiom schema of separation [4], (ii) follows immediately from countability of a set $\Im$.

**Proposition 4.3.3**. A set $\Re_c$ is inconsistent.

**Proof**. From formula (216) one obtains

$$\mathbf{Th} \vdash \Re_c \in \Re_c \Longleftrightarrow \mathbf{Pr_{Th}}([\Re_c \notin \Re_c]^c). \quad (217)$$

From formula (216) and definition 4.3.5 one obtains

$$\mathbf{Th} \vdash \Re_c \in \Re_c \Longleftrightarrow \Re_c \notin \Re_c \quad (218)$$

and therefore

$$\mathbf{Th} \vdash (\Re_c \in \Re_c) \wedge (\Re_c \notin \Re_c). \quad (219)$$

But this is a contradiction.

Thus finally we obtain:

**Theorem 4.3.2**. [5]**cite: Foukzon15**.$\neg Con(ZFC_2^{fss})$.

It well known that under $ZFC$ it can be shown that $\kappa$ is inaccessible iff $(V_\kappa, \in)$ is a model of $ZFC_2$ [12]. Thus finally we obtain.

**Theorem 4.3.3**. [5], [6]**cite: Foukzon15cite: Foukzon19**.$\neg Con(ZFC + \exists M_{st}^{ZFC}(M_{st}^{ZFC} = H_k))$.

# 6mm . 5.Discussion. How can we safe the set theory $ZFC + \exists M_{st}^{ZFC}$ 3mm

## 5mm . 5.1. The set theory.$ZFC_w$ with a weakened axiom of infinity 2mm

We remind that a major part of modern mathematical analysis and related areas based not only on set theory $ZFC$ but on strictly stronger set theory: $ZFC + \exists M_{st}^{ZFC}$. In order to avoid difficultnes which arises from $\neg Con(ZFC + \exists M_{st}^{ZFC})$ in this subsection we introduce the set theory $ZFC_w$ with a weakened axiom of infinity. Without loss of generality we consider second-order arithmetic $\mathbb{Z}_2$ with an restricted induction schema.

Second-order arithmetic $\mathbb{Z}_2$ includes, but is significantly stronger than, its first-order counterpart Peano arithmetic. Unlike Peano arithmetic, second-order arithmetic allows quantification over sets of natural numbers as well as numbers themselves. Because real numbers can be represented as (infinite) sets of natural numbers in well-known ways, and because second order arithmetic allows quantification over such sets, it is possible to formalize the real numbers in second-order arithmetic. For this reason, second-order arithmetic is sometimes called “analysis".

**Induction schema of second-order arithmetic $\mathbb{Z}_2$.**

If $\varphi(n)$ is a formula of second-order arithmetic $\mathbb{Z}_2$ with a free number variable $n$ and possible other free number or set variables (written $m$ and $X$), the induction axiom for $\varphi$ is the axiom:

$$\forall m \forall X((\varphi(0) \wedge \forall n(\varphi(n) \to \varphi(n+1))) \to \forall n \varphi(n)). \quad (220)$$

The (full) second-order induction scheme consists of all instances of this axiom, over all second-order formulas. One particularly important instance of the induction scheme is when $\varphi$ is the formula “$n \in X$” expressing the fact that n is a member of $X$ ($X$ being a free set variable): in this case, the induction axiom for $\varphi$ is

$$\forall X((0 \in X \wedge \forall n(n \in X \to n+1 \in X)) \to \forall n(n \in X)). \quad (221)$$

This sentence is called the second-order induction axiom.

**Comprehension schema of second-order arithmetic** $\mathbb{Z}_2$.

If $\varphi(n)$ is a formula with a free variable n and possibly other free variables, but not the variable $Z$, the comprehension axiom for $\varphi$ is the formula

$$\exists Z \forall n(n \in Z \leftrightarrow \varphi(n)). \tag{222}$$

This axiom makes it possible to form the set $Z = \{n | \varphi(n)\}$ of natural numbers satisfying $\varphi(n)$. There is a technical restriction that the formula $\varphi$ may not contain the variable $Z$.

**Designation 5.1.1**. Let $\mathbf{Wff}_k(\mathbb{Z}_2)$ be a set of the all k-place open wff's of the second-order arithmetic $\mathbb{Z}_2$ and let $\Re_k(\mathbb{Z}_2)$ be a set of the all primitive recursve k-place open wff's $\psi_{\Re_k}$ of the second-order arithmetic $\mathbb{Z}_2$. Let $\mathbf{\Sigma}_k(\mathbb{Z}_2)$ be a set of the all k-place open wff's $\psi_{\Sigma_k}$ of the second-order arithmetic $\mathbb{Z}_2$ such that

$$\Re_k \triangleq \Re_k(\mathbb{Z}_2) \subseteqq \mathbf{\Sigma}_k(\mathbb{Z}_2) \subsetneqq \mathbf{Wff}_k(\mathbb{Z}_2). \tag{223}$$

Let $\widetilde{\mathbf{Wff}}_{1,X}$ be a set of the all sets definable by 1-place open wff's $\psi(X) \in \mathbf{Wff}_{1,X}(\mathbb{Z}_2)$,
let $\widetilde{\mathbf{\Sigma}}_1$ be a set of the all sets definable by 1-place open wff's $\psi_{\mathbf{\Sigma}_1}(X) \in \mathbf{\Sigma}_1(\mathbb{Z}_2)$ and
let $\widetilde{\Re}_1$be a set of the all sets definable by 1-place open wff's $\psi_{\Re_1}(X) \in \Re_1(\mathbb{Z}_2)$

**Restricted induction schema of second-order arithmetic** $\mathbb{Z}_2^{\Sigma}$.

If $\varphi_{\Sigma_k}(n) \in \mathbf{\Sigma}_k \triangleq \mathbf{\Sigma}(\mathbb{Z}_2)$ is a formula of second-order arithmetic $\mathbb{Z}_2$ with a free number variable $n$ and possible other free number and set variables (written $m$ and $X$), the induction axiom for $\varphi_\Sigma$ is the axiom:

$$\forall m \forall X\left(X \in \widetilde{\mathbf{\Sigma}}_1\right)((\varphi_{\Sigma_k}(0) \wedge \forall n(\varphi_{\Sigma_k}(n) \to \varphi_{\Sigma_k}(n+1))) \to \forall n \varphi_{\Sigma_k}(n)). \tag{224}$$

The restricted second-order induction scheme consists of all instances of this axiom, over all second-order formulas. One particularly important instance of the induction scheme is when $\varphi_{\Sigma_k} \in \Sigma$ is the formula $(n \in X) \wedge \left(X \in \widetilde{\mathbf{\Sigma}}_1\right)$ expressing the fact that $n$ is a member of $X$ and $X \in \widetilde{\mathbf{\Sigma}}_1$ ($X$ being a free set variable): in this case, the induction axiom for $\varphi_{\Sigma_k}$ is

$$\forall X\left(X \in \widetilde{\mathbf{\Sigma}}_1\right)((0 \in X \wedge \forall n((n \in X) \to (n+1 \in X))) \to \forall n((n \in X))). \tag{225}$$

**Restricted comprehension schema of second-order arithmetic** $\mathbb{Z}_2^{\Sigma_k}$.

If $\varphi_{\Sigma_1}(n) \in \Sigma_1$ is a formula with a free variable $n$ and possibly other free variables, but not the variable $Z$, the comprehension axiom for $\varphi_{\Sigma_1}$ is the formula

$$\exists Z \forall n(n \in Z \leftrightarrow \varphi_{\Sigma_1}(n)). \tag{226}$$

**Remark 5.1.1.**Let $\widetilde{\mathbb{Z}}_2^{\Sigma_k}$ be a theory $\mathbb{Z}_2^{\Sigma_k} + \exists M_{st}\left[\mathbb{Z}_2^{\Sigma_k}\right]$ where $M_{st}\left[\mathbb{Z}_2^{\Sigma_k}\right]$ is a standard model of $\mathbb{Z}_2^{\Sigma}$.

We assume now that

$$Con\left(\mathbb{Z}_2^{\Sigma_k} + \exists M_{st}\left[\mathbb{Z}_2^{\Sigma_k}\right]\right). \tag{227}$$

**Definition 5.1.1**. Let $g(x) : \mathbb{R} \to \mathbb{R}$ be any real analytic function such that: (i)

$$g_{\mathbb{Q}}(x) = \sum_{n=0}^{\infty} a_n x^n, |x| < r, \tag{228}$$

where $\forall n(a_n \in \mathbb{Q})$ and where (ii) the sequence $\{a_n\}_{n \in \mathbb{N}} \in M_{st}[\mathbb{Z}_2^{\Sigma}]$ (in particular $\{a_n\}_{n \in \mathbb{N}} \in M_{st}\left[\mathbb{Z}_2^{\Re_1}\right]$ if the sequence $\{a_n\}_{n \in \mathbb{N}}$ is primitive recursive .

Then we will call any function given by Eq. (228) $\mathbb{Q}$-analytic $\Sigma$-function and denoted such functions by $g_{\mathbb{Q}}^{\Sigma}(x)$. In particular we will call any function $g_{\mathbb{Q}}^{\Re_1}(x)$ constructive $\mathbb{Q}$-analytic function.

**Definition 5.1.2**. A transcendental number $z \in \mathbb{R}$ is called $\Sigma$-transcendental number over field

$\mathbb{Q}$, if there does not exist $\mathbb{Q}$-analytic $\Sigma$-function $g_{\mathbb{Q}}^{\Sigma}(x)$ such that $g_{\mathbb{Q}}^{\Sigma}(z) = 0$.

In particular a transcendental number $z \in \mathbb{R}$ is called #-transcendental number over field $\mathbb{Q}$, if there does not exist constructive $\mathbb{Q}$-analytic function $g_{\mathbb{Q}}^{\Re_1}(x)$ such that $g_{\mathbb{Q}}(z) = 0$, i.e. for every constructive $\mathbb{Q}$-analytic function $g_{\mathbb{Q}}^{\Re_1}(x)$ the inequality $g_{\mathbb{Q}}^{\Re_1}(z) \neq 0$ is satisfied.

**Example 5.1.1**. Number $\pi$ is transcendental but number $\pi$ is not #-transcendental number over field $\mathbb{Q}$ since

(1) function $\sin x$ is a $\mathbb{Q}$-analytic and

(2) $\sin\left(\frac{\pi}{2}\right) = 1$, i.e.

$$-1 + \frac{\pi}{2} - \frac{\pi^3}{2^3 3!} + \frac{\pi^5}{2^5 5!} - \frac{\pi^7}{2^7 7!} + \ldots + \frac{(-1)^{2n+1}\pi^{2n+1}}{2^{2n+1}(2n+1)!} + \ldots = 0. \tag{229}$$

**Remark 5.1.2**. Note that a sequence $a_n = \frac{(-1)^{2n+1}}{2^{2n+1}(2n+1)!}, n = 0,1,2.\ldots.$ obviously is primitive recursive and therefore

$$\{a_n\}_{n\in\mathbb{N}} \in M_{st}\left[\mathbb{Z}_2^{\Re_1}\right], \tag{230}$$

since we assume $Con\left(\mathbb{Z}_2^{\Re_1} + \exists M_{st}\left[\mathbb{Z}_2^{\Re_1}\right]\right)$.

**Propostion 5.1.1**. Let $v_0 = 1$. For each $n > 0$ choose a rational number $v_n$ inductively such that

$$1 - \sum_{k=1}^{n-1} v_k e^k - (n!)^{-1} < v_n e^n < 1 - \sum_{k=1}^{n-1} v_k e^k. \tag{231}$$

The rational number $v_n$ exists because the rational numbers are dense in $\mathbb{R}$. Now the power series $f(x) = 1 - \sum_{n=1}^{\infty} v_n e^n$ has the radius of convergence $\infty$ and $f(e) = 0$. However any sequence $\{v_n\}_{n\in\mathbb{N}}$ obviously is not primitive recursive and therefore

$$\{v_n\}_{n\in\mathbb{N}} \notin M_{st}\left[\mathbb{Z}_2^{\Re_1}\right]. \tag{232}$$

**Theorem 5.1.1**. [22]cite: Foukzon15a Assume that $Con\left(\mathbb{Z}_2^{\Re_1} + \exists M_{st}\left[\mathbb{Z}_2^{\Re_1}\right]\right)$. Then number $e$ is #-transcendental over the field $\mathbb{Q}$.

**Theorem 5.1.2**. [22]cite: Foukzon15a Number $e^e$ is transcendental over the field $\mathbb{Q}$.

Proof. Immediately from Theorem 5.1.2.

**Theorem 5.1.3**. [22]cite: Foukzon15a Assume that $Con(\mathbb{Z}_2^{\Sigma} + \exists M_{st}[\mathbb{Z}_2^{\Sigma}])$. Then number $e$ is $\Sigma$-transcendental over the field $\mathbb{Q}$.

## 5mm . **5.2. The set theory.**$ZFC^{\#}$ **with a nonstandard axiom of infinity** 2mm

We remind that a major part of modern set theory involves the study of different models of $ZF$ and $ZFC$. It is crucial for the study of such models to know which properties of a set are absolute to different models [8]cite: Cohen66. It is common to begin with a fixed model of set theory and only consider other transitive models containing the same ordinals as the fixed model.

Certain fundamental properties are absolute to all transitive models of set theory, including the following: (i) $x$ is the empty set, (ii) $x$ is an ordinal, (iii) $x$ is a finite ordinal, (iv) $x = \omega$, (v) $x$ is (the graph of) a function. Other properties, such as countability, $x = 2^y$ are not absolute, see [8]cite: Cohen66.

**Remark 5.2.1**. Note that for nontransitive models the properties (ii)-(v) no longer holds.

Let $\langle M, \widetilde{\in}\rangle$ be a non standard model of $ZFC$. It follows from consideration above that any such model $\langle M, \widetilde{\in}\rangle$ is substantially non standard model of $ZFC$, i.e., there does not exist an standard model $\langle M_{st}, \in\rangle$ of $ZFC$ such that $M_{st} \subset M$, where

$$\widetilde{\in} M_{st} = \in M_{st}. \tag{233}$$

and

$$\omega \notin \langle M, \widetilde{\in}\rangle. \tag{234}$$

**Theorem 5.2.1**. [9]cite: Godel68. Let $\langle M, \widetilde{\in}\rangle$ be a non standard model of $ZF$. A necessary and sufficient condition for $\langle M, \widetilde{\in}\rangle$ to be isomorfic to a standard model $\langle M, \in\rangle$ is that there does not exist a countable sequence $\{x_n\}_{n\in\mathbb{N}}$ of elements in $M$ such that $x_{n+1}\widetilde{\in}x_n$.

**Definition 5.2.1**. Let $M_{Nst} = \langle M, \widetilde{\in}\rangle$ be a non standard model of $ZFC$. We will say that:

(i) element $z\widetilde{\in}M$ is a non standard relative to $\mathbb{N}$ and abbreviate $Nst_{\mathbb{N}}(z)$, if there exists a sequence $\{x_n\}_{n\in\mathbb{N}}$ of elements in $M$ such that $x_{n+1}\widetilde{\in}x_n$ and $z = x_0$ , and

(ii) element $z\widetilde{\in}M_{Nst}$ is a standard relative to $\mathbb{N}$ and abbreviate $st(z)$ if there does not exist a countable sequence $\{x_n\}_{n\in\mathbb{N}}$ of elements in $M$ such that $x_{n+1}\widetilde{\in}x_n$ and $z = x_0$, i.e.,

$st(z) \Leftrightarrow \neg Nst_{\mathbb{N}}(z)$.

**Remark 5.2.2**. We denote by $ZFC_{\widetilde{\in}}$ set theory which is obtained from set theory $ZFC$ by using wff's of $ZFC$ with quantifiers bounded on a non standard model $\langle M, \widetilde{\in}\rangle$. The first-order language corresponding to set theory $ZFC_{\widetilde{\in}}$ we denote by $\mathcal{L}_{\widetilde{\in}}$.

Let $\mathbf{Wff}(ZFC_{\widetilde{\in}})$ be a set of the all wff's of $ZFC_{\widetilde{\in}}$. Note that $st_{\mathbb{N}}(z), Nst_{\mathbb{N}}(z) \notin \mathbf{Wff}(ZFC_{\widetilde{\in}})$, i.e., predicates $st_{\mathbb{N}}(z)$ and $Nst_{\mathbb{N}}(z)$ are not well defined in $ZFC_{\widetilde{\in}}$ since $\mathbb{N} \notin M_{Nst}$.

**Definition 5.2.2**. In set theory, an ordinal number $\alpha$ is an admissible ordinal if $L_\alpha$ is an admissible set (that is, a transitive model of Kripke–Platek set theory); in other words, $\alpha$ is admissible when $\alpha$ is a limit ordinal and $L_\alpha \models \Sigma_0$-collection.

**Definition 5.2.3**. Let $\langle M, \widetilde{\in}\rangle$ be a non standard model of $ZF$. Assume that ordinal of $\langle M, \widetilde{\in}\rangle$ have a largest minimal segment isomorfic to some standard ordinal $\alpha \in M$, which is called the <u>standard</u> <u>part</u> of $\langle M, \widetilde{\in}\rangle$, see [14]-[15]cite: Friedman73cite: MagidorShelahStavi83. We shall

assume that $\alpha \subset M$, and that for $\beta < \alpha$ :

$$\widetilde{\in}\mathfrak{R}^M(\beta) = \in \mathfrak{R}^M(\beta), \tag{235}$$

where $\mathfrak{R}^M(\beta)$ is the set of all elements of $M$ with $M$ rank is less then $\beta$.

Which standard ordinal $\alpha$ can be standard part of $\langle M, \widetilde{\in}\rangle$? It well-known that a necessary condition is that $\alpha$ is admissible ordinal. A well-known Friedman theorem (see [14]-[15]cite: Friedman73cite: MagidorShelahStavi83) implies that for countable $\alpha$ the admissibility is also sufficient condition. Thus there is no admissible countable ordinal $\alpha$ in any non standard model of $ZFC$.

**Remark 5.2.3**. We introduce now in consideration an <u>conservative</u> extension of the theory $ZFC_{\widetilde{\in}}$ by adding to language $\mathcal{L}_{\widetilde{\in}}$ the atomic predicate $Nst(z)$ which satisfies the following condition

$$\forall z[Nst(z) \Rightarrow \exists x[(x\widetilde{\in}z) \wedge Nst(x)]]. \tag{236}$$

**1. Axioms of non standartness**

(a) There exists at least one non standard set

$$\exists z[Nst(z)]. \tag{237}$$

(b) There exists at least one non standard transitive set

$$\exists z[Nst(z) \wedge TR(z)], \tag{238}$$

where: $TR(z) \Leftrightarrow \forall x[(x\widetilde{\in}z \wedge \alpha\widetilde{\in}x) \Rightarrow \alpha\widetilde{\in}z]$.

**2. Axiom of extensionality**

$$\forall x\forall y[\forall z(z\widetilde{\in}x \Leftrightarrow z\widetilde{\in}y) \Rightarrow x = y]. \tag{239}$$

**3. Axiom of regularity**

$$\forall x\exists\alpha[(\exists\alpha(\alpha\widetilde{\in}x)) \Rightarrow \exists y(y\widetilde{\in}x) \wedge \neg\exists z[(z\widetilde{\in}y) \wedge (z\widetilde{\in}x)]]. \tag{240}$$

**4. Axiom schema of specification**

Let $\phi^{st}$ be any formula in the language of $ZFC_{\widetilde{\in}}$ such that (i) formula $\phi^{st}$ free from occurrence of the atomic predicate $Nst(z)$, i.e., $\phi^{st}$ <u>can not contain the atomic predicate $Nst(z)$</u> and

(ii) $\phi^{st}$ is a formula with all free variables among $x, z, w_1, \ldots, w_n$ ($y$ is not free in $\phi^{st}$ ). Then:

$$\forall z\forall w_1\ldots\forall w_n\exists y\forall x[x\widetilde{\in}y \iff (x\widetilde{\in}z)\wedge\phi^{st}(x,z,w_1,\ldots,w_n)]. \tag{241}$$

**4′. Axiom of empty set**

$$\exists x\forall y\left[\neg\left(y\widetilde{\in}\, x\right)\right]. \tag{242}$$

We will denote the empty set by $\widetilde{\varnothing}$.

**5. Axiom of pairing**

$$\forall x\forall y\exists z[x\widetilde{\in}z\wedge y\widetilde{\in}z]. \tag{243}$$

**6. Axiom of union**

$$\forall\mathcal{F}\exists A\forall Y\forall x\left[\left(\left(x\,\widetilde{\in}Y\right)\wedge\left(Y\,\widetilde{\in}\mathcal{F}\right)\right)\Rightarrow x\,\widetilde{\in}A\right]. \tag{244}$$

**7. Axiom schema of replacement**

The axiom schema of replacement asserts that the image of a set under any definable in $ZFC_{\widetilde{\in}}$ function will also fall inside a set.

Let $\phi^{st}$ be any formula in the language of $ZFC_{\widetilde{\in}}$ such that (i) formula $\phi^{st}$ free from occurrence of the atomic predicate $Nst(z)$, i.e., $\phi^{st}$ can not contain the atomic predicate $Nst(z)$ and

(ii) $\phi^{st}$ is a formula whose free variables are among $x,y,A,w_1,\ldots,w_n$, so that in particular $B$ is not free in $\phi^{st}$. Then:

$$\begin{gathered}\forall A\forall w_1\ldots\forall w_n[\forall x\left(x\,\widetilde{\in}A\Rightarrow\exists!y\phi^{st}(A,w_1,\ldots,w_n,x,y)\right)\Rightarrow\\ \exists B\forall x\left(x\,\widetilde{\in}A\Rightarrow\exists y(y\widetilde{\in}B\wedge\phi^{st}(A,w_1,\ldots,w_n,x,y))\right)].\end{gathered} \tag{245}$$

**8. Axiom of infinity**

Let $S_{\widetilde{\in}}(x)$ abbreviate $x\cup_{\widetilde{\in}}\{x\}_{\widetilde{\in}}$, where $w$ is some set. Then:

$$\exists\mathbf{I}\left[\left(\widetilde{\varnothing}\,\widetilde{\in}\mathbf{I}\right)\wedge\forall x\widetilde{\in}\mathbf{I}\left(x\cup_{\widetilde{\in}}\{x\}_{\widetilde{\in}}\,\widetilde{\in}\mathbf{I}\right)\right]. \tag{246}$$

Such a set as usually called an inductive set.

**Definition 5.2.4**. We will say that $x$ is a non standard set and abbreviate $x^{Nst}$ iff $x$ contain at least one non standard element, i.e.,

$$x^{Nst}\iff\exists\alpha[\alpha\widetilde{\in}x\wedge Nst(\alpha)]. \tag{247}$$

**Remark 5.2.4**. It follows from Axiom schema of specification and Axiom schema of replacement (245) we can not extract from a non standard set the standard and non standard elements separately, i.e. for any non standard set $x^{Nst}$ there is no exist a set $y$ and $z$ such that

$$x^{Nst}=y\cup_{\widetilde{\in}}z, \tag{248}$$

where $y$ contain only standard sets and $z$ contain only standard sets!

As it follows from Theorem 5.3.1 any inductive set is a non standard set.

Thus Axiom of infinity can be written in the following form

**8′. Axiom of infinity**

Let $S_{\widetilde{\in}}(x)$ abbreviate $x\cup_{\widetilde{\in}}\{x\}_{\widetilde{\in}}$, where $w$ is some set. Then:

$$\exists\mathbf{I}^{Nst}\left[\left(\widetilde{\varnothing}\,\widetilde{\in}\mathbf{I}^{Nst}\right)\wedge\forall x\widetilde{\in}\mathbf{I}^{Nst}\left(x\cup_{\widetilde{\in}}\{x\}_{\widetilde{\in}}\,\widetilde{\in}\mathbf{I}^{Nst}\right)\right]. \tag{249}$$

Such a set as usually called a non standard inductive set.

**9. Strong axiom of infinity**

Let $S_{\widetilde{\in}}(x)$ abbreviate $x\cup_{\widetilde{\in}}\{x\}_{\widetilde{\in}}$, where $w$ is some set. Then:

$$\exists \mathbf{I}^{Nst}\left\{[TR(\mathbf{I}^{Nst})] \wedge \left[\left(\widetilde{\varnothing} \widetilde{\in} \mathbf{I}^{Nst}\right) \wedge \forall x \widetilde{\in} \mathbf{I}^{Nst}\left(x \cup_{\widetilde{\in}} \{x\}_{\widetilde{\in}} \widetilde{\in} \mathbf{I}^{Nst}\right)\right]\right\}. \quad (250)$$

## 5mm . **5.3. Extracting the standard and nonstandard natural numbers from the infinite nonstandard set** $\mathbf{I}^{Nst}$ 2mm

**Definition 5.3.1**. We will say that $x^{Nst}$ is inductive if there is an formula $\Phi(x)$ of $ZFC_{\widetilde{\in}}$ that says: '$x^{Nst}$ is $\widetilde{\in}$-inductive'; i.e.

$$\Phi(x^{Nst}) = \left(\widetilde{\varnothing} \widetilde{\in} x^{Nst} \wedge \forall y\left(y \widetilde{\in} x^{Nst} \Rightarrow S_{\widetilde{\in}}(y) \widetilde{\in} x^{Nst}\right)\right). \quad (251)$$

Thus we wish to prove the existence of a unique non standard set $\widetilde{W}^{Nst}$ such that

$$\forall x\left[x \widetilde{\in} \widetilde{W}^{Nst} \Leftrightarrow \forall I^{Nst}(\Phi(I^{Nst}) \Rightarrow x \widetilde{\in} I^{Nst})\right]. \quad (252)$$

(1) For existence, we will use the Axiom of Infinity combined with the Axiom schema of specification. Let $I^{Nst}$ be an inductive (non standard) set guaranteed by the Axiom of Infinity. Then we use the Axiom Schema of Specification to define our set

$$\widetilde{W}^{Nst} = \{x \in I^{Nst} : \forall J^{Nst}(\Phi(J^{Nst}) \to x \widetilde{\in} J^{Nst})\}, \quad (253)$$

i.e. $\widetilde{W}^{Nst}$ is the set of all elements of $I^{Nst}$ which happen also to be elements of every other inductive set. This clearly satisfies the hypothesis of (5.3.2), since if $x \widetilde{\in} \widetilde{W}^{Nst}$, then $x$ is in every inductive set, and if $x$ is in every inductive set, it is in particular in $I^{Nst}$, so it must also be in $\widetilde{W}^{Nst}$.

(2) For uniqueness, first note that any set which satisfies (252) is itself inductive, since $\widetilde{\varnothing}$ is in all inductive sets, and if an element $x$ is in all inductive sets, then by the inductive property so is its successor. Thus if there were another set $\widetilde{W}_1^{Nst}$ which satisfied (252) we would have that $\widetilde{W}_1^{Nst} \subseteq_{\widetilde{\in}} \widetilde{W}^{Nst}$ since W is inductive, and $\widetilde{W}^{Nst} \subseteq_{\widetilde{\in}} \widetilde{W}_1^{Nst}$ since $\widetilde{W}_1^{Nst}$ is inductive.

Thus $\widetilde{W}_1^{Nst} =_{\widetilde{\in}} \widetilde{W}^{Nst}$. Let $\widetilde{\omega}$ denote this unique set.

(3) For nonstardarntness we assume that $\widetilde{\omega}$ is a standard set, i.e. there is no nonstandard element in $\widetilde{\omega}$. Then $\widetilde{\omega} =_{\widetilde{\in}} \mathbb{N}_{\widetilde{\in}}$ where $\mathbb{N}_{\widetilde{\in}}$ is isomorphic to $\mathbb{N}$, but this is a contradiction, since $\mathbb{N} \notin_{\widetilde{\in}} \langle M, \widetilde{\in} \rangle$.

**Theorem 5.3.1**. There exists unique nonstardard set $\widetilde{\omega}$ such that (252) holds, i.e.

$$\forall x\left[x \widetilde{\in} \widetilde{\omega} \Leftrightarrow \forall I^{Nst}(\Phi(I^{Nst}) \Rightarrow x \widetilde{\in} I^{Nst})\right]. \quad (254)$$

**Definition 5.3.2**. We will say that a set S is $\widetilde{\in}$-finite if every surjective $\widetilde{\in}$-function from $S$ onto itself is one-to-one.

**Theorem 5.3.2**. There exist $\widetilde{\in}$-finite nonstandard natural numbers in $\widetilde{\omega}$.

Proof. Assuming that any nonstandard natural number is not $\widetilde{\in}$-finite one obviously obtains a contradiction.

**Remark 5.3.1.**Assuming that $\widetilde{\omega}$ is ф standard set then this method mentioned above produce system which satisfy the axioms of second-order arithmetic $Z_2^{fss}$, since the axiom of power set allows us to quantify over the power set of $\widetilde{\omega}$ , as in second-order logic. Thus it completely determine isomorphic systems, and since they are isomorphic under the identity map, they must in fact be equal.

# 6mm . **6. Conclusion** 3mm

In this paper we have proved that the second order *ZFC* with the full second-order semantic is inconsistent, i.e. $\neg Con(ZFC_2^{fss})$. Main result is: let $k$ be an inaccessible cardinal and $H_k$ is a set of all sets having hereditary size less then $k$, then $\neg Con(ZFC + \exists M_{st}^{ZFC}(M_{st}^{ZFC} = H_k))$. This result also was obtained in [3], [4], [5]cite: Lemhoff16cite: Foulzon17cite: Foukzon15 essentially another approach. For the first time this result has been declared to AMS in [23], [24]cite: Foulzon13a cite: Foukzon13b

## Acknowledgments

Reviewers provided important clarifications.
3mm
We thank the reviewers for their comments.